\documentclass{amsart}
\usepackage{amsmath,amssymb,amscd,latexsym}
\usepackage[all]{xy}

\newcommand{\N}{{\mathbb{N}}}

\newcommand{\R}{{\mathbb{R}}}

\newcommand{\Bh}{{\mathcal B}}
\newcommand{\Ch}{{\mathcal C}}

\newcommand{\Eh}{{\mathcal E}}
\newcommand{\Fh}{{\mathcal F}}
\newcommand{\Gh}{{\mathcal G}}

\newcommand{\Zh}{{\mathcal Z}}

\newcommand{\be}{\mathbf{1}}

\newcommand{\dimnuc}{\mathrm{dim}_{\mathrm{nuc}}}

\newcommand{\dr}{\mathrm{dr}\,}

\newcommand{\halb}{\frac{1}{2}}
\newcommand{\her}{\mathrm{her}}

\newcommand{\id}{\mathrm{id}}

\newcounter{number}[section]

\newenvironment{nummer}{\refstepcounter{number}{\noindent\arabic{section}.\arabic{number}}}{}

\newcommand{\bn}{\noindent \begin{nummer} \rm}
\newcommand{\en}{\end{nummer}}

\newenvironment{ntheorem}{\noindent {\sc Theorem:} \it}{}
\newenvironment{nlemma}{\noindent {\sc Lemma:} \it}{}
\newenvironment{nprop}{\noindent {\sc Proposition:} \it}{}
\newenvironment{ndefn}{\noindent {\sc Definition:} \it}{}
\newenvironment{ncor}{\noindent {\sc Corollary:} \it}{}

\newenvironment{nremark}{\noindent {\sc Remark:}}{}
\newenvironment{nremarks}{\noindent {\sc Remarks:}}{}

\newenvironment{nnotation}{\noindent {\sc Notation:} }{}

\newenvironment{nproof}{\noindent {\sc Proof:}}{\mbox{}\hfill 
\rule[-.2ex]{.25em}{1.8ex}}

\begin{document}

\title[Nuclear dimension and $\mathcal{Z}$-stability of pure $\mathrm{C}^{*}$-algebras]{{\sc Nuclear dimension and $\mathcal{Z}$-stability\\ of pure $\mathrm{C}^{*}$-algebras}}

\dedicatory{In Erinnerung an meinen Vater}

\author{Wilhelm Winter}
\address{School of Mathematical Sciences\\
University of Nottingham\\
United Kingdom}

\email{wilhelm.winter@nottingham.ac.uk}

\date{\today}
\thanks{Partially supported by EPSRC First Grant EP/G014019/1}

\setcounter{section}{-1}

\begin{abstract}
In this article I  study a number of topological and algebraic dimension type properties of simple $\mathrm{C}^{*}$-algebras  and their interplay. In particular,  a simple $\mathrm{C}^{*}$-algebra is defined to be (tracially) $(m,\bar{m})$-pure, if it has (strong tracial) $m$-comparison and is (tracially) $\bar{m}$-almost divisible. These notions are related to each other, and to nuclear dimension.

The main result says that if a separable, simple, nonelementary, unital $\mathrm{C}^{*}$-algebra with locally finite nuclear dimension is $(m,\bar{m})$-pure, then it absorbs the Jiang--Su algebra $\Zh$ tensorially.  It follows that a separable, simple, nonelementary, unital $\mathrm{C}^{*}$-algebra with locally finite nuclear dimension is $\mathcal{Z}$-stable if and only if it has the Cuntz semigroup of a $\mathcal{Z}$-stable $\mathrm{C}^{*}$-algebra. The result may be regarded as a  version of Kirchberg's celebrated theorem that separable, simple, nuclear, purely infinite $\mathrm{C}^{*}$-algebras absorb the Cuntz algebra $\mathcal{O}_{\infty}$ tensorially.  

As a corollary we obtain that finite nuclear dimension implies $\mathcal{Z}$-stability for separable, simple, nonelementary, unital $\mathrm{C}^{*}$-algebras; this settles an important case of a conjecture by Toms and the author. 

The main result also has a number of consequences for Elliott's program to classify nuclear $\mathrm{C}^{*}$-algebras by their $\mathrm{K}$-theory data. In particular, it completes the classification of simple, unital, approximately homogeneous algebras with slow dimension growth by their Elliott invariants, a question left open in the Elliott--Gong--Li classification of simple AH algebras. 

Another consequence is that for simple, unital, approximately subhomogeneous algebras, slow dimension growth and $\mathcal{Z}$-stability are equivalent. In the case where projections separate traces,  this completes the classification of simple, unital, approximately subhomogeneous algebras with slow dimension growth by their ordered $\mathrm{K}$-groups. 
\end{abstract}

\maketitle

\section{Introduction}

\noindent
In the structure and classification theory of nuclear $\mathrm{C}^{*}$-algebras, the interplay between algebraic and topological regularity properties is becoming increasingly important. Algebraic regularity occurs in many different ways; it might be related to the number of generators of certain modules, to tensorial absorption of suitable touchstone algebras, or  to structural properties of homological invariants. Topological regularity properties usually stem from the interpretation of $\mathrm{C}^{*}$-algebras as noncommutative topological spaces. They might   manifest themselves as properties of the primitive ideal space (e.g.\ separation  properties),  but there is also a rich supply of noncommutative versions of topological dimension, which carry nontrivial information even when the spectrum doesn't. Examples of such concepts are stable and real rank (\cite{Rfl:sr} and \cite{BroPed:realrank}), or (in increasing order of generality) approximately homogeneous dimension, approximately subhomogeneous dimension, decomposition rank, and nuclear dimension (cf.\ \cite{Ror:encyc}, \cite{KirWinter:dr} and \cite{WinterZac:dimnuc}). Of these, the former are linked by design to certain algebraic regularity properties; the latter concepts behave much more like topological invariants, and their connections to algebraic regularity are subtle.   

In \cite{TomsWinter:VI}, \cite{Winter:dr-Z-stable} and \cite{WinterZac:dimnuc}, some links between such regularity properties have been proven, and some others have been formally established as conjectures (see also \cite{EllToms:BullAMS}). In particular, A.~S.~Toms and the author have conjectured that finite decomposition rank, $\mathcal{Z}$-stability and strict comparison are equivalent, at least for   separable, simple, nonelementary, unital, nuclear and stably finite $\mathrm{C}^{*}$-algebras. The main result of \cite{Winter:dr-Z-stable} establishes one implication of this conjecture: a  separable, simple, nonelementary, unital $\mathrm{C}^{*}$-algebra with finite decomposition rank is $\mathcal{Z}$-stable, i.e., absorbs the Jiang--Su algebra $\mathcal{Z}$ tensorially. By \cite{Ror:Z-absorbing}, this also implies strict comparison of positive elements (i.e., the order structure of the Cuntz semigroup is determined by lower semicontinuous dimension functions, which in turn are given by tracial states). The result has been particularly useful for Elliott's  classification program for nuclear $\mathrm{C}^{*}$-algebras; for example it has led to the complete classification of $\mathrm{C}^{*}$-algebras associated to uniquely ergodic, minimal, finite-dimensional dynamical systems, in  \cite{TomsWinter:minhom} and \cite{TomsWinter:PNAS}.

The conjecture of Toms and the author was formulated in the not necessarily finite case in \cite{WinterZac:dimnuc}, replacing decomposition rank by nuclear dimension. This generalization is interesting since on the one hand it allows a unified view on the purely infinite and the purely finite branches of the classification program; on the other hand, nuclear dimension is highly accessible to applications, in particular in dynamical systems, but also in coarse geometry.

In the present paper  some more regularity properties for simple $\mathrm{C}^{*}$-algebras are studied, $m$-comparison and $m$-almost divisibility. These have an algebraic flavour and can be formulated at the level of the Cuntz semigroup, but they might as well be interpreted as dimension type conditions. I also introduce and study tracial versions of these conditions, (strong) tracial $m$-comparison and tracial $m$-almost divisibility. I  then call a simple $\mathrm{C}^{*}$-algebra (tracially) $(m,\bar{m})$-pure, if it has (strong tracial) $m$-comparison and (tracial) $\bar{m}$-almost divisibility; it is shown that tracial $(m,\bar{m})$-pureness implies $(\tilde{m},\bar{m})$-pureness (for a suitable $\tilde{m}$ depending on $m$ and $\bar{m}$). 

L.~Robert has shown in \cite{Rob:dimnuc-comparison} that if a $\mathrm{C}^{*}$-algebra has nuclear dimension at most $m$, then it has $m$-comparison in the sense of \cite{Rob:dimnuc-comparison} (in the simple case, Robert's notion of $m$-comparison agrees with the definition used here). I will use a result of \cite{WinterZac:dimnuc} on Kirchberg's covering number to show that finite nuclear dimension (at least in the simple case)  implies $\bar{m}$-almost divisibility for some $\bar{m} \in \mathbb{N}$;  tracial $(m,\bar{m})$-pureness will also be derived in this situation. This generalizes results of \cite[Section~3]{Winter:dr-Z-stable} to the case of nuclear dimension; \cite[Lemma~3.10 and Corollary~3.12]{Winter:dr-Z-stable}  essentially yield tracial $(m,\bar{m})$-pureness for simple, unital $\mathrm{C}^{*}$-algebras with finite decomposition rank, and provide crucial ingredients for the main result of \cite{Winter:dr-Z-stable}.

While the subtle interplay between nuclear dimension and (tracial) $(m,\bar{m})$-pure\-ness is interesting by itself, these relations become particularly useful when combined with the main technical result, Theorem~\ref{A}.  This says that a simple, unital, nonelementary $\mathrm{C}^{*}$-algebra with locally finite nuclear dimension is $\mathcal{Z}$-stable, if it is tracially $(m,\bar{m})$-pure. The proof of Theorem~\ref{A} is quite technical; it will to some extent follow ideas from \cite{Winter:dr-Z-stable}, isolating and generalizing some of the abstract machinery of that paper. 

As a consequence of these results we obtain that a separable, simple, unital, nonelementary $\mathrm{C}^{*}$-algebra with finite nuclear dimension is $\mathcal{Z}$-stable, Corollary~\ref{findimnuc-Z-stable}, thus establishing one of the implications of \cite[Conjecture~9.3]{WinterZac:dimnuc}. As a second consequence, the statement that a pure, separable, simple and unital $\mathrm{C}^{*}$-algebra with locally finite nuclear dimension is $\mathcal{Z}$-stable  (Corollary~\ref{Acor}) partially verifies another implication of \cite[Conjecture~9.3]{WinterZac:dimnuc}. (The reverse implication was established by M.~R{\o}rdam, who has shown that pureness is a necessary condition for $\mathcal{Z}$-stability.) Remarkably, then, if $A$ is simple and unital with locally finite nuclear dimension,  $\mathcal{Z}$-stability can be read off from the Cuntz semigroup, and, in fact, $A$ is $\mathcal{Z}$-stable if and only if the Cuntz semigroups of $A$ and $A \otimes \mathcal{Z}$ agree, see Corollary~\ref{W-Z-stable}.

Corollaries~\ref{Acor} and \ref{findimnuc-Z-stable} are both nontrivial, and both are useful in the structure theory and in the classification program for nuclear $\mathrm{C}^{*}$-algebras. For example, the result on finite nuclear dimension provides an alternative (to some extent more direct, but not necessarily simpler) access to  \cite[Theorem~0.2]{TomsWinter:minhom} (still using \cite[Theorem~0.3]{TomsWinter:minhom}). It will also be a crucial ingredient of one of the main results of \cite{HirWinZac:Rokhlin-dimension}, which says that if $A$ is a separable, simple, unital $\mathrm{C}^{*}$-algebra with finite nuclear dimension, then there is a dense $G_{\delta}$ subset of the space of automorphisms of $A$ such that  the respective crossed product $\mathrm{C}^{*}$-algebras are again simple with finite nuclear dimension. 

The result on $\mathcal{Z}$-stability of pure $\mathrm{C}^{*}$-algebras in combination with a theorem of A.~S.~Toms (\cite{Toms:ash-sdg})  can be used to make progress on the classification and structure theory for  simple, unital, approximately (sub)homogeneous $\mathrm{C}^{*}$-algebras. Corollary~\ref{ASH-sdg-Z} says that for such algebras slow dimension growth is equivalent to $\mathcal{Z}$-stability.  In the approximately subhomogeneous case with slow dimension growth, this confirms the Elliott conjecture provided projections separate traces -- Corollary~\ref{ASH-sdg-classification}. In the approximately homogeneous case with slow dimension growth,  we obtain classification without any additional trace space condition -- see Corollary~\ref{AH-sdg-classification}; this answers an important technical question left open in \cite{EllGongLi:simple_AH}, where classification was proven under the more restrictive (and less natural) condition of \emph{very} slow dimension growth.  (It should be noted, however, that the result of \cite{EllGongLi:simple_AH} is still a necessary step.)

In the absence of traces, a pure, separable, simple, unital $\mathrm{C}^{*}$-algebra is purely infinite in the sense of \cite{Cuntz:On}. Therefore, the statement that pureness implies $\mathcal{Z}$-stability (at least under the additional  structural hypothesis of locally finite nuclear dimension) may be regarded as a partial generalization of E.~Kirchberg's celebrated result that purely infinite, separable, simple, nuclear $\mathrm{C}^{*}$-algebras absorb the Cuntz algebra $\mathcal{O}_{\infty}$; see \cite{Kir:ICM}. (Kirchberg's theorem later became one of the cornerstones of the Kirchberg--Phillips classification of purely infinite, simple, nuclear $\mathrm{C}^{*}$-algebras.)

It should be mentioned that it remains an open question how restrictive the hypothesis of locally finite nuclear dimension really is. In fact, at this point we do not know of \emph{any} nuclear $\mathrm{C}^{*}$-algebra which does not have locally finite nuclear dimension.

I also wish to point out that the methods of this article do not in any way involve a dichotomy between the purely finite (i.e., finite and $\mathcal{Z}$-stable) and the purely infinite (i.e., $\mathcal{O}_{\infty}$-stable or, equivalently, infinite and $\mathcal{Z}$-stable) situation. This may be interpreted as promising evidence towards an approach to classification which would unify the purely finite and the purely infinite case.

The paper is organized as follows. Section~{\ref{preliminaries}} introduces some notation, recalls the notions of order zero maps, decomposition rank and nuclear dimension, as well as a useful criterion for $\mathcal{Z}$-stability. In Section~{\ref{pure-finiteness}},  the interplay of several versions of $m$-comparison and $m$-almost divisibility for $\mathrm{C}^{*}$-algebras is studied and  the notion of  (tracially) $(m,\bar{m})$-pure $\mathrm{C}^{*}$-algebras is introduced. In Section~{\ref{dimnuc-pureness}} it is shown that finite nuclear dimension implies tracial $(m,\bar{m})$-pureness, at least in the simple and unital case. Section~{\ref{tracial-matrix-cone-absorption}} generalizes \cite[Section~3]{Winter:dr-Z-stable} (it provides tracially large almost central order zero maps into tracially $m$-almost divisible $\mathrm{C}^{*}$-algebras with locally finite nuclear dimension); in a somewhat similar manner, Section~{\ref{almost-central-dimension-drop-embeddings}} generalizes \cite[Section~4]{Winter:dr-Z-stable}. The main result and its corollaries are established in Section~{\ref{main-result}}. 

I am indebted to G.\ Elliott and to the referee for a number of comments on  earlier versions of the paper.


\section{Preliminaries}
\label{preliminaries}

\noindent
Below we recall the definition of and some facts about order zero maps, the Cuntz semigroup and $\mathcal{Z}$-stability; we also fix some notation which will be used frequently throughout the paper.

\bn
Recall from \cite{WinZac:order-zero} that a completely positive (c.p.)  map $\varphi:A \to B$ between    $\mathrm{C}^{*}$-algebras   $A$ and $B$ is said to have order zero if it preserves orthogonality. 
We collect below some  facts about  
order zero maps from \cite{WinZac:order-zero} (see \cite[Proposition~3.2(a)]{Win:cpr}  and 
\cite[1.2]{Winter:fintopdim}  for the case of finite-dimensional domains). 

\begin{nprop}
\label{order-zero-facts}
Let $A$ and $B$ be  $\mathrm{C}^{*}$-algebras, and let 
$\varphi \colon A \to B$ be a completely positive contractive (c.p.c.) order zero map.
\begin{enumerate}
\item[(i)] There is a unique $*$-homomorphism 
\[
\tilde{\varphi} \colon 
\mathcal{C}_{0}((0,1]) \otimes A \to B
\]
such that 
\[
\varphi(x) = \tilde{\varphi}( \iota\otimes x)
\]
for all $x \in A$, where $\iota(t) = t$.
\item[(ii)] There are a $*$-homomorphism 
\[
\pi \colon
A \to \mathcal{M}(\mathrm{C}^{*}(\varphi(A))) \subset B''
\]
and a positive element 
\[
0 \le h  \in \mathcal{M}(\mathrm{C}^{*}(\varphi(A))) 
\]
(both uniquely determined) such that 
\begin{equation*}
\label{order-zero-facts-w1}
\varphi(x)= \pi(x) h = h \pi(x) 
\end{equation*}
for all $x \in A$.  If $A$ is unital, then $h = \varphi(\be_{A})$.
\item[(iii)] If, for some  $u  \in B$ with $\|u\|\le 1$, the element
  $u^{*}u$ commutes  
with $\pi(A)$, then the map
\[
\varphi_{u} \colon A \to B
\]
given by 
\[
 \varphi_{u}(x) = u \pi(x) u^{*}
 \] 
 for $x \in A$, is a  c.p.c.\ order zero map.
\end{enumerate}
\end{nprop}
\en

\bn
\label{order-zero-notation}
\begin{nnotation}
The map $\pi$ in \ref{order-zero-facts}(ii) above will be called the \emph{canonical supporting 
$*$-homomorphism} of $\varphi$. Following \cite{WinZac:order-zero} (also cf.\ \cite[1.3]{Winter:dr-Z-stable}), whenever $\varphi:A \to B$ is a c.p.c.\ order zero map and $f \in \Ch_{0}((0,1])$ is a positive function of norm at most one,  we define a c.p.c.\ order zero map 
\[
f(\varphi): A \to B
\]
by setting
\begin{equation*}
\label{oznw2}
f(\varphi)(x) := \pi(x) f(h).
\end{equation*}
\end{nnotation}
\en

\bn
\label{order-zero-lifts}
\begin{nremark}
Let $F$ be finite-dimensional, $\pi:A \to B$ a surjective $*$-homomorphism and $\varphi:F \to B$ a c.p.c.\ order zero map. Then, $\varphi$ lifts to a c.p.c.\ order zero map $\bar{\varphi}: F \to A$, see \cite{Win:cpr}.
\end{nremark}
\en

\bn
It will be convenient to introduce the following notation.

\begin{nnotation}
For positive numbers $0\le \eta < \epsilon \le 1$ define continuous 
functions $$f_{\eta,\epsilon}, g_{\eta,\epsilon} \colon [0,1] \to [0,1 ]$$ by
\[
g_{\eta,\epsilon}(t)= \begin{cases}0 & \mbox{ if } t \le \eta,\\ 1 & \mbox{ if } \epsilon \le t \le 1,
\\ \text{linear} &  \text{else,}\end{cases} 
\]
and
\[
f_{\eta,\epsilon}(t)= \begin{cases}0 & \mbox{ if } t \le \eta,\\ t & \mbox{ if } \epsilon \le t \le 1,
\\ \text{linear} & \text{else.}\end{cases} 
\]
\end{nnotation}
\en

\bn
\label{standard-partition}
\begin{nnotation}
For a given $L \in \mathbb{N}$, by a standard partition of unity for the unit interval $[0,1]$ we mean a set of sawtooth functions
\[
h^{(0)}_{0}, \ldots, h^{(0)}_{L}, h^{(1)}_{1}, \ldots, h^{(1)}_{L} \in \mathcal{C}([0,1]) 
\]
defined by
\[
h^{(i)}_{l}(t)= \begin{cases}0 & \mbox{ if } t \le \frac{2l-i-1}{2L},\\ 1 & \mbox{ if }  t = \frac{2l-i}{2L},
\\ 
0 & \mbox{ if } t \ge \frac{2l-i+1}{2L},\\
\text{linear} & \text{else.}\end{cases} 
\]
Note that the functions $h^{(0)}_{1}, \ldots, h^{(0)}_{L}, h^{(1)}_{1}, \ldots, h^{(1)}_{L}$ are in $\mathcal{C}_{0}((0,1])$; we will slightly abuse notation and call these a standard partition of unity for the interval $(0,1]$.
\end{nnotation}
\en

\bn
\label{traces-not}
\begin{nnotation}
If $A$ is a unital $\mathrm{C}^{*}$-algebra, we denote by $T(A)$ the space of tracial states  and by $QT(A)$ the space of normalized quasitraces, see \cite{Bla:encyc}. $T(A)$ and $QT(A)$ are compact (and metrizable if $A$ is separable). If $A$ is simple, then any $0 \neq a \in A_{+}$ is    strictly positive and continuous on  $QT(A)$, and there is $\eta>0$ such that $\tau(a) > \eta$ for all $\tau \in QT(A)$. 
\end{nnotation}
\en

\bn
\label{d-tau-notation1}
\begin{nremark}
We associate to any $\tau \in QT(A)$ a lower semicontinuous dimension function  $d_{\tau}$ on $A_{+}$ by 
\[
d_{\tau}(a):= \lim_{n \to \infty} \tau(a^{\frac{1}{n}});
\]
it was shown in \cite{BlaHan:quasitrace} that every lower semicontinuous dimension function on $A$ is induced by a quasitrace. One checks 
\[
d_{\tau}(a) = \lim_{\epsilon \searrow 0} \tau(g_{0,\epsilon}(a))  = \lim_{\epsilon \searrow 0} \tau(g_{\epsilon,2\epsilon}(a))
\]
and
\begin{equation}
\label{dtnw1}
d_{\tau}((a-\beta)_{+}) \le \tau(g_{\beta/2,\beta}(a))
\end{equation}
whenever $\beta>0$. 

If $f$ is a positive continuous function on $QT(A)$ and $a \in A_{+}$ satisfies 
\[
d_{\tau}(a) > f(\tau)
\]
for all $\tau \in QT(A)$, then there are $n \in \mathbb{N}$ and $\eta>0$ such that 
\[
\tau(a^{\frac{1}{n}})>f(\tau)
\]
and
\[
\tau(g_{\eta,2\eta}(a))>f(\tau)
\]
for all $\tau \in QT(A)$ (we still asume $A$ to be unital, so that $QT(A)$ is compact). The argument is essentially the same as for Dini's Theorem, the details are omitted.
\end{nremark}
\en

\bn
\label{order-zero-traces}
\begin{nremark}
It was shown in \cite[Corollary~3.4]{WinZac:order-zero} that the composition of a trace with a c.p.\ order zero map is again a trace; the respective statement holds for quasitraces and for (lower semicontinuous) dimension functions.
\end{nremark}
\en

\bn
\begin{nnotation}
If $A$ is a $\mathrm{C}^{*}$-algebra, we write 
\[
A_{\infty}:= \prod_{\mathbb{N}} A / \bigoplus_{\mathbb{N}}A;
\]
we regard $A$ as being embedded in $A_{\infty}$ as constant sequences. 

We denote by $T_{\infty}(A)$ those tracial states on $A_{\infty}$ which are of the form
\[
[(a_{l})_{l \in \N}] \mapsto \lim_{\omega} \tau_{l}(a_{l})
\]
for a free ultrafilter $\omega $ on $ \N$ and a sequence $(\tau_{l})_{l \in \N} \subset T(A)$. Note that any $\tau \in T_{\infty}(A)$ restricts to a tracial state on $A$.  
\end{nnotation}
\en

\bn
\label{d-tau-notation2}
\begin{nnotation}
If $A$ is a $\mathrm{C}^{*}$-algebra we denote by $W(A)$ its Cuntz semigroup; this consists of equivalence classes of  elements in 
\[
M_{\infty}(A)_{+} = \bigcup_{k=1}^{\infty} M_{k}(A)_{+}.
\]
Two elements $a,b \in M_{\infty}(A)_{+}$  are equivalent, $a \sim b$, if they are Cuntz subequivalent to each other, $a \precsim b$ and $b \precsim a$.  The elements $a$ and $b$ are Murray--von Neumann equivalent, $a \approx b$, if there is $x \in M_{\infty}(A)$ such that $x^{*}x=a$ and $xx^{*}=b$.

Any  dimension function $d_{\tau}$ on $A$ as in \ref{d-tau-notation1} induces a positive real-valued character (also denoted by $d_{\tau}$) on the Cuntz semigroup $W(A)$ with its natural order (cf.\ \cite{Cuntz:dimension}); if $A$ is unital and $\tau$ is a normalized quasitrace, then $d_{\tau}:W(A) \to \R_{+}$ is a state. 

The reader is referred to \cite{Ror:Z-absorbing} and \cite{Ror:UHFII} for more detailed information and for  background material on Cuntz comparison of positive elements and on the Cuntz semigroup. 
\end{nnotation}
\en

\bn
\label{cc-sub}
\begin{nprop}
Let $A$ be a $\mathrm{C}^{*}$-algebra, $a,d \in A$ positive contractions and $0<\beta<1$ a number.  If $d \precsim a$, then there is $x \in A$ such that
\[
x^{*}x = g_{\beta/2,\beta}(b) \mbox{ and } xx^{*} \in \her(a) \subset A.
\]
\end{nprop}

\begin{nproof}
Since $d \precsim a$, there is $y \in A$ such that 
\[
\|d-y^{*} a y\|< \frac{\beta}{4}.
\]
By \cite[Lemma~2.2]{KirRor:pi2}, there is $z \in A$ with 
\[
(d-\beta/4)_{+} = z^{*}y^{*} a y z.
\]
By functional calculus there is $t \in A$ such that 
\[
t^{*} (d-\beta/4)_{+} t = g_{\beta/2,\beta}(d).
\]
Now 
\[
x:= a^{\frac{1}{2}} y z t
\]
will have the desired properties.
\end{nproof}
\en

\bn
\label{cc-order-zero}
\begin{nprop}
Let $A$ be a $\mathrm{C}^{*}$-algebra, $a,b \in A$ positive contractions, $k \in \mathbb{N}$, $0<\beta<1$,  and $f \in \mathcal{C}_{0}((\beta,1])$ a positive function of norm at most one. If 
\[
k \cdot \langle b \rangle \le \langle a \rangle 
\]
in $W(A)$, then there is a c.p.c.\ order zero map 
\[
\varphi:M_{k} \to \her(a) \subset A
\]
such that
\[
\varphi(e_{11}) \approx f(b).
\]
\end{nprop}

\begin{nproof}
Set 
\[
d:= \be_{k} \otimes b \in M_{k} \otimes A,
\]
then 
\[
\langle d \rangle = k \cdot \langle b \rangle \le \langle a \rangle 
\]
in $W(A)$. By Proposition~\ref{cc-sub}, there is $x \in M_{k} \otimes A$ such that 
\[
x^{*}x = g_{\beta/2,\beta}(d) \mbox{ and } xx^{*} \in \her(a).
\]
Define a c.p.c.\  map
\[
\varphi:M_{k} \to  \her(a)
\]
by
\[
\varphi(y):= x(y \otimes f(b))x^{*} \mbox{ for } y \in M_{k}.
\]
Using Proposition~\ref{order-zero-facts}, it is straightforward to check that $\varphi$ has the desired properties.
\end{nproof}
\en

\bn
\label{Z-intro}
\begin{nremark}
Recall from \cite{JiaSu:Z} that the Jiang--Su algebra $\mathcal{Z}$ is the uniquely determined simple and  monotracial inductive limit of so-called prime dimension drop $\mathrm{C}^{*}$-algebras, 
\[
\mathcal{Z}= \lim_{\to} Z_{p_{\nu},q_{\nu}},
\]
where
\[
Z_{p_{\nu},q_{\nu}} = \{ f \in \Ch([0,1], M_{p_{\nu}} \otimes M_{q_{\nu}}) \mid f(0) \in M_{p_{\nu}} \otimes \be_{M_{q_{\nu}}} \mbox{ and }  f(1) \in \be_{M_{p_{\nu}}}  \otimes M_{q_{\nu}}\}
\]
for $\nu \in \mathbb{N}$, and $p_{\nu}$ and $q_{\nu}$ are relatively prime. 

It was shown in \cite{JiaSu:Z} that $\mathcal{Z}$ is strongly self-absorbing in the sense of \cite{TomsWin:ssa}, and that it is $\mathrm{KK}$-equivalent to the complex numbers. In \cite{RorWin:Z-revisited}, several alternative characterizations of the Jiang--Su algebra were given; see also \cite{DadToms:Z} and \cite{Win:ssa-Z-stable}.
\end{nremark}
\en
 
\bn
\label{Z-stable-relations} 
In order to show that a separable and unital $\mathrm{C}^{*}$-algebra $A$ is $\mathcal{Z}$-stable, i.e., $A \cong A \otimes \mathcal{Z}$, it will suffice to construct approximately central unital $*$-homomorphisms from (arbitrarily large) prime dimension drop intervals to $A$; cf.\ \cite[Proposition~2.2]{TomsWin:ZASH}. Using \cite[Proposition~5.1]{RorWin:Z-revisited}, to this end it will be enough to realize the generators and relations of \cite[Proposition~5.1(iii)]{RorWin:Z-revisited} approximately -- and in an approximately central way. This was formalized in terms of the relations $\mathcal{R}(n, \Fh,\eta)$ in \cite[Section~2]{Winter:dr-Z-stable}. For our purposes the following combination of \cite[Propositions~2.3 and 2.4]{Winter:dr-Z-stable} will be most useful.

\begin{nprop}
Let $A$ be a  separable and unital $\mathrm{C}^{*}$-algebra. Suppose that, for any $n \in \N$, any finite subset  $\Fh \subset A$, and for any  numbers  $0<\delta<1$ and $0<\zeta<1$, there are a  c.p.c.\ order zero map
\[
\varphi': M_{n} \to A
\]
and
\[
v' \in A
\]
satisfying
\begin{enumerate}
\item $\|(v')^{*} v' - (\be_{A} - \varphi'(\be_{M_{n}}))\| < \delta$,
\item $v' (v')^{*} \in \overline{(\varphi'(e_{11})- \zeta)_{+}A (\varphi'(e_{11})- \zeta)}_{+} $,
\item $\|[\varphi'(x),a]\| \le \delta \|x\|$ for all $ x \in M_{n}, \, a \in \Fh$,
\item $\|[v',a]\| < \delta$ for all $a \in \Fh$.
\end{enumerate}
Then, $A$ is $\mathcal{Z}$-stable.
\end{nprop}
\en

\section{pure $\mathrm{C}^{*}$-algebras}
\label{pure-finiteness}

\noindent
In this section we define and relate several versions of comparison and almost divisibility of a $\mathrm{C}^{*}$-algebra (or its Cuntz semigroup, respectively). The notion of a pure $\mathrm{C}^{*}$-algebra is based on these concepts.

\bn
\label{strict-comparison}
\begin{ndefn}
Let $A$ be a separable, simple, unital $\mathrm{C}^{*}$-algebra and $m \in \mathbb{N}$.
\begin{enumerate}
\item  We say $A$ has $m$-comparison of positive elements, if for any nonzero positive contractions $a,b_{0},\ldots,b_{m} \in M_{\infty}(A)$ we have 
\[
a \precsim b_{0} \oplus \ldots \oplus b_{m}
\]
whenever 
\[
d_{\tau}(a) < d_{\tau}(b_{i})
\]
for every $\tau \in QT(A)$ and $i=0,\ldots,m$. 
\item $A$ has tracial $m$-comparison  of positive elements, if for any nonzero positive contractions $a,b_{0},\ldots,b_{m} \in M_{\infty}(A)$ we have 
\[
a \precsim  b_{0} \oplus \ldots \oplus b_{m}
\]
whenever 
\[
d_{\tau}(a) < \tau(b_{i})
\]
for every $\tau \in QT(A)$ and $i=0,\ldots,m$. 
\item $A$ has strong tracial $m$-comparison  of positive elements, if for any nonzero positive contractions $a,b \in M_{\infty}(A)$ we have 
\[
a \precsim  b
\]
whenever 
\[
d_{\tau}(a) < \frac{1}{m+1} \cdot \tau(b)
\]
for every $\tau \in QT(A)$. 
\end{enumerate}
\end{ndefn}
\en

\bn
\label{comparison-tracial-comparison-rem}
\begin{nremarks}
Let $A$ and $m$ be as above.
\begin{enumerate}
\item Note that $0$-comparison is the same as strict comparison as defined in \cite{Ror:Z-absorbing};  we will occasionally just write `comparison' for brevity.  Note also that tracial comparison and strong tracial comparison (i.e., in the case $m=0$) agree, and that strong tracial $m$-comparison implies tracial $m$-comparison. (We will derive a partial converse below.) 
\item The definition of $m$-comparison above differs slightly from that of \cite{Rob:dimnuc-comparison} since here we restrict ourselves to the case of simple $\mathrm{C}^{*}$-algebras. In this situation, it is not hard to show that Definition~\ref{strict-comparison}(i) and  \cite[Definition~2]{Rob:dimnuc-comparison} are in fact equivalent (using \cite[Lemma~1]{Rob:dimnuc-comparison} and an argument similar to that of Proposition~\ref{tau-dtau-eta} below).
\item I regard $m$-comparison as an algebraic interpretation of dimension. This is the reason why I write `$m$-comparison' when the  definition involves $m+1$ positive elements $b_{i}$, as this notation is more coherent with what is common in dimension theory. The same comment applies to (strong) tracial $m$-comparison and to $m$-almost  divisibility, see Definition~\ref{almost-divisible} below. 
\end{enumerate}
\end{nremarks}
\en

\bn
\label{comparison-tracial-comparison}
\begin{nprop}
If $A$ is separable, simple and unital, then $A$ has  $m$-comparison if and only if $A$ has tracial $m$-comparison. 
\end{nprop}

\begin{nproof}
We have $\tau(b) \le d_{\tau}(b)$ for every $\tau \in QT(A)$ and whenever $b $ is a positive contraction, so clearly $m$-comparison implies tracial $m$-comparison.

For the converse, let $a,b_{0},\ldots, b_{m} \in M_{\infty}(A)$ be positive contractions such that
\[
d_{\tau}(a) < d_{\tau}(b_{i})
\]
for every $\tau \in QT(A)$ and $i=0,\ldots,m$. By passing to a matrix algebra, we may assume that $a,b_{i} \in A$.

Let $\epsilon>0$, then 
\begin{eqnarray*}
d_{\tau}(f_{2\epsilon,3\epsilon}(a)) & \le & \tau(f_{\epsilon,2\epsilon}(a)) \\
& \le & d_{\tau}(a) \\
& < & d_{\tau}(b_{i}) \\
& = & \lim_{\eta \searrow 0} \tau(g_{\eta,2\eta}(b_{i})) 
\end{eqnarray*}
for every $\tau \in QT(A)$ and $i=0,\ldots,m$. By Dini's Theorem (cf.\ \ref{d-tau-notation1}) there is $\eta_{\epsilon}>0$ such that 
\[
\tau(f_{\epsilon,2\epsilon}(a)) < \tau(g_{\eta_{\epsilon},2\eta_{\epsilon}}(b_{i}))
\]
for every $\tau \in QT(A)$ and $i=0,\ldots,m$.
We obtain 
\[
d_{\tau}(f_{2\epsilon,3\epsilon}(a)) < \tau(g_{\eta_{\epsilon},2\eta_{\epsilon}}(b_{i}))
\]
for every $\tau \in QT(A)$ and $i=0,\ldots,m$, whence
\[
f_{2\epsilon,3\epsilon}(a) \precsim g_{\eta_{\epsilon},2\eta_{\epsilon}}(b_{0}) \oplus \ldots \oplus g_{\eta_{\epsilon},2\eta_{\epsilon}}(b_{m}) \precsim  b_{0} \oplus \ldots \oplus b_{m}
\]
by tracial $m$-comparison. Since $\epsilon>0$ was arbitrary, this yields 
\[
a \precsim  b_{0} \oplus \ldots \oplus b_{m},
\]
as desired.
\end{nproof}
\en

\bn
\label{tau-dtau-eta}
\begin{nprop}
Let $A$ be a separable, simple, unital $\mathrm{C}^{*}$-algebra and $m \in \mathbb{N}$. Suppose that, whenever $a,b \in M_{\infty}(A)$ are positive contractions and $\eta>0$ is a number such that 
\[
d_{\tau}(a) <\frac{1}{m+1} \cdot \tau(b) - \eta
\]
for every $\tau \in QT(A)$, then $a \precsim b$. 

Then, $A$ has strong tracial $m$-comparison.
\end{nprop}

\begin{nproof}
Let $a,b \in M_{\infty}(A)$ be positive contractions such that 
\[
d_{\tau}(a) < \frac{1}{m+1} \cdot \tau(b)
\]
for every $\tau \in QT(A)$. 

If $a$ is Cuntz equivalent to a projection, we may as well assume that $a$ is a projection itself. In this case, 
\[
\tau(a) = d_{\tau}(a) < \frac{1}{m+1} \cdot \tau(b)
\]
for every $\tau \in QT(A)$. Since $A$ is simple and unital, the assignment
\[
\tau \mapsto \frac{1}{m+1} \cdot \tau(b) - \tau(a)
\]
is a strictly positive continuous function on the compact space $QT(A)$, whence there is $\eta>0$ such that 
\[
\frac{1}{m+1} \cdot \tau(b) - d_{\tau}(a) = \frac{1}{m+1} \cdot \tau(b) - \tau(a) > \eta
\]
for all $\tau \in QT(A)$ (cf.\ \ref{traces-not}); we now have $a \precsim b$ by our hypothesis. 

If $a$ is not equivalent to a projection, then $0$ is an accumulation point of its spectrum and $a- f_{\epsilon,2\epsilon}(a)$ is positive and nonzero for any $\epsilon>0$. But then 
\[
\tau(a- f_{\epsilon,2\epsilon}(a)) > \eta_{\epsilon}
\]
for some $\eta_{\epsilon}>0$ and for all $\tau \in QT(A)$ (again cf.\ \ref{traces-not}).

We now estimate
\begin{eqnarray*}
d_{\tau}(f_{2\epsilon,3\epsilon}(a)) + \eta_{\epsilon} & < & d_{\tau}(f_{2\epsilon,3\epsilon}(a)) + \tau(a- f_{\epsilon,2\epsilon}(a)) \\
& \le & d_{\tau}(f_{2\epsilon,3\epsilon}(a)) + d_{\tau}(a- f_{\epsilon,2\epsilon}(a)) \\
& \le & d_{\tau}(a) \\
& < & \frac{1}{m+1} \cdot \tau(b),
\end{eqnarray*}
using the fact that $f_{2\epsilon,3\epsilon}(a)$ and $a- f_{\epsilon,2\epsilon}(a)$ are orthogonal, with sum dominated by $a$. By our hypothesis we then have $f_{2\epsilon,3\epsilon}(a) \precsim b$; since $\epsilon>0$ was arbitrary, this implies $a \precsim b$. 
\end{nproof}
\en

\bn
\label{almost-divisible}
\begin{ndefn}
Let  $A$ be a unital $\mathrm{C}^{*}$-algebra and $m \in \mathbb{N}$. 
\begin{enumerate}
\item $A$ has $m$-almost divisible Cuntz semigroup, if for any positive contraction $a \in M_{\infty}(A)$ and $0 \neq k \in \mathbb{N}$ there is $x \in W(A)$ such that 
\[
k \cdot x \le \langle a \rangle \le (k+1) (m+1) \cdot x.
\]
\item We say $A$ is tracially $m$-almost divisible, if for any positive contraction $a \in M_{\infty}(A)$, $\epsilon>0$ and $0 \neq k \in \mathbb{N}$ there is a c.p.c.\ order zero map 
\[
\psi:M_{k} \to \her(a) \subset  M_{\infty}(A)
\]
such that 
\[
\tau(\psi(\be_{k})) \ge \frac{1}{m+1}\cdot  \tau(a) - \epsilon 
\]
for all $\tau \in QT(A)$.
\end{enumerate}
\end{ndefn}
\en

\bn
\label{purely-finite}
\begin{ndefn}
Let $A$ be a separable, simple, unital $\mathrm{C}^{*}$-algebra and $m,\bar{m} \in \mathbb{N}$. 
\begin{enumerate}
\item We say $A$ is $(m,\bar{m})$-pure, if it has $m$-comparison of positive elements and if it is  $\bar{m}$-almost divisible.  We say $A$ is pure, if it is $(0,0)$-pure.
\item  We say $A$ is tracially $(m,\bar{m})$-pure, if it has strong tracial $m$-comparison of positive elements and if it is tracially $\bar{m}$-almost divisible.  We say $A$ is tracially pure, if it is tracially $(0,0)$-pure.
\end{enumerate}
\end{ndefn}
\en

\bn
\label{W-Z-purely-finite}
For  convenience, let us note the following combination of results by R{\o}rdam explicitly.

\begin{nprop}
A separable, simple, unital and $\mathcal{Z}$-stable $\mathrm{C}^{*}$-algebra is pure.
\end{nprop}

\begin{nproof}
By \cite[Theorem~4.5]{Ror:Z-absorbing},  a $\mathcal{Z}$-stable $\mathrm{C}^{*}$-algebra  $A$ has almost unperforated  Cuntz semigroup, hence comparison (see, for example, \cite[Proposition~2.1]{OrtPerRor:cfp-dr}). 

To show almost divisibility, let $0 \neq k \in \mathbb{N}$ and let $a \in M_{\infty}(A)$ be a positive contraction. By \cite[Lemma~4.2]{Ror:Z-absorbing},  there is a positive contraction $e \in \mathcal{Z}$ such that 
\[
k \cdot \langle e \rangle \le \be_{\mathcal{Z}} \le (k+1) \cdot \langle e \rangle. 
\]
Now $a \otimes e \in M_{\infty}(A \otimes \mathcal{Z}) $ will satisfy
\begin{equation}
\label{www-2-1}
k \cdot \langle a \otimes e \rangle \le a \otimes \be_{\mathcal{Z}} \le (k+1) \cdot \langle a \otimes e \rangle 
\end{equation}
in  $W(A \otimes \mathcal{Z})$. 

By \cite{TomsWin:ssa} there is an isomorphism 
\[
\varphi:A \to A \otimes \mathcal{Z}
\]
which is approximately unitarily equivalent to the embedding 
\[
\id_{A} \otimes \be_{\mathcal{Z}}: A \to A \otimes \mathcal{Z}.
\]
This in turn implies that  $\id_{A} \otimes \be_{\mathcal{Z}}$ and $\varphi$ both induce the same isomorphism between Cuntz semigroups, say
\[
\psi:W(A) \stackrel{\cong}{\longrightarrow} W(A \otimes \mathcal{Z}),
\]
whence
\begin{equation}
\label{www-2-2}
\psi^{-1}(\langle a \otimes \be_{\mathcal{Z}} \rangle) = \langle a \rangle.
\end{equation}
It follows from \eqref{www-2-1} and \eqref{www-2-2} that 
\[
x:= \psi^{-1}(\langle a \otimes e \rangle) \in W(A)
\]
satisfies 
\[
k\cdot x \le \langle a \rangle \le (k+1) \cdot x,
\]
as desired.
\end{nproof}
\en

\bn
\label{divisible-tracially-divisible}
\begin{nprop}
Let $A$ be a separable, simple, unital $\mathrm{C}^{*}$-algebra and $m \in \mathbb{N}$. 

If $A$ has $m$-almost divisible Cuntz semigroup, then $A$ is tracially $m$-almost divisible.
\end{nprop}

\begin{nproof}
Let $a \in M_{\infty}(A)$ a positive contraction, $\epsilon>0$ and $0 \neq k \in \mathbb{N}$ be given. By passing to a matrix algebra over $A$, we may as well assume that $a \in A$. 

Choose $\bar{k} \in \mathbb{N}$ such that $k$ divides $\bar{k}$ and such that 
\[
\frac{1}{\bar{k}} < \frac{\epsilon}{2}.
\]  
By $m$-almost divisibility, there is $x \in W(A)$ such that 
\begin{equation}
\label{w2-1}
\bar{k} \cdot x \le \langle a \rangle \le (\bar{k} +1) (m+1) \cdot x.
\end{equation}
Suppose that $x = \langle b \rangle$ for some positive contraction $b \in M_{r}(A)$ and $r \in \mathbb{N}$. We then have
\[
\frac{1}{m+1} \cdot \tau(a) \le \frac{1}{m+1} \cdot d_{\tau}(a) \stackrel{\eqref{w2-1}}{\le} (\bar{k}+1) \cdot d_{\tau}(b)
\]
for every $\tau \in QT(A)$. 

Since 
\[
d_{\tau}(b) = \lim_{\eta \searrow 0} \tau(g_{\eta,2\eta}(b))
\]
for every $\tau$, by  Dini's Theorem (cf.\ \ref{d-tau-notation1}) there is $\eta>0$ such that
\begin{equation}
\label{w2-2}
\frac{1}{m+1} \cdot \tau(a) \le (\bar{k}+1) \cdot \tau(g_{\eta,2\eta}(b)) +\frac{\epsilon}{2}
\end{equation}
for every  $\tau \in QT(A)$. We also have 
\begin{equation}
\label{w2-3}
\tau(g_{\eta,2\eta}(b)) \le d_{\tau}(b) \stackrel{\eqref{w2-1}}{\le} \frac{1}{\bar{k}} \cdot d_{\tau}(a) \le \frac{1}{\bar{k}} < \frac{\epsilon}{2},
\end{equation}
whence
\begin{equation}
\label{w2-4}
\frac{1}{m+1} \cdot \tau(a)  \stackrel{\eqref{w2-2},\eqref{w2-3}}{<} \bar{k} \cdot \tau(g_{\eta,2\eta}(b)) +\epsilon
\end{equation}
for every $\tau \in QT(A)$. 

Since 
\[
\bar{k}\cdot x \le \langle a \rangle, 
\] 
by Proposition~\ref{cc-order-zero} there is a c.p.c.\ order zero map
\[
\bar{\psi}:M_{\bar{k}} \to \her(a) \subset A
\]
such that 
\[
\bar{\psi}(e_{11}) \approx g_{\eta,2\eta}(b),
\]
which in turn implies that 
\[
\tau(\bar{\psi}(\be_{\bar{k}})) = \bar{k} \cdot \tau(\bar{\psi}(e_{11})) = \bar{k} \cdot \tau(g_{\eta,2\eta}(b)) \stackrel{\eqref{w2-4}}{>} \frac{1}{m+1} \cdot \tau(a) - \epsilon
\]
for every $\tau \in QT(A)$. 

Upon composing $\bar{\psi}$ with a unital embedding of $M_{k}$ into $M_{\bar{k}}$ ($k$ divides $\bar{k}$), we obtain $\psi$ as desired. 
\end{nproof}
\en

\bn
\label{pure-comparison}
\begin{nprop}
Given $m,\bar{m} \in \mathbb{N}$, there is $\tilde{m} \in \mathbb{N}$ such that the following holds: 

If $A$ is a separable, simple, unital  $\mathrm{C}^{*}$-algebra which has tracial $m$-comparison and tracial $\bar{m}$-almost divisibility, then $A$ has strong tracial $\tilde{m}$-comparison, i.e., $A$ is $(\tilde{m},\bar{m})$-pure. 
\end{nprop}

\begin{nproof}
Take 
\begin{equation}
\label{www-2-3}
\tilde{m} := (m+1) (\bar{m}+1) -1.
\end{equation}
Now suppose we have positive contractions $a,b \in A$ and some $\eta>0$ with
\begin{equation}
\label{pure-comparison1}
d_{\tau}(a) < \frac{1}{\tilde{m}+1} \cdot \tau(b) - \eta
\end{equation}
for all $\tau \in QT(A)$. By Proposition~\ref{tau-dtau-eta} it will suffice to show that $a \precsim b$. (We do not have to consider matrix algebras, since any matrix algebra over $A$ satisfies the hypotheses of the proposition just as $A$ does.) 

By tracial $\bar{m}$-almost divisibility, there is a c.p.c.\ order zero map 
\[
\psi:M_{m+1} \to \her(b) \subset A
\] 
such that
\begin{equation}
\label{pure-comparison2}
\tau(\psi(\be_{m+1})) \ge \frac{1}{\bar{m}+1} \cdot \tau(b) - \frac{\tilde{m}+1}{\bar{m}+1} \cdot \eta
\end{equation}
for all $\tau \in QT(A)$. Then \eqref{pure-comparison1} and \eqref{pure-comparison2} together give
\begin{eqnarray*}
d_{\tau}(a) & < & \frac{\bar{m} +1}{\tilde{m}+1} \cdot \tau(\psi(\be_{m+1})) \\
&  \stackrel{\ref{order-zero-traces}}{=} & \frac{(\bar{m} +1)(m+1)}{(\tilde{m}+1)} \cdot \tau(\psi(e_{11})) \\
& = & \tau(\psi(e_{11})) 
\end{eqnarray*}
for all $\tau \in QT(A)$; by tracial $m$-comparison and since $\psi$ has order zero, we obtain 
\[
a \precsim \psi(\be_{m+1}) \precsim b,
\]
as desired.
\end{nproof}
\en

\bn
\label{pure-tracially-pure-cor}
Let us explicitly note the following corollary, in which $\tilde{m}$ comes from Proposition~\ref{pure-comparison}.

\begin{ncor}
If $A$ is separable, simple, unital and $(m,\bar{m})$-pure, then $A$ is tracially $(\tilde{m},\bar{m})$-pure.
\end{ncor}

\begin{nproof}
This combines Propositions~\ref{comparison-tracial-comparison}, \ref{divisible-tracially-divisible} and \ref{pure-comparison}.
\end{nproof}
\en

\bn
The final result of this section says that pureness and tracial pureness are in fact equivalent. 
At this point I do not know whether Corollary~\ref{pure-tracially-pure-cor} admits a more general  converse.

\begin{ntheorem}
Let $A$ be a separable, simple and unital $\mathrm{C}^{*}$-algebra.

Then, $A$ is pure if and only if $A$ is tracially pure.
\end{ntheorem}

\begin{nproof}
The forward implication is a special case of Corollary~\ref{pure-tracially-pure-cor}, observing that Proposition~\ref{pure-comparison} yields $\tilde{m}=0$ if $m=\bar{m}=0$ by \eqref{www-2-3}.

So let us assume that $A$ is tracially pure, i.e., $A$ has strong tracial comparison and is tracially almost divisible.

But strong tracial comparison is just tracial comparison, which implies comparison by Proposition~\ref{comparison-tracial-comparison}.  It therefore remains to show that  $W(A)$ is almost divisible.  

To this end, let $a \in M_{\infty}(A)$ be a (nonzero) positive contraction and $0\neq k \in \mathbb{N}$ a number; we need to find $x \in W(A)$ such that 
\[
k \cdot x \le \langle a \rangle \le (k+1) \cdot x.
\]
By passing to a matrix algebra over $A$, we may clearly assume that $a \in A$. 

We first consider the case when $a$ is Cuntz equivalent to a projection, so that we may as well assume that $a$ itself is a projection. But then 
\begin{equation}
\label{w2-7}
\tau(a) = d_{\tau}(a)
\end{equation}
for all $\tau \in QT(A)$; by \ref{traces-not} there is $\delta>0$ such that 
\begin{equation}
\label{w2-6}
d_{\tau}(a) = \tau(a)>\delta
\end{equation}
for every $\tau \in QT(A)$. 

By tracial almost divisibility there is a c.p.c.\ order zero map 
\[
\psi:M_{k} \to \her(a) \subset A
\]
such that
\begin{equation}
\label{w2-5}
\tau(\psi(\be_{k})) \ge \tau(a) - \frac{\delta}{k+1} = d_{\tau}(a) - \frac{\delta}{k+1} 
\end{equation}
for every $\tau \in QT(A)$. But then we have
\begin{eqnarray}
(k+1) \cdot d_{\tau}(\psi(e_{11})) & \ge & (k+1) \cdot \tau(\psi(e_{11})) \nonumber \\
& \stackrel{\ref{order-zero-traces}}{=} &  \frac{(k+1)}{k} \cdot \tau(\psi(\be_{k})) \nonumber \\
& \stackrel{\eqref{w2-5}}{\ge} & \frac{1}{k} ((k+1) \cdot \tau(a) - \delta) \nonumber \\
& \stackrel{\eqref{w2-7}}{=} & d_{\tau}(a) + \frac{1}{k} \cdot (d_{\tau}(a) - \delta) \nonumber \\
& \stackrel{\eqref{w2-6}}{>} & d_{\tau}(a)  \label{w2-8}
\end{eqnarray}
for every $\tau \in QT(A)$. Setting 
\[
x:= \langle \psi(e_{11}) \rangle,
\]
we obtain 
\[
d_{\tau}((k+1) \cdot x) = (k+1) \cdot d_{\tau}(x) \stackrel{\eqref{w2-8}}{>} d_{\tau}(a)
\]
for all $\tau$, whence 
\[
\langle a \rangle \le (k+1) \cdot x
\]
by comparison. We have 
\[
k \cdot x \le \langle a \rangle
\]
since $\psi(\be_{k}) \in \her (a)$. 

We now consider the case when $a$ is not Cuntz equivalent to a projection. But then $0$ is an accumulation point of the spectrum, and (using \ref{traces-not}) it is straightforward to construct a strictly decreasing sequence $(\eta_{l})_{l \in \mathbb{N}}$ of strictly positive numbers such that the elements
\[
a_{l}:= g_{\eta_{l},2\eta_{l}}(a) \in A
\]
form a strictly increasing sequence of positive contractions, satisfying
\begin{equation}
\label{w2-9}
\tau(a_{l}) > d_{\tau}(a_{l-1}) > \tau(a_{l-1})
\end{equation}
for all $l \ge 1$ and $\tau \in QT(A)$; note that
\begin{equation}
\label{w2-11}
d_{\tau}(a_{l}) \stackrel{l \to \infty}{\longrightarrow} d_{\tau}(a)
\end{equation}
for all $\tau$ by lower semicontinuity. 

For $1 \le l \in \mathbb{N}$, set 
\begin{equation}
\label{w2-13}
\epsilon_{l}:= \frac{1}{2} \cdot (\tau(a_{l}) - d_{\tau}(a_{l-1})) \stackrel{\eqref{w2-9}}{>} 0.
\end{equation}
By tracial almost divisibility, for each $l \in \mathbb{N}$ there is a c.p.c.\ order zero map
\begin{equation}
\label{w2-10}
\psi_{l}: M_{k} \to \her(a_{l})
\end{equation}
such that 
\begin{equation}
\label{w2-12}
\tau(\psi_{l}(\be_{k})) \ge \tau(a_{l}) - \epsilon_{l}
\end{equation}
for every $\tau \in QT(A)$. We then have
\begin{eqnarray}
d_{\tau}(a) & \stackrel{\eqref{w2-11},\eqref{w2-10}}{>} & d_{\tau}(\psi_{l}(\be_{k})) \nonumber \\
& \ge & \tau(\psi_{l}(\be_{k})) \nonumber \\
& \stackrel{\eqref{w2-12}}{\ge} & \tau(a_{l}) - \epsilon_{l} \nonumber \\
& \stackrel{\eqref{w2-13}}{>} & d_{\tau}(a_{l-1}) \nonumber \\
& \stackrel{\eqref{w2-10}}{\ge} & d_{\tau}(\psi_{l-1}(\be_{k}))  \label{w2-14}
\end{eqnarray}
for all $ l \ge 1$ and $\tau \in QT(A)$. 

Set 
\[
x_{l}:= \langle \psi_{l}(e_{11}) \rangle \in W(A),
\]
then 
\begin{equation}
\label{w2-15}
d_{\tau}(k \cdot x_{l}) = k \cdot d_{\tau}(x_{l}) \stackrel{\ref{order-zero-traces}}{=} d_{\tau}(\psi_{l}(\be_{k}))
\end{equation}
for all $l$ and $\tau$; by \eqref{w2-14} we have 
\[
d_{\tau}(x_{l}) < d_{\tau}(x_{l+1})
\]
for all $l$ and $\tau$, so comparison yields
\[
x_{l} \le x_{l+1} \le \ldots
\]
and
\[
k \cdot x_{l} \le k \cdot x_{l+1} \le \ldots \stackrel{\eqref{w2-14},\eqref{w2-15}}{\le} \langle a \rangle.
\]
Let 
\[
x:= \sup \{x_{l} \} \in W(A \otimes \mathcal{K}) 
\]
be the supremum of the $x_{l}$, cf.\ \cite{CEI:Cu} or \cite{BroPerToms:cuntz-semigroup}; it follows from \cite{BRTTW:Cu} that in fact 
\[
x \in W(A) \subset W(A \otimes \mathcal{K}).
\]
By \cite{CEI:Cu} we also have 
\[
k \cdot x = \sup\{k \cdot x_{l}\} \le \langle a \rangle.
\]
Moreover, we have
\[
d_{\tau}(x) = \frac{1}{k} \cdot d_{\tau}(k \cdot x) \ge \frac{1}{k} \cdot d_{\tau}(\psi_{l}(\be_{k})) \stackrel{l \to \infty}{\longrightarrow} \frac{1}{k} \cdot d_{\tau}(a)>0,
\]
whence
\[
d_{\tau}(a) < (k+1) \cdot d_{\tau}(x) = d_{\tau}((k+1) \cdot x)
\]
for all $\tau \in QT(A)$. Comparison yields
\[
\langle a \rangle \le (k+1) \cdot x.
\]
We have now found $x \in W(A)$ such that 
\[
k \cdot x \le \langle a \rangle \le (k+1) \cdot x,
\]
regardless of whether $a$ is equivalent to a projection or not. 
\end{nproof}
\en

\section{pureness and nuclear dimension}
\label{dimnuc-pureness}

\noindent
Corollary~\ref{findimnuc-Z-stable} below says that separable, simple, nonelementary, unital $\mathrm{C}^{*}$-algebras with finite nuclear dimension are $\mathcal{Z}$-stable;  together with Proposition~\ref{W-Z-purely-finite} this entails   pureness.

As an intermediate step towards $\mathcal{Z}$-stability we prove  in this section that, for separable, simple, unital $\mathrm{C}^{*}$-algebras, finite nuclear dimension implies tracial $(m,\bar{m})$-pureness (for suitable $m, \bar{m} \in \mathbb{N}$). The key features here are $m$-comparison as established in \cite{Rob:dimnuc-comparison}, and finite covering number of the central sequence algebra (see \cite{Kir:CentralSequences}) of a $\mathrm{C}^{*}$-algebra of nuclear dimension at most $m$.

\bn
The following definitions are collected from  \cite{KirWinter:dr}, \cite{Winter:lfdr} and \cite{WinterZac:dimnuc}. Note that in (iii), we may additionally ask for the composition $\varphi \circ \psi$ to be contractive by \cite[Remark~2.2(iv)]{WinterZac:dimnuc}. 

\label{d-dr} 
\begin{ndefn}
Let $A$ be a $\mathrm{C}^{*}$-algebra, $m \in \mathbb{N}$.
\begin{itemize}
\item[(i)] A c.p.\ map $\varphi : F \to A$  is $m$-decomposable (where $F$ is a finite-dimensional $\mathrm{C}^{*}$-algebra), if there is a decomposition 
\[
F=F^{(0)} \oplus \ldots \oplus F^{(m)}
\]
such that the restriction $\varphi^{(i)}$ of $\varphi$ to $F^{(i)}$ has   order zero for each $i \in \{0, \ldots, m\}$; we say $\varphi$ is $m$-decomposable with respect to $F=F^{(0)} \oplus \ldots \oplus F^{(m)}$.
\item[(ii)] $A$ has decomposition rank $m$, $\dr A = m$, if $m$ is the least integer such that the following holds: For any finite subset $\Gh  \subset A$ and $\varepsilon > 0$, there is a finite-dimensional c.p.c.\ approximation $(F, \psi, \varphi)$ for $\Gh$ to within $\varepsilon$ (i.e., $F$ is finite-dimensional $\psi:A \to F$ and $\varphi:F \to A$ are c.p.c.\ and $\|\varphi \psi (b) - b\| < \varepsilon \; \forall \, b \in \Gh$) such that $\varphi$ is $m$-decomposable. If no such $m$ exists, we write $\dr A = \infty$.  
\item[(iii)] $A$ has nuclear dimension $m$, $\dimnuc A = m$, if $m$ is the least integer such that the following holds: For any finite subset $\Gh  \subset A$ and $\varepsilon > 0$, there is a finite-dimensional c.p.\  approximation $(F, \psi, \varphi)$ for $\Gh$ to within $\varepsilon$ (i.e., $F$ is finite-dimensional, $\psi:A \to F$ and $\varphi:F \to A$ are c.p.\  and $\|\varphi \psi (b) - b\| < \varepsilon \; \forall \, b \in \Gh$) such that $\psi$ is c.p.c., and $\varphi$ is $m$-decomposable with c.p.c.\ order zero components $\varphi^{(i)}$. If no such $m$ exists, we write $\dimnuc A = \infty$. 
\item[(iv)] $A$ has locally finite decomposition rank (or locally finite nuclear dimension, respectively), if for any finite subset $\Gh  \subset A$ and $\varepsilon > 0$, there is a $\mathrm{C}^{*}$-subalgebra $B \subset A$ such that $\dr B < \infty$ (or $\dimnuc B < \infty$, respectively) and $\mathcal{G} \subset_{\epsilon} B$.
\end{itemize}
\end{ndefn}
\en

\bn
\label{central-cutout}
The following is contained in the proof of \cite[Proposition~4.3]{WinterZac:dimnuc}. An explicit proof is extracted here for  convenience of the reader.

\begin{nprop}
Let $A$ be a unital and separable $\mathrm{C}^{*}$-algebra with $\dimnuc A \le m < \infty$.  

Then, there is a system of finite-dimensional $m$-decomposable c.p.\ approximations for $A$
\[
(F_{p}=F_{p}^{(0)} \oplus \ldots \oplus F_{p}^{(m)}, \psi_{p},\varphi_{p})_{p \in \mathbb{N}}
\]
such that the 
\[
\psi_{p}:A \to F_{p}
\]
are c.p.c.\ maps,  the 
\[
\varphi^{(i)}_{p}:F^{(i)}_{p} \to A
\]
are c.p.c.\ order zero maps, and such that 
\begin{equation}
\label{central-cutout-1}
\|\varphi_{p} \psi_{p}(b) - b\| \stackrel{p \to \infty}{\longrightarrow} 0 ,
\end{equation}
\begin{equation}
\label{central-cutout-2}
\|[\varphi_{p} \psi_{p}(b) , \varphi^{(i)}_{p}(\psi^{(i)}_{p}(\be_{A})) ]\| \stackrel{p \to \infty}{\longrightarrow} 0
\end{equation}
and
\begin{equation}
\label{central-cutout-3}
\|b \varphi^{(i)}_{p}(\psi^{(i)}_{p}(\be_{A})) - \varphi^{(i)}_{p} \psi^{(i)}_{p}(b) \| \stackrel{p \to \infty}{\longrightarrow} 0
\end{equation}
for all $b \in A$, $i \in \{0, \ldots,m\}$. 
\end{nprop}

\begin{nproof}
Choose a system 
\[
\left( F_{p}= F_{p}^{(0)}\oplus \ldots \oplus F_{p}^{(m)}, \psi_{p}, \varphi_{p}= \sum_{i=0}^{m} \varphi^{(i)}_{p} \right)_{p \in \mathbb{N}}
\]
of  $m$-decomposable c.p.\ approximations for $A$ such that \eqref{central-cutout-1} holds and the maps 
\[
\psi_{p}, \varphi_{p}^{(i)}, \varphi_{p} \psi_{p}
\] 
are all c.p.\ contractions. 

For each $p \in \mathbb{N}$, we define
\[
\hat{\psi}_{p}( \, . \,) := \psi_{p}(\be_{A})^{-\frac{1}{2}} \psi_{p}(\, . \,) \psi_{p}(\be_{A})^{-\frac{1}{2}},
\]   
where the inverses are taken in the hereditary subalgebras 
\[
\tilde{F}_{p} := \her(\psi_{p}(\be_{A})) \subset F_{p}
\]
and 
\[
\hat{\varphi}_{p}( \, . \,) := \varphi_{p}(\psi_{p}(\be_{A})^{\frac{1}{2}} \, . \, \psi_{p}(\be_{A})^{\frac{1}{2}}),
\]   
then
\[
\hat{\varphi}_{p} \hat{\psi}_{p} = \varphi_{p} \psi_{p};
\]
moreover, the 
\[
\hat{\psi}_{p}: A \to \tilde{F}_{p}
\] 
are unital c.p.\  and the 
\[
\hat{\varphi}_{p}:\tilde{F}_{p} \to A
\]
are c.p.c.\ maps.

From \cite[Lemma~3.6]{KirWinter:dr}, we see that for $i \in \{0, \ldots,m\}$,  $p \in \mathbb{N}$ and any  projection $q_{p}  \in \tilde{F}_{p}$ we have
\begin{eqnarray*}
\lefteqn{ \|\hat{\varphi}_{p}(q_{p} \hat{\psi}_{p}(a)) - \hat{\varphi}_{p}(q_{p}) \hat{\varphi}_{p} \hat{\psi}_{p} (a)\| } \\
& \le &  3 \cdot \max \{\| \hat{\varphi}_{p} \hat{\psi}_{p} (a) - a\|, \, \| \hat{\varphi}_{p} \hat{\psi}_{p} (a^{2}) - a^{2}\| \}^{\frac{1}{2}},
\end{eqnarray*}
from which follows that 
\[
\| \varphi^{(i)}_{p} \psi^{(i)}_{p} (a) - \varphi^{(i)}_{p} \psi^{(i)}_{p}(\be_{A}) \varphi_{p} \psi_{p}(a) \| \stackrel{p \to \infty}{\longrightarrow} 0
\]
for any $a \in A$, i.e., \eqref{central-cutout-3} holds. From this we  obtain \eqref{central-cutout-2} by taking adjoints. 
\end{nproof}
\en

\bn
\label{quasitraces-traces}
\begin{nremark}
A celebrated result of U.\ Haagerup says that every quasitrace on an exact $\mathrm{C}^{*}$-algebra is a trace, see \cite{Haa:quasi}. In \cite{BroWinter:dimnuc-quasitraces}, N.~Brown and the author  give a short  proof of this fact in the case of $\mathrm{C}^{*}$-algebras with locally finite nuclear dimension.
\end{nremark}
\en

\bn
\label{Kirchberg-cov}
In \cite[Proposition~4.3]{WinterZac:dimnuc} it was shown that if $A$ is separable, simple, nonelementary and unital with $\dimnuc A\le m$, then $A_{\omega} \cap A'$ (where $\omega$ is a free ultrafilter on $\mathbb{N}$) has Kichberg covering number (in the sense of \cite{Kir:CentralSequences}) at most $(m+1)^{2}$.  If $A$ is unital, then the proof   shows that in fact $A_{\infty} \cap A'$ has covering number at most $(m+1)^{2}$. We note here a straightforward generalization of that result. The proof involves a standard diagonal sequence argument, which is spelled out explicitly since it will appear repeatedly in the sequel. 

\begin{nprop}
Let $A$ be a separable, simple, nonelementary and unital $\mathrm{C}^{*}$-algebra with $\dimnuc A \le m < \infty$; let $X \subset A_{\infty}$ be a separable subspace. 

Then, for any $k \in \mathbb{N}$ there are c.p.c.\ order zero maps 
\begin{equation}
\label{Kirchberg-cov-1}
\psi^{(1)},\ldots, \psi^{((m+1)^{2})}:M_{k} \oplus M_{k+1} \to A_{\infty} \cap A' \cap X'
\end{equation}
satisfying 
\begin{equation}
\label{Kirchberg-cov-2}
\sum_{j=1}^{(m+1)^{2}} \psi^{(j)}(\be_{k} \oplus \be_{k+1}) \ge \be_{A_{\infty}}.
\end{equation}
\end{nprop}

\begin{nproof}
By (the proof of) \cite[Proposition~4.3]{WinterZac:dimnuc}, there are c.p.c.\ order zero maps
\begin{equation}
\label{Kirchberg-cov-3}
\varphi^{(1)},\ldots, \varphi^{((m+1)^{2})}:M_{k} \oplus M_{k+1} \to A_{\infty} \cap A' 
\end{equation}
satisfying 
\begin{equation}
\label{Kirchberg-cov-4}
\sum_{j=1}^{(m+1)^{2}} \varphi^{(j)}(\be_{k} \oplus \be_{k+1}) \ge \be_{A_{\infty}}.
\end{equation}
For each $j$, let 
\[
\bar{\varphi}^{(j)}: M_{k} \oplus M_{k+1} \to \prod_{\mathbb{N}} A
\]
be a c.p.c.\ order zero lift of $\varphi^{(j)}$ (cf.\ \cite[Remark~2.4]{KirWinter:dr}) with c.p.c.\ order zero components
\[
\bar{\varphi}^{(j)}_{l}: M_{k} \oplus M_{k+1} \to A, \, l \in \mathbb{N}. 
\]
Let 
\[
\{x_{i} \mid i \in \mathbb{N}\} \subset X 
\]
be a countable dense subset; for each $i \in \mathbb{N}$ let 
\[
\bar{x}_{i}= (\bar{x}_{i,r})_{r \in \mathbb{N}} \in \prod_{\mathbb{N}} A
\]  
be an isometric lift of $x_{i}$.  

Since $\bar{\varphi}^{(j)}$ lifts $\varphi^{(j)}$, using \eqref{Kirchberg-cov-3} it is straightforward to find a strictly  increasing sequence $(l_{r})_{r \in \mathbb{N}} \subset \mathbb{N}$ such that for each $r \in \mathbb{N}$, $i \le r$ and $y \in M_{k} \oplus M_{k+1}$ we have
\begin{equation}
\label{Kirchberg-cov-5}
\| [ \bar{\varphi}^{(j)}_{l_{r}}(y), \bar{x}_{i,r} ]\| \le \frac{1}{r+1} \|y\|.
\end{equation}
We obtain c.p.c.\ order zero maps 
\[
\bar{\psi}^{(j)}:= (\bar{\varphi}^{(j)}_{l_{r}})_{r \in \mathbb{N}}: M_{k} \oplus M_{k+1} \to \prod_{\mathbb{N}} A
\] 
which induce c.p.c.\ order zero maps
\[
\psi^{(j)}: M_{k} \oplus M_{k+1} \to A_{\infty}, \, j = 1, \ldots, (m+1)^{2}. 
\]
The images of the $\psi^{(j)}$ commute with $X$ by \eqref{Kirchberg-cov-5}, they commute with $A$ by \eqref{Kirchberg-cov-3};  \eqref{Kirchberg-cov-2} follows directly  from \eqref{Kirchberg-cov-4}.
\end{nproof}
\en

\bn
\label{dimnuc-comparison}
It was shown in \cite[Lemma~6.1]{TomsWinter:VI} that decomposition rank at most $m$ implies $m$-comparison (in the simple and unital case; it was observed in \cite{OrtPerRor:cfp-dr} that essentially the same proof holds also in the general situation). Let us recall  a result of L.\ Robert, \cite[Theorem~1]{Rob:dimnuc-comparison}, which establishes the respective statement for nuclear dimension.  Note that the theorem indeed applies since our definition of $m$-comparison agrees with that of \cite{Rob:dimnuc-comparison} in the simple situation, cf.\ \ref{comparison-tracial-comparison-rem}(ii).

\begin{ntheorem} {\rm (L.\ Robert)}
If $A$ is a separable $\mathrm{C}^{*}$-algebra with $\dimnuc A \le m$, then $A$ has $m$-comparison.
\end{ntheorem}
\en

\bn
\label{2-d-i}
The following proposition nicely illustrates a technique which will be used frequently in the sequel, in particular in Section~{\ref{tracial-matrix-cone-absorption}}.

\begin{nprop}
Given $m,k \in \mathbb{N}$ (with $k\ge 1$), there is $1 \le L_{m,k} \in \mathbb{N}$ such that the following holds:

If $A$ is a separable, simple, nonelementary and unital $\mathrm{C}^{*}$-algebra with $\dimnuc A \le m$, and if $X \subset (A_{\infty})_{+}$ is a separable subspace, then there are pairwise orthogonal positive contractions 
\[
d^{(1)}, \ldots, d^{(k)} \in A_{\infty} \cap A' \cap X'
\]
such that 
\begin{equation}
\label{w3-15}
\tau(d^{(i)}b ) \ge \frac{1}{L_{m,k}} \cdot \tau(b)
\end{equation}
for all $b \in \mathrm{C}^{*}(A,X)_{+} \subset A_{\infty}$ and $\tau \in T_{\infty}(A)$. 
\end{nprop}

\begin{nproof}
For $k=1$ we may take $L_{m,1}:= 1$ and $d^{(1)}:= \be_{A_{\infty}}$, so let us  prove the statement for $k=2$.

Take 
\begin{equation}
\label{ww3-1}
L_{m,2}:= 2 (m+1)^{2} +1;
\end{equation}
for notational convenience set
\begin{equation}
\label{w3-4}
\bar{k}:= 2 (m+1)^{2}.
\end{equation}
By Proposition~\ref{Kirchberg-cov} there are c.p.c.\ order zero maps 
\[
\psi^{(1)},\ldots, \psi^{((m+1)^{2})}:M_{\bar{k}} \oplus M_{\bar{k}+1} \to A_{\infty} \cap A' \cap X'
\]
satisfying 
\[
\sum_{j=1}^{(m+1)^{2}} \psi^{(j)}(\be_{\bar{k}} \oplus \be_{\bar{k}+1}) \ge \be_{A_{\infty}}.
\]
We then have
\begin{equation}
\label{w3-1}
\tau\left(\sum_{j=1}^{(m+1)^{2}} \psi^{(j)}(e_{11}\oplus e_{11})b \right) \ge \frac{1}{\bar{k}+1} \cdot \tau(b)
\end{equation}
for all $\tau \in T_{\infty}(A)$ and $b \in \mathrm{C}^{*}(A,X)_{+}$ (using the fact that $\psi^{(j)}(\, . \,)b$ is a c.p.c.\ order zero map, and that the composition of an order zero map with a trace is again a trace by  \ref{order-zero-traces}).

For $\eta>0$, define
\begin{equation}
\label{w3-2}
d^{(0)}_{\eta}:= g_{\eta,2\eta}\left( \sum_{j=1}^{(m+1)^{2}} \psi^{(j)}(e_{11} \oplus e_{11})\right)
\end{equation}
and
\begin{equation}
\label{w3-3}
d^{(1)}_{\eta}:= \be_{A_{\infty}} - g_{0,\eta}\left( \sum_{j=1}^{(m+1)^{2}} \psi^{(j)}(e_{11} \oplus e_{11})\right).
\end{equation}
We have
\[
d^{(0)}_{\eta} \perp d^{(1)}_{\eta},
\]
\[
d^{(0)}_{\eta}, d^{(1)}_{\eta} \in A_{\infty} \cap A' \cap X'
\]
and
\begin{equation}
\label{w3-5}
\tau(d^{(0)}_{\eta} b) \stackrel{\eqref{w3-1},\eqref{w3-2}}{\ge} \frac{1}{\bar{k}+1} \cdot \tau(b) - \eta
\end{equation}
for all $\tau \in T_{\infty}(A)$ and all $b \in \mathrm{C}^{*}(A,X)_{+}$ of norm at most one.  Moreover, for all $\tau \in T_{\infty}(A)$ and all $b \in \mathrm{C}^{*}(A,X)_{+}$ of norm at most one we estimate
\begin{eqnarray*}
 \tau((\be_{A_{\infty}} - d^{(1)}_{\eta})b) & \stackrel{\eqref{w3-3}}{=}  &   \tau \left(g_{0,\eta} \left(\sum_{j=1}^{(m+1)^{2}}\psi^{(j)}(e_{11} \oplus e_{11}) \right) b \right) \\
& \le &   \lim_{l \to \infty}   \tau\left( \left(\sum_{j=1}^{(m+1)^{2}} \psi^{(j)}(e_{11} \oplus e_{11}) \right)^{1/l} b\right) \\
& \le & \sum_{j=1}^{(m+1)^{2}} \lim_{l \to \infty}   \tau((\psi^{(j)}(e_{11} \oplus e_{11}) )^{1/l} b) \\
& \stackrel{\ref{order-zero-notation}}{=} & \sum_{j=1}^{(m+1)^{2}} \lim_{l \to \infty}   \tau((\psi^{(j)})^{1/l}(e_{11} \oplus e_{11})  b)  \\
& \le & \sum_{j=1}^{(m+1)^{2}} \lim_{l \to \infty}  \frac{1}{\bar{k}} \cdot  \tau((\psi^{(j)})^{1/l}(\be_{\bar{k}} \oplus \be_{\bar{k}+1})  b) \\
& \stackrel{\eqref{w3-4}}{\le} & \frac{(m+1)^{2}}{2(m+1)^{2}} \cdot \tau(b)  \\
& = & \frac{1}{2} \cdot \tau(b),
\end{eqnarray*}
whence 
\begin{equation}
\label{w3-6}
\tau(d_{\eta}^{(1)}b) \ge \left(1- \frac{1}{2} \right) \cdot \tau(b) = \frac{1}{2} \cdot \tau(b).
\end{equation}
Here, for the second estimate we have used that the map 
\[
a \mapsto \lim_{l \to \infty} \tau(a^{1/l} b)
\]
is a dimension function on 
\[
\mathrm{C}^{*}(\psi^{(j)}(e_{11}\oplus e_{11})\mid j=1, \ldots,(m+1)^{2}) \subset A_{\infty} \cap A' \cap X',
\]
and as such respects the order structure of the Cuntz semigroup, cf.\ \cite{BlaHan:quasitrace}; for the third inequality we have used that 
\[
\tau((\psi^{(j)})^{1/l}(\, .\,)b):M_{\bar{k}} \oplus M_{\bar{k}+1} \to \mathbb{C}
\]
is a  trace by \ref{order-zero-traces} (cf.\ the argument for \eqref{w3-1}).

Taking $\frac{1}{l+1}$ for $\eta$, we then obtain sequences 
\begin{equation}
\label{w3-9}
\left(d^{(0)}_{\frac{1}{l+1}}\right)_{l \in \mathbb{N}}, \left(d^{(1)}_{\frac{1}{l+1}}\right)_{l \in \mathbb{N}} \subset A_{\infty} \cap A' \cap X'
\end{equation}
of positive contractions satisfying
\begin{equation}
\label{w3-8}
d^{(0)}_{\frac{1}{l+1}} \perp d^{(1)}_{\frac{1}{l+1}}
\end{equation}
and (by \eqref{w3-5} and  \eqref{w3-6})
\begin{equation}
\label{2-d-i-1}
\tau\left(d^{(i)}_{\frac{1}{l+1}} b\right) \ge \frac{1}{\bar{k}+1} \cdot \tau(b) - \frac{1}{l+1}, \,  i=0,1,
\end{equation}
for each $l \in \mathbb{N}$, $\tau \in T_{\infty}(A)$ and $b \in \mathrm{C}^{*}(A,X)_{+}$ of norm at most one.

For each $l \in \mathbb{N}$ and $i=0,1$, let 
\[
(\bar{d}^{(i)}_{l,r})_{r \in \mathbb{N}} \in \prod_{\mathbb{N}} A
\]
be a sequence of positive contractions lifting $d^{(i)}_{\frac{1}{l+1}}$. Let 
\[
(\mathcal{G}_{l})_{l \in \mathbb{N}}
\]
be a nested sequence of finite subsets in the unit ball of $(\prod_{\mathbb{N}} A)_{+}$ such that the image of $\bigcup_{l\in \mathbb{N}} \mathcal{G}_{l} \subset \prod A$ in $A_{\infty}$ is a dense subset of the unit ball of  $\mathrm{C}^{*}(A,X)_{+}$.

We claim that for each $l \in \mathbb{N}$ there is $K_{l} \in \mathbb{N}$ such that 
\begin{equation}
\label{w3-10}
\tau(\bar{d}^{(i)}_{l,r} f_{r}) \ge \frac{1}{\bar{k}+1} \cdot \tau(f_{r}) - \frac{2}{l+1}
\end{equation}
for $i=0,1$, for each $\tau \in T(A)$, each $(f_{r'})_{r'\in \mathbb{N}} \in \mathcal{G}_{l}$ and for each $r \ge K_{l}$.

Indeed, if this was not true, then for some $l \in \mathbb{N}$ there were $i \in \{0,1\}$, $(f_{r})_{r \in \mathbb{N}} \in \mathcal{G}_{l}$, an increasing sequence $(r_{s})_{s\in \mathbb{N}} \subset \mathbb{N}$ and a sequence $(\tau_{s})_{s\in \mathbb{N}} \subset T(A)$ with 
\begin{equation}
\label{w3-7}
\tau_{s}(\bar{d}^{(i)}_{l,r_{s}} f_{r_{s}}) < \frac{1}{\bar{k}+1} \cdot \tau_{s}(f_{r_{s}}) - \frac{2}{l+1}.
\end{equation}
Let $\omega \in \beta\mathbb{N} \setminus \mathbb{N}$ be a free ultrafilter, then 
\[
[(a_{r})_{r \in \mathbb{N}}] \mapsto \lim_{\omega} \tau_{s} (a_{r_{s}})
\]
is a trace in $T_{\infty}(A)$. Taking the limit along $\omega$ on both sides of \eqref{w3-7}, we obtain a contradiction to \eqref{2-d-i-1}, hence the claim holds.

On increasing the $K_{l}$, if necessary, we may also assume that they form a strictly increasing sequence, that
\begin{equation}
\label{w3-11}
\|\bar{d}^{(0)}_{l,r} \bar{d}^{(1)}_{l,r} \| < \frac{1}{l+1}
\end{equation}
and that
\begin{equation}
\label{w3-12}
\|[\bar{d}^{(i)}_{l,r},f_{r}]\| < \frac{1}{l+1},\,  i=0,1,
\end{equation}
for any $r \ge K_{l}$ and any $(f_{r'})_{r' \in \mathbb{N}} \in \mathcal{G}_{l}$ (using \eqref{w3-8} and \eqref{w3-9}).

We now define  sequences 
\[
(\tilde{d}^{(i)}_{r})_{r \in \mathbb{N}} \in \prod_{\mathbb{N}} A, \,i=0,1,
\]
of positive contractions by
\begin{equation}
\label{w3-13}
\tilde{d}^{(i)}_{r} :=
\begin{cases}
\bar{d}^{(i)}_{l,K_{l}} & \mbox{ if } K_{l} \le r < K_{l+1} \\
0 & \mbox{ if } r<K_{0}.
\end{cases}
\end{equation}
By combining \eqref{w3-13} with \eqref{w3-11}, \eqref{w3-12}, \eqref{w3-10} and the facts that the  $\mathcal{G}_{l}$ are nested and the $K_{l}$ are strictly increasing, it is straightforward to check the following properties:
\begin{itemize}
\item[(a)] $\|\tilde{d}_{r}^{(0)} \tilde{d}^{(1)}_{r}\| \stackrel{r \to \infty}{\longrightarrow} 0$
\item[(b)] $\|[\tilde{d}_{r}^{(i)},f_{r}]\| \stackrel{r \to \infty}{\longrightarrow} 0$, $i=0,1$, for every $(f_{r})_{r \in \mathbb{N}} \in \bigcup_{l \in \mathbb{N}} \mathcal{G}_{l}$
\item[(c)] for any $\delta>0$ and any $(f_{r})_{r \in \mathbb{N}} \in \bigcup_{l \in \mathbb{N}} \mathcal{G}_{l}$ there is $\bar{r} \in \mathbb{N}$ such that 
\[
\tau(\tilde{d}^{(i)}_{r} f_{r}) \ge \frac{1}{\bar{k}+1} \cdot \tau(f_{r}) - \delta, \, i=0,1,
\]
for any $r \ge \bar{r}$ and any $\tau \in T(A)$.
\end{itemize} 
Then, 
\[
d^{(i)}:= [(\tilde{d}^{(i)}_{r})_{r \in \mathbb{N}} ] \in A_{\infty}, \, i=0,1,
\]
will have the desired properties: They are positive contractions because all the $\tilde{d}^{(i)}_{r}$ are; they are orthogonal because of (a) and they commute with $\mathrm{C}^{*}(A,X)$ because of (b). Finally, (c) shows that for any $\tau \in T_{\infty}(A)$ and any $b \in \mathrm{C}^{*}(A,X)_{+}$ we have 
\[
\tau(d^{(i)}b) \ge \frac{1}{\bar{k}+1} \cdot \tau(b) \stackrel{\eqref{ww3-1}}{=} \frac{1}{L_{m,2}} \cdot \tau(b), \, i=0,1.
\]
This proves the proposition for $k = 2$. 

Let us now use induction to handle the cases when $k=2^{l}$ for  $1 \le l \in \mathbb{N}$. The anchor step $l=1$ has already been taken care of, so let us assume that the statement of the proposition holds for $k=2^{l}$ for some $l \ge 1$, with constant $L_{m,2^{l}}$ and positive elements 
\[
\hat{d}^{(1)},\ldots, \hat{d}^{(2^{l})} \subset A_{\infty} \cap A' \cap X'.
\]
Set
\[
L_{m,2^{l+1}}:= L_{m,2} \cdot L_{m,2^{l}}
\] 
and apply the first part of the proof (with 
\[
\tilde{X}:= X \cup \{\hat{d}^{(1)}, \ldots, \hat{d}^{(2^{l})}\}
\]
in place of $X$) to obtain positive elements 
\[
\tilde{d}^{(0)}, \tilde{d}^{(1)} \in A_{\infty} \cap A' \cap \tilde{X}'.
\]
Set
\[
d^{(j)}:=
\begin{cases}
\tilde{d}^{(0)} \hat{d}^{(j)} & \mbox{ if } 1 \le j \le 2^{l} \\
\tilde{d}^{(1)} \hat{d}^{(j-2^{l})} & \mbox{ if } 2^{l} < j \le 2^{l+1}; 
\end{cases}
\]
it is straightforward to check that these indeed satisfy the assertion of the proposition. 

This settles the case $k=2^{l}$ for $1 \le l \in \mathbb{N}$.  The case of general $k \in \mathbb{N}$ is an immediate consequence.
\end{nproof}
\en

\bn
\label{dimnuc-tracial-divisibility}
\begin{nprop}
Given $m \in \mathbb{N}$, there is $\tilde{m} \in \mathbb{N}$ such that the following holds:

If $A$ is a separable, simple, nonelementary and unital $\mathrm{C}^{*}$-algebra with $\dimnuc A \le m$, then $A$ has tracial $\tilde{m}$-almost divisibility. 
\end{nprop}

\begin{nproof}
Invoke Proposition~\ref{2-d-i} to define
\begin{equation}
\label{w3-20}
\tilde{m}:= L_{m,(m+1)^{2}} - 1.
\end{equation}
Let $b \in A$  a positive contraction, $k \in \mathbb{N}$ and $\epsilon >0$ be given. (Since any matrix algebra over $A$ will satisfy the hypotheses of the proposition just as $A$ does, to test  tracial $\tilde{m}$-almost divisibility, it will suffice to consider a positive contraction $b$ in $A$, as opposed to $M_{\infty}(A)$.) 

Choose $\bar{k} \in \mathbb{N}$ such that 
\begin{equation}
\label{w3-18}
\frac{1}{\bar{k}+1} < \frac{\epsilon}{4}
\end{equation}
and such that $k$ divides $\bar{k}$. Use Proposition~\ref{Kirchberg-cov} to find c.p.c.\ order zero maps 
\[
\varphi^{(1)}, \ldots, \varphi^{((m+1)^{2})}:M_{\bar{k}} \oplus M_{\bar{k}+1} \to A_{\infty} \cap A'
\]
such that 
\begin{equation}
\label{w3-16}
\sum_{r=1}^{(m+1)^{2}} \varphi^{(r)}(\be_{\bar{k}} \oplus \be_{\bar{k}+1}) \ge \be_{A_{\infty}}.
\end{equation}
Take pairwise orthogonal positive contractions 
\[
d^{(1)}, \ldots, d^{((m+1)^{2})} \in A_{\infty} \cap A' \cap \left(\bigcup_{r=1}^{(m+1)^{2}} \varphi^{(r)} (M_{\bar{k}} \oplus M_{\bar{k}+1})\right)'
\]
as in Proposition~\ref{2-d-i} (with $(m+1)^{2}$ in place of $k$ and $\bigcup_{r=1}^{(m+1)^{2}} \varphi^{(r)} (M_{\bar{k}} \oplus M_{\bar{k}+1})$ in place of $X$).

We then obtain a c.p.\ map 
\[
\bar{\varphi}: M_{\bar{k}} \oplus M_{\bar{k}+1} \to A_{\infty} \cap A'
\]
by setting
\begin{equation}
\label{w3-14}
\bar{\varphi}(x):=  \sum_{r=1}^{(m+1)^{2}} \varphi^{(r)}(x) d^{(r)} ;
\end{equation}
since the $d^{(j)}$ are pairwise orthogonal and commute with the images of the c.p.c.\ order zero maps $\varphi^{(j)}$, the map $\bar{\varphi}$ is itself c.p.c.\ order zero.  Note that we have
\begin{eqnarray}
\tau(b \bar{\varphi}(e_{jj} \oplus e_{jj})) & = & \tau(b \bar{\varphi}(e_{11} \oplus e_{11})) \nonumber \\ 
& \stackrel{\eqref{w3-14}}{=} &  \tau \left( \sum_{r=1}^{(m+1)^{2}} b \varphi^{(r)}(e_{11} \oplus e_{11}) d^{(r)}  \right) \nonumber \\
& \stackrel{\eqref{w3-15}}{\ge} &    \frac{1}{L_{m,(m+1)^{2}}} \cdot \tau\left( \sum_{r=1}^{(m+1)^{2}} b \varphi^{(r)}(e_{11} \oplus e_{11})  \right) \nonumber \\
& \ge &      \frac{1}{L_{m,(m+1)^{2}} (\bar{k}+1)} \cdot \tau\left(\sum_{r=1}^{(m+1)^{2}} b \varphi^{(r)}(\be_{\bar{k}} \oplus \be_{\bar{k}+1})  \right) \nonumber \\
& \stackrel{\eqref{w3-16}}{\ge} &      \frac{1}{L_{m,(m+1)^{2}} (\bar{k}+1)} \cdot \tau(b  ) \label{w3-17}
\end{eqnarray}
for $j=1,\ldots,\bar{k}$ and for all $\tau \in T_{\infty}(A)$ (for the first equality and for the second inequality we use Remark~\ref{order-zero-traces}, cf.\ the  argument for \eqref{w3-1}). From this we obtain
\begin{eqnarray}
\lefteqn{\tau\left( f_{\epsilon/4,\epsilon/2}(b) \bar{\varphi}\left( \sum_{l=1}^{\bar{k}} e_{ll} \oplus e_{ll}  \right)  \right) }\nonumber  \\
&  \ge&  \tau\left(b \bar{\varphi}\left( \sum_{l=1}^{\bar{k}} e_{ll} \oplus e_{ll}  \right)  \right)  - \frac{\epsilon}{4} \nonumber\\
& \stackrel{\eqref{w3-17}}{\ge} & \frac{\bar{k}}{L_{m,(m+1)^{2}} \cdot (\bar{k}+1)} \cdot \tau(b) - \frac{\epsilon}{4} \nonumber \\
& \ge &  \frac{1}{L_{m,(m+1)^{2}}} \cdot \tau(b) - \frac{1}{\bar{k}+1} -\frac{\epsilon}{4} \nonumber \\
& \stackrel{\eqref{w3-18}}{\ge} & \frac{1}{L_{m,(m+1)^{2}}} \cdot \tau(b) - \frac{\epsilon}{2}
\label{dimnuc-tracial-divisibility1}
\end{eqnarray}
for all $\tau \in T_{\infty}(A)$.

Let 
\begin{equation}
\label{w3-7-1}
\tilde{\varphi} = (\tilde{\varphi}_{l})_{l \in \mathbb{N}}: M_{\bar{k}} \oplus M_{\bar{k}+1} \to \overline{b \left(\prod A \right) b} \subset \prod_{\mathbb{N}}A
\end{equation}
be a c.p.c.\ order zero lift (cf.\ \cite[Remark~2.4]{KirWinter:dr}) of the c.p.c.\ order zero map 
\[
x \mapsto f_{\epsilon/4,\epsilon/2}(b) \bar{\varphi}(x).
\] 

We claim that there is $\bar{l} \in \mathbb{N}$ such that 
\begin{equation}
\label{w3-18b}
\tau\left( \tilde{\varphi}_{\bar{l}} \left( \sum_{j=1}^{\bar{k}} e_{jj} \oplus e_{jj}   \right) \right) \ge \frac{1}{L_{m,(m+1)^{2}}} \cdot \tau(b) - \epsilon
\end{equation}
for all $\tau \in T(A)$. 

If not, then for every $l \in \mathbb{N}$ there is some $\tau_{l} \in T(A)$ such that
\[
\tau_{l}\left( \tilde{\varphi}_{l} \left( \sum_{j=1}^{\bar{k}} e_{jj} \oplus e_{jj}   \right) \right) < \frac{1}{L_{m,(m+1)^{2}}} \cdot \tau_{l}(b) - \epsilon;
\]
let $\omega \in \beta\mathbb{N} \setminus \mathbb{N}$ be a free ultrafilter, then 
\[
\lim_{\omega} \tau_{l}\left( \tilde{\varphi}_{l} \left( \sum_{j=1}^{\bar{k}} e_{jj} \oplus e_{jj}   \right) \right) \le \frac{1}{L_{m,(m+1)^{2}}} \cdot \lim_{\omega} \tau_{l}(b) - \epsilon,
\]
a contradiction to \eqref{dimnuc-tracial-divisibility1}, thus proving the claim.

Let 
\[
\iota_{1}:M_{k} \to M_{\bar{k}}
\]
be a unital embedding ($k$ divides $\bar{k}$) and let 
\[
\iota_{2}:M_{\bar{k}} \to M_{\bar{k}+1}
\]
be the upper left corner embedding; set
\[
\iota:= (\id_{M_{\bar{k}}} \oplus \iota_{2}) \circ \iota_{1}: M_{k} \to M_{\bar{k}} \oplus M_{\bar{k}+1}
\]
and 
\begin{equation}
\label{w3-19}
\psi:= \tilde{\varphi}_{\bar{l}} \circ \iota,
\end{equation}
then
\[
\psi:M_{k} \to \her(b) \subset A
\]
is c.p.c.\ order zero (cf.\ \eqref{w3-7-1}) and
\[
\tau(\psi(\be_{k})) \stackrel{\eqref{w3-19},\eqref{w3-18b}}{\ge} \frac{1}{L_{m,(m+1)^{2}}} \cdot \tau(b) - \epsilon \stackrel{\eqref{w3-20}}{=} \frac{1}{\tilde{m}+1} \cdot \tau(b) - \epsilon
\]
for all $\tau \in T(A)$.
\end{nproof}
\en

\bn
\label{dimnuc-strong-tracial-comparison}
\begin{nprop}
Given $m \in \mathbb{N}$, there is $\bar{m} \in \mathbb{N}$ such that the following holds:

If $A$ is a separable, simple, unital $\mathrm{C}^{*}$-algebra with $\dimnuc A \le m$, then $A$ has strong tracial $\bar{m}$-comparison.
\end{nprop}

\begin{nproof}
Take 
\begin{equation}
\label{w3-25}
\bar{m}:= L_{m,(m+1)^{2}} \cdot (m+2) -1,
\end{equation}
where $ L_{m,(m+1)^{2}}$ comes from Proposition~\ref{2-d-i}.

Since any matrix algebra over $A$ will satisfy the hypotheses of the proposition just as $A$ does, to test strong tracial $\bar{m}$-comparison it will suffice to consider positive contractions $a$ and $b$ in $A$ (as opposed to $M_{\infty}(A)$). In view of Proposition~\ref{tau-dtau-eta} it will suffice to show that, if  
\begin{equation}
\label{dimnuc-strong-tracial-comparison1}
d_{\tau}(a) < \frac{1}{\bar{m}+1} \cdot \tau(b) - \eta
\end{equation}
for some $\eta>0$ and for all $\tau \in T(A)$, then
\[
f_{\epsilon,2\epsilon}(a) \precsim b
\]
for any given $\epsilon >0$. (We may use traces as opposed to quasitraces  by Remark~\ref{quasitraces-traces}.)

So, let $a,b, \eta$ and $\epsilon$ as above be given. We then have
\begin{equation}
\label{w3-26}
\tau(g_{\epsilon/2,\epsilon}(a)) \le d_{\tau}(a) \stackrel{\eqref{dimnuc-strong-tracial-comparison1}}{<} \frac{1}{\bar{m}+1} \cdot \tau(b) - \eta \le  \frac{1}{\bar{m}+1} \cdot \tau(g_{\eta/2,\eta}(b)) - \frac{\eta}{2}
\end{equation}
for all $\tau \in T(A)$. Use Proposition~\ref{Kirchberg-cov} to find c.p.c.\ order zero maps
\begin{equation}
\label{w3-21}
\varphi^{(1)},\ldots,\varphi^{((m+1)^{2})}:M_{m+1} \oplus M_{m+2} \to A_{\infty} \cap A'
\end{equation}
such that 
\begin{equation}
\label{w3-24}
\sum_{j=1}^{(m+1)^{2}} \varphi^{(j)}(\be_{m+1} \oplus \be_{m+2}) \ge \be_{A_{\infty}}.
\end{equation}
Take pairwise orthogonal positive contractions 
\begin{equation}
\label{w3-22}
d^{(1)}, \ldots, d^{((m+1)^{2})} \in A_{\infty} \cap A' \cap \left(  \bigcup_{j=1}^{(m+1)^{2}} \varphi^{(j)}(M_{m+1} \oplus M_{m+2})  \right)'
\end{equation}
as in Proposition~\ref{2-d-i} (with $(m+1)^{2}$ in place of $k$). We then have
\begin{equation}
\label{w3-23}
\tau(c \varphi^{(j)}(e_{11} \oplus e_{11}) d^{(j)}) \stackrel{\eqref{w3-15}}{\ge} \frac{1}{L_{m,(m+1)^{2}}} \cdot \tau(c \varphi^{(j)} (e_{11} \oplus e_{11}))
\end{equation}
for all $j = 1, \ldots,(m+1)^{2}$, $c \in A_{+}$ and $\tau \in T_{\infty}(A)$.  Let 
\[
\bar{\varphi}:M_{m+1} \oplus M_{m+2} \to \overline{b A_{\infty}b}
\]
be the c.p.\ map given by 
\begin{equation}
\label{w3-8-1}
\bar{\varphi}(\, .\,) := \sum_{j=1}^{(m+1)^{2}} g_{\eta/2,\eta}(b) \varphi^{(j)}(\, .\,) d^{(j)};
\end{equation}
using \eqref{w3-21} and \eqref{w3-22} and mutual orthogonality of the $d^{(j)}$ one checks that $\bar{\varphi}$ in fact is a c.p.c.\ order zero map. We compute
\begin{eqnarray}
\lefteqn{\tau(\bar{\varphi}(e_{11} \oplus e_{11})) } \nonumber \\
& \stackrel{\eqref{w3-8-1}}{=} & \tau \left( \sum_{j=1}^{(m+1)^{2}} g_{\eta/2,\eta}(b) \varphi^{(j)}(e_{11} \oplus e_{11}) d^{(j)} \right) \nonumber \\
& \stackrel{\eqref{w3-23}}{\ge} & \frac{1}{L_{m,(m+1)^{2}}} \cdot \tau \left( \sum_{j=1}^{(m+1)^{2}} g_{\eta/2,\eta}(b) \varphi^{(j)}(e_{11} \oplus e_{11})  \right) \nonumber \\
& \ge & \frac{1}{L_{m,(m+1)^{2}} \cdot (m+2)} \cdot \tau \left( \sum_{j=1}^{(m+1)^{2}} g_{\eta/2,\eta}(b) \varphi^{(j)}(\be_{m+1} \oplus \be_{m+2})  \right) \nonumber \\
& \stackrel{\eqref{w3-24}}{\ge} & \frac{1}{L_{m,(m+1)^{2}} \cdot (m+2)} \cdot \tau (  g_{\eta/2,\eta}(b) ) \nonumber \\
& \stackrel{\eqref{w3-25}}{=} & \frac{1}{\bar{m}+1} \cdot \tau (  g_{\eta/2,\eta}(b) ) \nonumber \\
& \stackrel{\eqref{w3-26}}{>} & \tau(g_{\epsilon/2,\epsilon}(a)) +\frac{\eta}{2} \label{w3-28}
\end{eqnarray}
for all $\tau \in T_{\infty}(A)$. For the second inequality we have once more used \ref{order-zero-traces} as in the argument for \eqref{w3-1}; for the last inequality we restrict $\tau$ to $A$, so that \eqref{w3-26} indeed applies. 

Take a c.p.c.\ order zero lift  (cf.\ \cite[Remark~2.4]{KirWinter:dr})
\begin{equation}
\label{w3-8-2}
\Phi = (\Phi_{l})_{l \in \mathbb{N}}: M_{m+1} \oplus M_{m+2}  \to \ \overline{b\left(\prod A\right) b} \subset\prod_{\mathbb{N}} A
\end{equation}
of the c.p.c.\ order zero map $\bar{\varphi}$. 

We claim that there is $\bar{l} \in \mathbb{N}$ such that
\begin{equation}
\label{w3-27}
\tau(\Phi_{\bar{l}} (e_{11} \oplus e_{11}) ) \ge \tau(g_{\epsilon/2,\epsilon}(a)) + \frac{\eta}{4}
\end{equation}
for all $\tau \in T(A)$. If not, for every $l \in \mathbb{N}$ there is $\tau_{l} \in T(A)$ such that
\[
\tau_{l}(\Phi_{l}(e_{11} \oplus e_{11})) < \tau_{l}(g_{\epsilon/2,\epsilon}(a)) +\frac{\eta}{4};
\]
for a free ultrafilter $\omega \in \beta\mathbb{N} \setminus \mathbb{N}$ we obtain
\[
\lim_{\omega} \tau_{l}(\Phi_{l}(e_{11} \oplus e_{11})) \le \lim_{\omega} \tau_{l}(g_{\epsilon/2,\epsilon}(a)) + \frac{\eta}{4},
\]
a contradiction to \eqref{w3-28}.

We obtain
\[
d_{\tau}(\Phi_{\bar{l}}(e_{11} \oplus e_{11})) \ge \tau(\Phi_{\bar{l}}(e_{11} \oplus e_{11}))   \stackrel{\eqref{w3-27}}{\ge} \tau(g_{\epsilon/2,\epsilon}(a)) +\frac{\eta}{4} > d_{\tau}(f_{\epsilon,2\epsilon}(a))
\]
for all $\tau \in T(A)$. Now by $m$-comparison (see Theorem~\ref{dimnuc-comparison}) and since $\Phi_{\bar{l}}$ is c.p.c.\ of order zero, we have
\[
f_{\epsilon,2\epsilon}(a) \precsim \sum_{j=1}^{m+1} \Phi_{\bar{l}}(e_{jj} \oplus e_{jj}) \stackrel{\eqref{w3-8-2}}{\precsim} b,
\]
as desired.
\end{nproof}
\en

\bn
\label{dimnuc-tracially-pure}
We note the following direct consequence of Propositions~\ref{dimnuc-tracial-divisibility} and \ref{dimnuc-strong-tracial-comparison} explicitly.

\begin{ncor}
If $A$ is a separable, simple, nonelementary and unital $\mathrm{C}^{*}$-algebra with finite nuclear dimension, then $A$ is tracially $(m,\bar{m})$-pure for suitable $m,\bar{m} \in \mathbb{N}$.
\end{ncor}
\en

\section{Tracial matrix cone absorption}
\label{tracial-matrix-cone-absorption}

\noindent
We now set out to prove a refined version of \cite[Lemma~3.10]{Winter:dr-Z-stable}, Lemma~\ref{3-10prime} below. The strategy is to isolate some of the abstract machinery in \cite[Section~3]{Winter:dr-Z-stable} so that we do not have to ask for finite decomposition rank (or finite nuclear dimension) beforehand.

\bn
\label{Hnew}
\begin{nprop}
Let $A$ be a separable, simple, unital $\mathrm{C}^{*}$-algebra with  tracial $m$-almost divisibility, let $0 \neq k \in \mathbb{N}$ and $d \in A_{\infty}$ be a positive contraction.

Then, for any $\eta>0$ there is a c.p.c.\ order zero map 
\[
\varphi:M_{k} \to \overline{dA_{\infty}d}
\]
such that
\[
\tau(\varphi(\be_{k})) \ge \frac{1}{m+1} \cdot \tau(d) - \eta
\]
for all $\tau \in T_{\infty}(A)$.
\end{nprop}

\begin{nproof}
Let 
\[
\bar{d}= (\bar{d}_{n})_{n \in \mathbb{N}} \in \prod_{\mathbb{N}} A
\]
be a positive contraction lifting $d$. Set
\[
\tilde{d}:= g_{\eta/2,\eta}(\bar{d}) = (g_{\eta/2,\eta}(\bar{d}_{n}))_{n \in \mathbb{N}}.
\]
For each $n \in \mathbb{N}$, by  tracial $m$-almost divisibility there is a c.p.c.\ order zero map
\begin{equation}
\label{ww4-1}
\varphi_{n}:M_{k} \to \her(g_{\eta/2,\eta}(\bar{d}_{n})) \subset A
\end{equation}
such that 
\begin{eqnarray}
\tau(\varphi_{n}(\be_{k})) & \ge & \frac{1}{m+1} \cdot \tau(g_{\eta/2,\eta}(\bar{d}_{n})) - \frac{\eta}{2} \nonumber \\
& \ge &  \frac{1}{m+1} \cdot \tau(\bar{d}_{n}) - \eta \label{ww4-2}
\end{eqnarray}
for all $\tau \in T(A)$. Let
\[
\bar{\varphi}:M_{k} \to \prod_{\mathbb{N}} A
\]
and
\[
\varphi:M_{k} \to A_{\infty}
\]
be the induced c.p.c.\ maps; it is clear (from \eqref{ww4-1}) that 
\[
g_{0,\eta/2}(\bar{d}) \bar{\varphi}(\be_{k}) = \bar{\varphi}(\be_{k}),
\]
whence
\[
g_{0,\eta/2}(d) \varphi(\be_{k}) = \varphi(\be_{k})
\]
and 
\[
\varphi(M_{k}) \subset \overline{d A_{\infty} d} \subset A_{\infty}.
\]
Moreover, by \eqref{ww4-2} we have
\[
\tau(\varphi(\be_{k})) \ge  \frac{1}{m+1} \cdot \tau(d) - \eta
\]
for all $\tau \in T_{\infty}(A)$.
\end{nproof}
\en

\bn
\label{Gnew}
\begin{nprop}
Let $A$ be a separable, simple, unital $\mathrm{C}^{*}$-algebra with tracial $m$-almost divisibility, let $0 \neq k \in \mathbb{N}$ and $d \in A_{\infty}$ be a positive contraction.

Then, there are c.p.c.\ order zero maps 
\[
\psi^{(0)},\psi^{(1)}:M_{k} \to A_{\infty} \cap \{d\}'
\]
satisfying
\[
\tau((\psi^{(0)}(\be_{k}) + \psi^{(1)}(\be_{k})  ) f(d)) \ge  \frac{1}{m+1} \cdot \tau(f(d))
\]
for all $\tau \in T_{\infty}(A)$ and all $f \in \mathcal{C}_{0}((0,1])_{+}$.
\end{nprop}

\begin{nproof}
Let 
\[
\bar{d} = (\bar{d}_{q})_{q \in \mathbb{N}} \in \prod_{\mathbb{N}} A
\]
be a positive contraction lifting $d$. Let $(\mathcal{G}_{q})_{q \in \mathbb{N}}$ be a nested sequence of finite subsets of the unit ball of $\mathcal{C}_{0}((0,1])_{+}$ with dense union. For each $q \in \mathbb{N}$, it is straightforward to choose $0 \neq L_{q} \in \mathbb{N}$ such that the sawtooth functions $f^{(i)}_{q,j} \in \mathcal{C}_{0}((0,1])_{+}$, given by
\[
f^{(1)}_{q,j}(t) = \left\{
\begin{array}{ll}
0 & \mbox{if } 0 \le t \le \frac{j-1}{L_{q}} \\
1 & \mbox{if } t = \frac{j-1/2}{L_{q}} \\
0 & \mbox{if } \frac{j}{L_{q}} \le t \le 1\\
\mbox{linear} & \mbox{else}
\end{array}
\right.
\]
and 
\[
f^{(0)}_{q,j}(t) = \left\{
\begin{array}{ll}
0 & \mbox{if } 0 \le t \le \frac{1}{2L_{q}}\\
f^{(1)}_{q,j}(t- \frac{1}{2L_{q}}) & \mbox{else}
\end{array}
\right.
\]
for $j=1,\ldots,L_{q}$, satisfy 
\begin{equation}
\label{Gnew1}
\left\|f \cdot g(f_{q,j}^{(i)}) - f(j/L_{q}) \cdot g(f^{(i)}_{q,j}) \right\| \le \frac{1}{q+1}
\end{equation}
for $j=1, \ldots,L_{q}$, $i=0,1$, $f \in \mathcal{G}_{q}$ and $g$ in the unit ball of $\mathcal{C}_{0}((0,1])_{+}$, and
\begin{equation}
\label{Gnew5}
\left\| \sum_{i,j} f(j/L_{q}) \cdot f^{(i)}_{q,j} - f \right\| \le \frac{1}{q+1} 
\end{equation}
for $f \in \mathcal{G}_{q}$.

By tracial $m$-almost divisibility, for each $q \in \mathbb{N}$, $i \in \{0,1\}$, $ j \in \{1,\ldots,L_{q}\}$ there exists a c.p.c.\ order zero map 
\begin{equation}
\label{Gnew2}
\psi^{(i)}_{q,j}:M_{k} \to \her(f^{(i)}_{q,j}(\bar{d}_{q})) \subset A
\end{equation}
such that 
\begin{equation}
\label{Gnew6}
\tau(\psi^{(i)}_{q,j}(\be_{k})) \ge \frac{1}{m+1} \cdot \tau(f^{(i)}_{q,j}(\bar{d}_{q})) - \frac{1}{L_{q} \cdot (q+1)}
\end{equation}
for all $\tau \in T(A)$. Since for fixed $q \in \mathbb{N}$ and $i \in \{0,1\}$ the functions $f^{(i)}_{q,j}, j=1,\ldots,L_{q}$, are pairwise orthogonal, the maps
\[
\psi^{(i)}_{q}:= \bigoplus_{j=1}^{L_{q}} \psi^{(i)}_{q,j} :M_{k} \to A
\]
are c.p.c.\ order zero (by \eqref{Gnew2}). 

Using \eqref{Gnew1} and \eqref{Gnew2} it is straightforward to estimate 
\begin{equation}
\label{Gnew7}
\left\| f(\bar{d}_{q}) \psi^{(i)}_{q}(x) - \sum_{j=1}^{L_{q}} f(j/L_{q}) \cdot \psi^{(i)}_{q}(x) \right\| \le \frac{1}{q+1}
\end{equation}
for $q \in \mathbb{N}$, $i \in \{0,1\}$, $x$ in the unit ball of $(M_{k})_{+}$ and $f \in \mathcal{G}_{q}$, whence 
\begin{equation}
\label{Gnew3}
\| [f(\bar{d}_{q}), \psi^{(i)}_{q}(x)]\| \le \frac{2}{q+1}.
\end{equation}
Moreover, we have 
\begin{eqnarray}
\lefteqn{\tau \left( \sum_{i=0,1} f(\bar{d}_{q}) \psi^{(i)}_{q} (\be_{k})\right)} \nonumber \\
 & \stackrel{\eqref{Gnew7}}{\ge} & \tau \left( \sum_{i=0,1;j=1,\ldots,L_{q}} f(j/L_{q}) \cdot \psi^{(i)}_{q} (\be_{k}) \right) - \frac{2}{q+1} \nonumber \\
& \stackrel{\eqref{Gnew6}}{\ge} &  \frac{1}{m+1} \cdot  \tau \left( \sum_{i=0,1;j=1,\ldots,L_{q}} f(j/L_{q}) \cdot f^{(i)}_{q,j} (\bar{d}_{q}) \right) - \frac{4}{q+1}  \nonumber \\ 
& \stackrel{\eqref{Gnew5}}{\ge} &  \frac{1}{m+1} \cdot \tau(f(\bar{d}_{q})) - \frac{5}{q+1}
\label{Gnew4}
\end{eqnarray}
for all $q \in \mathbb{N}$, $f \in \mathcal{G}_{q}$ and $\tau \in T(A)$. 

For $i=0,1$ let 
\[
\psi^{(i)}:M_{k} \to A_{\infty}
\]
be the maps induced by the $\psi^{(i)}_{q}$; it follows directly from \eqref{Gnew3} and \eqref{Gnew4}  that these have the desired properties.
\end{nproof}
\en

\bn
\label{Fnew}
The following will be a crucial ingredient for Propositions~\ref{Anew} and \ref{Bnew}, which in turn are  variations of Proposition~\ref{2-d-i}.

\begin{nprop}
Let $A$ be a separable, simple, unital  $\mathrm{C}^{*}$-algebra with  tracial $m$-almost divisibility,  and let $d \in A_{\infty}$ be a positive contraction.

Then, there are orthogonal positive contractions
\[
d_{0},d_{1} \in A_{\infty} \cap \{d\}'
\]
satisfying
\begin{equation}
\label{w4-12}
\tau(d_{i} f(d)) \ge \frac{1}{4(m+1)} \cdot \tau(f(d)) 
\end{equation}
for all $\tau \in T_{\infty}(A)$, all $f \in \mathcal{C}_{0}((0,1])_{+}$ and $i=0,1$.
\end{nprop}

\begin{nproof}
By Proposition~\ref{Gnew} there are c.p.c.\ order zero maps 
\[
\psi_{1},\psi_{2}:M_{4} \to A_{\infty} \cap \{d\}'
\]
satisfying 
\[
\tau((\psi_{1}(e_{11}) + \psi_{2}(e_{11})) f(d)) \ge \frac{1}{4(m+1)} \cdot \tau(f(d))
\]
for all $\tau \in T_{\infty}(A)$ and $f \in \mathcal{C}_{0}((0,1])_{+}$ (note that $\tau((\psi_{1}+\psi_{2})(\, . \,) f(d))$ is a trace on $M_{4}$). 

For $\eta>0$, define positive contractions
\[
d^{(0)}_{\eta}:= g_{\eta,2\eta}(\psi_{1}(e_{11})+ \psi_{2}(e_{11}))
\]
and
\begin{equation}
\label{w4-3-1}
d^{(1)}_{\eta}:= \be_{A_{\infty}} - g_{0,\eta}(\psi_{1}(e_{11})+ \psi_{2}(e_{11})).
\end{equation}
We have
\[
d^{(0)}_{\eta} \perp d^{(1)}_{\eta},
\]
\[
d^{(0)}_{\eta}, d^{(1)}_{\eta} \in A_{\infty} \cap\{d\}'
\]
and
\begin{equation}
\label{w4-1}
\tau(d^{(0)}_{\eta} f(d)) \ge \frac{1}{4(m+1)} \cdot \tau(f(d)) - \eta
\end{equation}
for all $\tau \in T_{\infty}(A)$ and all $f \in \mathcal{C}_{0}((0,1])_{+}$ of norm at most one.  Moreover, for all $\tau \in T_{\infty}(A)$ and all $f \in \mathcal{C}_{0}((0,1])_{+}$ we estimate
\begin{eqnarray*}
\tau((\be_{A_{\infty}} - d^{(1)}_{\eta})f(d)) & \stackrel{\eqref{w4-3-1}}{=}  & \tau(g_{0,\eta}(\psi_{1}(e_{11}) + \psi_{2}(e_{11})) f(d)) \\
& \le & \lim_{l \to \infty}   \tau((\psi_{1}(e_{11}) + \psi_{2}(e_{11}))^{1/l} f(d)) \\
& \le & \lim_{l \to \infty}   \tau((\psi_{1}(e_{11}) )^{1/l} f(d)) +    \lim_{l \to \infty}   \tau((\psi_{2}(e_{11}))^{1/l} f(d)) \\
& \stackrel{\ref{order-zero-notation}}{=} & \lim_{l \to \infty}   \tau(\psi_{1}^{1/l}(e_{11})  f(d)) +    \lim_{l \to \infty}   \tau(\psi_{2}^{1/l}(e_{11}) f(d)) \\
& = & \lim_{l \to \infty}  \frac{1}{4} \cdot  \tau(\psi_{1}^{1/l}(\be_{4})  f(d)) \\
& &+    \lim_{l \to \infty}   \frac{1}{4} \cdot \tau(\psi_{2}^{1/l}(\be_{4}) f(d)) \\
& \le & \frac{1}{4} \cdot \tau(f(d)) + \frac{1}{4} \cdot \tau(f(d)) \\
& = & \frac{1}{2} \cdot \tau(f(d)),
\end{eqnarray*}
whence 
\begin{equation}
\label{w4-2}
\tau(d_{\eta}^{(1)}f(d)) \ge \frac{1}{2} \cdot \tau(f(d)).
\end{equation}
Here, for the second estimate we have used that the map $a \mapsto \lim_{l \to \infty} \tau(a^{1/l} f(d))$ is a dimension function on $\mathrm{C}^{*}(\psi_{1}(e_{11}),\psi_{2}(e_{11})) \subset A_{\infty} \cap \{d\}'$, and as such respects the order structure of the Cuntz semigroup; for the third equality we have used that $\psi_{i}^{1/l}(\, .\,)f(d):M_{k} \to A_{\infty}$ is a c.p.c.\ order zero map, and that the composition of an order zero map with a trace is again a trace (cf.\ \ref{order-zero-traces}).

Taking $\frac{1}{l+1}$ for $\eta$, we then obtain sequences 
\begin{equation}
\label{w4-3-3}
\left(d^{(0)}_{\frac{1}{l+1}}\right)_{l \in \mathbb{N}}, \left(d^{(1)}_{\frac{1}{l+1}}\right)_{l \in \mathbb{N}} \subset A_{\infty} \cap\{d\}'
\end{equation}
of positive contractions satisfying
\begin{equation}
\label{w4-3-2}
d^{(0)}_{\frac{1}{l+1}} \perp d^{(1)}_{\frac{1}{l+1}}
\end{equation}
and
\begin{equation}
\label{Fnew1}
\tau\left(d^{(i)}_{\frac{1}{l+1}} f(d)\right) \stackrel{\eqref{w4-1},\eqref{w4-2}}{\ge} \frac{1}{4(m+1)} \cdot \tau(f(d)) - \frac{1}{l+1}, \,  i=0,1,
\end{equation}
for each $l \in \mathbb{N}$, $\tau \in T_{\infty}(A)$ and $f \in \mathcal{C}_{0}((0,1])_{+}$ of norm at most one.

For each $l \in \mathbb{N}$ and $i=0,1$, let 
\[
(\bar{d}^{(i)}_{l,r})_{r \in \mathbb{N}} \in \prod_{\mathbb{N}} A
\]
be a sequence of positive contractions lifting $d^{(i)}_{\frac{1}{l+1}}$, and let 
\[
(d_{r})_{r \in \mathbb{N}} \in \prod_{\mathbb{N}} A
\]
be a sequence of positive contractions lifting $d$. Let $(\mathcal{G}_{l})_{l \in \mathbb{N}}$ be a nested sequence of finite subsets in the unit ball of $\mathcal{C}_{0}((0,1])_{+}$ with dense union.

We claim that for each $l \in \mathbb{N}$ there is $K_{l} \in \mathbb{N}$ such that 
\begin{equation}
\label{w4-16}
\tau(\bar{d}^{(i)}_{l,r} f(d_{r})) \ge \frac{1}{4(m+1)} \cdot \tau(f(d_{r})) - \frac{2}{l+1}
\end{equation}
for $i=0,1$, for each $\tau \in T(A)$, each $f \in \mathcal{G}_{l}$ and for each $r \ge K_{l}$.

Indeed, if this was not true, then for some $l \in \mathbb{N}$ there were $i \in \{0,1\}$, $f \in \mathcal{G}_{l}$, an increasing sequence $(r_{s})_{s\in \mathbb{N}} \subset \mathbb{N}$ and a sequence $(\tau_{s})_{s\in \mathbb{N}} \subset T(A)$ with 
\begin{equation}
\label{w4-3}
\tau_{s}(\bar{d}^{(i)}_{l,r_{s}} f(d_{r_{s}})) < \frac{1}{4(m+1)} \cdot \tau_{s}(f(d_{r_{s}})) - \frac{2}{l+1}.
\end{equation}
If $\omega \in \beta \mathbb{N} \setminus \mathbb{N}$ is a free ultrafilter, then 
\[
[(a_{r})_{r\in \mathbb{N}}] \mapsto \lim_{\omega} \tau_{s} (a_{r_{s}})
\]
is a trace in $T_{\infty}(A)$. Taking the limit along  $\omega $ on both sides of \eqref{w4-3}, this yields a contradiction to \eqref{Fnew1}, hence the claim holds.

On increasing the $K_{l}$, if necessary,  by \eqref{w4-3-3} and \eqref{w4-3-2} we may also assume that they form a strictly increasing sequence, that
\begin{equation}
\label{w4-4}
\|\bar{d}^{(0)}_{l,r} \bar{d}^{(1)}_{l,r} \| < \frac{1}{l+1}
\end{equation}
and that
\begin{equation}
\label{w4-5}
\|[\bar{d}^{(i)}_{l,r},d_{r}]\| < \frac{1}{l+1},\,  i=0,1,
\end{equation}
for any $r \ge K_{l}$.

We now define sequences 
\[
(\tilde{d}^{(i)}_{r})_{r \in \mathbb{N}} \in \prod_{\mathbb{N}} A, \,i=0,1, 
\]
of positive contractions by
\[
\tilde{d}^{(i)}_{r}:= 
\begin{cases}
d^{(i)}_{l,r} & \mbox{ if } K_{l} \le r < K_{l+1} \\
0 & \mbox{ if } r< K_{0}.
\end{cases}
\]
The following properties are straightforward to check, using \eqref{w4-4}, \eqref{w4-5} and \eqref{w4-16}:
\begin{itemize}
\item[(a)] $\|\tilde{d}_{r}^{(0)} \tilde{d}^{(1)}_{r}\| \stackrel{r \to \infty}{\longrightarrow} 0$
\item[(b)] $\|[\tilde{d}_{r}^{(i)},d_{r}]\| \stackrel{r \to \infty}{\longrightarrow} 0$, $i=0,1$
\item[(c)] for any $\delta>0$ and any $f \in \bigcup_{l \in \mathbb{N}} \mathcal{G}_{l}$ there is $\bar{r} \in \mathbb{N}$ such that 
\[
\tau(\tilde{d}^{(i)}_{r} f(d_{r})) \ge \frac{1}{4(m+1)} \cdot \tau(f(d_{r})) - \delta, \, i=0,1,
\]
for any $r \ge \bar{r}$ and any $\tau \in T(A)$.
\end{itemize} 
Then, 
\[
d_{i}:= [(\tilde{d}^{(i)}_{r})_{r \in \mathbb{N}} ] \in A_{\infty}, \, i=0,1,
\]
will have the desired properties: They are positive contractions because all the $\tilde{d}^{(i)}_{r}$ are; they are orthogonal because of (a) and they commute with $d$ because of (b). Finally, (c) shows that for any $\tau \in T_{\infty}(A)$ and any $f \in \mathcal{C}_{0}((0,1])_{+}$ we have 
\[
\tau(d_{i} f(d)) \ge \frac{1}{4(m+1)} \cdot \tau(f(d)), \, i=0,1.
\]
\end{nproof}
\en

\bn
The following will play a similar role for the proof of Lemma~\ref{3-10prime} as \cite[Proposition~3.7]{Winter:dr-Z-stable} does for \cite[Lemma~3.10]{Winter:dr-Z-stable}

\label{Enew}
\begin{nprop}
Let $A$ be a separable, simple, unital  $\mathrm{C}^{*}$-algebra with  tracial $m$-almost divisibility, let $0 \neq k \in \mathbb{N}$ and $d \in A_{\infty}$ be a positive contraction.

Then, there is a c.p.c.\ order zero map 
\[
\Phi:M_{k} \to A_{\infty} \cap \{d\}'
\]
satisfying
\begin{equation}
\label{w4-4-1}
\tau(\Phi(\be_{k}) d) \ge \frac{1}{4(m+1)^{2}}\cdot \tau(d)
\end{equation}
for all $\tau \in T_{\infty}(A)$.
\end{nprop}

\begin{nproof}
For a given $\eta>0$, choose $L \in \mathbb{N}$ with 
\begin{equation}
\label{w4-8}
1/L < \eta
\end{equation} 
and a standard partition of unity of the interval $(0,1]$
\begin{equation}
\label{w4-6}
h^{(i)}_{1}, \ldots, h^{(i)}_{L} \in \mathcal{C}_{0}((0,1]), \; i=0,1,
\end{equation}
as in \ref{standard-partition}.

Choose orthogonal positive contractions 
\begin{equation}
\label{w4-6a}
d_{0},d_{1} \in A_{\infty} \cap \{d\}'
\end{equation} 
as in Proposition~\ref{Fnew}. By Proposition~\ref{Hnew}, there are c.p.c.\ order zero maps 
\begin{equation}
\label{w4-7}
\varphi^{(i)}_{\eta,l}: M_{k} \to \overline{d_{i} h_{l}^{(i)}(d) A_{\infty} d_{i} h_{l}^{(i)}(d)}
\end{equation}
for $i=0,1$ and $l=1, \ldots,L$ satisfying
\begin{equation}
\label{w4-11}
\tau(\varphi^{(i)}_{\eta,l}(\be_{k})) \ge \frac{1}{m+1} \cdot \tau(d_{i} h_{l}^{(i)}(d)) - \frac{\eta}{L} 
\end{equation}
for all $\tau \in T_{\infty}(A)$.

The elements $d_{i} h_{l}^{(i)}(d) $ are pairwise orthogonal (using \eqref{w4-6a}, and since $d_{0} \perp d_{1}$,  and since the $h^{(i)}_{l}$ are pairwise orthogonal for fixed $i$), and so  the images of the $\varphi^{(i)}_{\eta,l}$ are pairwise orthogonal by \eqref{w4-7}, whence the $\varphi^{(i)}_{\eta,l}$ add up to a c.p.c.\ order zero map
\begin{equation}
\label{w4-9}
\varphi_{\eta}:= \bigoplus_{i=0,1;l=1,\ldots,L} \varphi^{(i)}_{\eta,l}: M_{k} \to A_{\infty}.
\end{equation}
One checks that
\begin{equation}
\label{w4-13}
\left\|  \sum_{i,l}   f(h^{(i)}_{l} (d)) d - \sum_{i,l} \frac{l}{L} \cdot f( h^{(i)}_{l}(d))  \right\| \le \frac{1}{L} \stackrel{\eqref{w4-8}}{<} \eta
\end{equation}
for any positive $f \in \mathcal{C}((0,1])$ of norm at most one, 
from which follows  that 
\begin{equation}
\label{w4-10}
\left\| \bigoplus_{i,l} \varphi^{(i)}_{\eta,l}(x)  d - \bigoplus_{i,l} \frac{l}{L} \cdot \varphi^{(i)}_{\eta,l} (x)\right\| \le \eta \|x\| 
\end{equation}
for all $x \in M_{k}$. We obtain
\begin{equation}
\label{w4-15}
\|[\varphi_{\eta}(x), d] \| \le 2 \eta \|x\| 
\end{equation}
for all $x \in M_{k}$.

Moreover, we have 
\begin{eqnarray}
\tau(\varphi_{\eta}(\be_{k})d)&\stackrel{\eqref{w4-9}}{=} & \tau\left( \sum_{i,l} \varphi_{\eta,l}^{(i)}(\be_{k}) d   \right) \nonumber \\
& \stackrel{\eqref{w4-10}}{\ge} & \sum_{i=0,1} \tau\left( \sum_{l=1,\ldots,L}  \frac{l}{L} \cdot \varphi^{(i)}_{\eta,l} (\be_{k})    \right) -  \eta \nonumber \\
& \stackrel{\eqref{w4-11}}{\ge} & \sum_{i=0,1} \frac{1}{m+1} \cdot \tau\left( \sum_{l=1,\ldots,L}  \frac{l}{L} \cdot d_{i} h_{l}^{(i)}(d)    \right) - 3 \eta \nonumber \\
& \stackrel{\eqref{w4-12}}{\ge} & \sum_{i=0,1} \frac{1}{4(m+1)^{2}} \cdot \tau\left( \sum_{l=1,\ldots,L}  \frac{l}{L} \cdot  h_{l}^{(i)}(d)    \right) - 3 \eta \nonumber \\
& \stackrel{\eqref{w4-13},\eqref{w4-6}}{\ge} & \frac{1}{4(m+1)^{2}} \cdot \tau(d) - 4 \eta \label{w4-14}
\end{eqnarray}
for any $\tau \in T_{\infty}(A)$. 

For $l \in \mathbb{N}$, we may now take $\frac{1}{l+1}$ in place of $\eta$, and apply the preceding construction to obtain c.p.c.\ order zero maps 
\[
\varphi_{\frac{1}{l+1}} :M_{k} \to A_{\infty}
\]
such that
\begin{equation}
\label{Enew1}
\tau(\varphi_{\frac{1}{l+1}} (\be_{k}) d) \stackrel{\eqref{w4-14}}{\ge} \frac{1}{4(m+1)^{2}} \cdot \tau(d) - \frac{4}{l+1} 
\end{equation}
for all $\tau \in T_{\infty}(A)$ and
\begin{equation}
\label{w4-17}
\|[\varphi_{\frac{1}{l+1}}(x), d ] \| \stackrel{\eqref{w4-15}}{\le} \frac{2}{l+1} \|x\| 
\end{equation}
for all $x \in M_{k}$.

Let 
\[
(d_{r})_{r \in \mathbb{N}} \in \prod_{\mathbb{N}} A
\]
be a sequence of positive contractions lifting $d$, and let
\[
(\bar{\varphi}_{l,r})_{r \in \mathbb{N}}: M_{k} \to \prod_{\mathbb{N}} A
\] 
be a  c.p.c.\ order zero map lifting $\varphi_{\frac{1}{l+1}}$. 

We claim that, for any $l \in \mathbb{N}$, there is $K_{l} \in \mathbb{N}$ such that
\begin{equation}
\label{w4-18}
\tau(\bar{\varphi}_{l,r}(\be_{k}) d_{r}) \ge \frac{1}{4(m+1)^{2}} \cdot \tau(d_{r}) - \frac{5}{l+1}
\end{equation}
for each $\tau \in T(A)$ and each $r \ge K_{l}$. 

If not, for some $l \in \mathbb{N}$ there were an increasing sequence $(r_{s})_{s \in \mathbb{N}} \subset \mathbb{N}$ and a sequence $(\tau_{s})_{s \in \mathbb{N}} \subset T(A)$ satisfying
\[
\tau_{s}(\bar{\varphi}_{l,r_{s}}(\be_{k}) d_{r_{s}}) < \frac{1}{4(m+1)^{2}} \cdot \tau(d_{r_{s}})  - \frac{5}{l+1},
\] 
a contradiction to \eqref{Enew1}, hence proving the claim (cf.\ the argument for \eqref{w4-16}).

Upon increasing the $K_{l}$, if necessary, we may also assume that they form a strictly increasing sequence and that 
\begin{equation}
\label{www4-4-1}
\|[\bar{\varphi}_{l,r}(x),d_{r} ]\| \stackrel{\eqref{w4-17}}{\le} \frac{3}{l+1}  \|x\| 
\end{equation}
for all  $x \in M_{k}$  and $r \ge K_{l}$.
 
We now define a sequence 
\[
(\tilde{\varphi}_{r}:M_{k} \to A)_{r \in \mathbb{N}} 
\]
of c.p.c.\ order zero maps by
\[
\tilde{\varphi}_{r}:= 
\begin{cases}
\bar{\varphi}_{l,r} & \mbox{ if } K_{l} \le r < K_{l+1} \\
0 & \mbox{ if } r<K_{0}.
\end{cases}
\]
From \eqref{www4-4-1} and \eqref{w4-18} we obtain
\begin{itemize}
\item[(a)] $\|[\tilde{\varphi}_{r}(x), d_{r} ]\| \stackrel{r \to \infty}{\longrightarrow} 0$ for all $x \in M_{k}$
\item[(b)] for any $\delta>0$ there is $\bar{r} \in \mathbb{N}$ such that 
\[
\tau(\tilde{\varphi}_{r}(\be_{k})d_{r}) \ge \frac{1}{4(m+1)^{2}} \cdot \tau(d_{r}) - \delta
\]
for any $r \ge \bar{r}$ and $\tau \in T(A)$.
\end{itemize}
Define 
\[
\Phi:M_{k} \to A_{\infty}
\]
by 
\[
\Phi(x):= [(\tilde{\varphi}_{r}(x))_{r \in \mathbb{N}} ],
\]
then $\Phi$ is a c.p.c.\ order zero map (since all the $\tilde{\varphi}_{r}$ are) satisfying
\[
\|[\Phi(M_{k}),d ] \| = 0
\]
(by (a) above) and 
\[
\tau(\Phi(\be_{k}) d) \ge \frac{1}{4(m+1)^{2}} \cdot \tau(d)
\]
for any $\tau \in T_{\infty}(A)$ (by (b) above).
\end{nproof}
\en

\bn
The next proposition generalizes \cite[3.8]{Winter:dr-Z-stable}.

\label{Dnew}
\begin{nprop}
Let $A$ be a separable, simple, unital  $\mathrm{C}^{*}$-algebra with  tracial $m$-almost divisibility, and let $0 \neq k \in \mathbb{N}$. Suppose further that $F$ is a finite-dimensional $\mathrm{C}^{*}$-algebra and let 
\[
\varphi:F \to A_{\infty}
\]
be a c.p.c.\ order zero map.

Then, there is a c.p.c.\ order zero map 
\[
\Phi:M_{k} \to A_{\infty} \cap \varphi(F)'
\]
satisfying
\begin{equation}
\label{w4-22}
\tau(\Phi(\be_{k}) \varphi(x)) \ge \frac{1}{4(m+1)^{2}}\cdot \tau(\varphi(x))
\end{equation}
for all $x \in F_{+}$ and $\tau \in T_{\infty}(A)$.
\end{nprop}

\begin{nproof}
Write 
\[
F = M_{r_{1}} \oplus \ldots \oplus M_{r_{s}}
\]
and let 
\[
\varphi^{(i)}:M_{r_{i}} \to A_{\infty}
\]
denote the (pairwise orthogonal) components of $\varphi$; let $\pi^{(i)}$ denote the respective supporting $*$-homomorphisms. 

For  $i = 1, \ldots, s$, apply Proposition~\ref{Enew} to obtain c.p.c.\ order zero maps
\begin{equation}
\label{w4-5-2}
\Phi^{(i)}:M_{k} \to A_{\infty}\cap \varphi^{(i)}(e_{11}^{(i)})'
\end{equation}
satisfying
\[
\tau(\Phi^{(i)} (\be_{k}) \varphi^{(i)} (e_{11}^{(i)}) ) \stackrel{\eqref{w4-4-1}}{\ge} \frac{1}{4(m+1)^{2}} \cdot \tau(\varphi^{(i)} (e_{11}^{(i)}))
\]
for all $\tau \in T_{\infty}(A)$, where $\{e_{pq}^{(i)}\}$ denotes a set of matrix units for $M_{r_{i}}$. 

For $\eta>0$, we may define pairwise orthogonal maps 
\begin{equation}
\label{w4-5-3}
\hat{\Phi}^{(i)}_{\eta} :M_{k} \to \overline{(g_{0,\eta}(\varphi^{(i)} ) (\be_{k})) A_{\infty}   (g_{0,\eta}(\varphi^{(i)} ) (\be_{k})) }
\end{equation}
by setting 
\begin{equation}
\label{w4-5-1}
\hat{\Phi}^{(i)}_{\eta}(y) := \sum_{p=1}^{r_{i}} \pi^{(i)}(e_{p1}^{(i)}) (g_{0,\eta}(\varphi^{(i)}) (e_{11}^{(i)})) \Phi^{(i)} (y) \pi^{(i)}(e_{1p}^{(i)}) 
\end{equation}
for $y \in M_{k}$. One checks that 
\[
\hat{\Phi}_{\eta} := \bigoplus_{i=1,\ldots,s} \hat{\Phi}^{(i)}_{\eta} 
\]
is c.p.c.\ order zero, that
\begin{equation}
\label{w4-5-5}
[\hat{\Phi}_{\eta}(M_{k}),\varphi(F)] = 0,
\end{equation}
and that 
\begin{equation}
\label{w4-5-4}
\tau(\hat{\Phi}_{\eta} (\be_{k}) \varphi(x)) \ge \frac{1}{4(m+1)^{2}} \cdot \tau(\varphi(x)) - \eta
\end{equation}
for all $x \in F_{+}$ of norm at most one and for all $\tau \in T_{\infty}(A)$ (using \eqref{w4-5-2}, \eqref{w4-5-3}, \eqref{w4-5-1}). 

For $ l \in \mathbb{N}$ we may now take $\frac{1}{l+1}$ in place of $\eta$ and choose c.p.c.\ order zero maps 
\[
\bar{\Phi}_{l,r}:M_{k} \to A,\, r \in \mathbb{N},
\]
such that
\[
[(\bar{\Phi}_{l,r}(y))_{r \in \mathbb{N}} ] = \hat{\Phi}_{\frac{1}{l+1}}(y) 
\]
for $y \in M_{k}$.

Let 
\[
(\bar{\varphi}_{r}:F \to A)_{r \in \mathbb{N}}
\]
be a sequence of c.p.c.\ order zero maps lifting $\varphi$. Just as in the proof of \ref{Enew} (see the argument for \eqref{w4-18}), one checks by contradiction  to \eqref{w4-5-4} that for each $l \in \mathbb{N}$ there is $K_{l} \in \mathbb{N}$ such that 
\begin{equation}
\label{w4-5-6}
\tau(\bar{\Phi}_{l,r}(\be_{k}) \bar{\varphi}_{r}(x)) \ge \frac{1}{4(m+1)^{2}} \cdot \tau(\bar{\varphi}_{r}(x)) - \frac{2}{l+1}
\end{equation}
for all $x \in F_{+}$ of norm at most one and $r \ge K_{l}$. Upon increasing the $K_{l}$ if necessary, we may assume the sequence $(K_{l})_{l \in \mathbb{N}}$ to be strictly increasing, and (by \eqref{w4-5-5}) that
\begin{equation}
\label{w4-5-7}
\| [\bar{\Phi}_{l,r}(y) ,\bar{\varphi}_{r}(x) ] \| \le \frac{1}{l+1}
\end{equation}
for normalized elements $y \in M_{k}$ and $x \in F$, if $r \ge K_{l}$. 

Again as in the proof of \ref{Enew}, we may define a sequence 
\[
(\tilde{\Phi}_{r}:M_{k} \to A)_{r \in \mathbb{N}} 
\]
of c.p.c.\ order zero maps by
\[
\tilde{\Phi}_{r}:= 
\begin{cases}
\bar{\Phi}_{l,r} & \mbox{ if } K_{l}\le r < K_{l+1}\\
0 & \mbox{ if } r< K_{0};
\end{cases}
\]
using \eqref{w4-5-7} and \eqref{w4-5-6} one checks that
\begin{itemize}
\item[(a)] $\|[\tilde{\Phi}_{r}(y), \varphi_{r}(x) ]\| \stackrel{r \to \infty}{\longrightarrow} 0$ for all $x \in F$, $y \in M_{k}$
\item[(b)] for any $\delta>0$ there is $\bar{r} \in \mathbb{N}$ such that 
\[
\tau(\tilde{\Phi}_{r}(\be_{k}) \bar{\varphi}_{r}(x)) \ge \frac{1}{4(m+1)^{2}} \cdot \tau(\bar{\varphi}_{r}(x)) - \delta
\]
for any $r \ge \bar{r}$, $x \in F_{+}$ of norm at most one, and $ \tau \in T(A)$. 
\end{itemize}
Define 
\[
\Phi:M_{k} \to A_{\infty}
\]
by 
\[
\Phi(x):= [(\tilde{\Phi}_{r}(x))_{r \in \mathbb{N}} ],
\]
then $\Phi$ is a c.p.c.\ order zero map since the $\tilde{\Phi}_{r}$ are; $\Phi(M_{k})$ commutes with $\varphi(F)$ by (a) above and
\[
\tau(\Phi(\be_{k})\varphi(x)) \ge \frac{1}{4(m+1)^{2}} \cdot \tau(\varphi(x))
\]
for all $x \in F_{+}$ and $ \tau \in T_{\infty}(A)$ by (b) above.
\end{nproof}
\en

\bn
The next proposition will play a role for \ref{3-10prime} similar to that of \cite[Proposition~3.9]{Winter:dr-Z-stable} for the proof of \cite[Lemma~3.10]{Winter:dr-Z-stable}. The subsequent Propositions~\ref{Anew} through \ref{3-9prime} will get this analogy to work.  

\label{Dprime}
\begin{nprop}
Let $A$ be a separable, simple, unital  $\mathrm{C}^{*}$-algebra with  tracial $m$-almost divisibility, and let $0 \neq k \in \mathbb{N}$. Suppose further that $F$ is a finite-dimensional $\mathrm{C}^{*}$-algebra and let 
\[
\varphi:F \to A_{\infty}
\]
be a c.p.c.\ order zero map. Let 
\begin{equation}
\label{w4-6-1}
d \in A_{\infty} \cap \varphi(F)'
\end{equation} 
be a positive contraction.

Then, there are c.p.c.\ order zero maps 
\[
\Phi^{(0)}, \ldots, \Phi^{(3)}:M_{k} \to A_{\infty} \cap \varphi(F)' \cap \{d\}'
\]
satisfying
\[
\tau((\Phi^{(0)}(\be_{k}) + \ldots + \Phi^{(3)}(\be_{k})) f(d) \varphi(x)) \ge \frac{1}{4(m+1)^{2}}\cdot \tau(f(d) \varphi(x))
\]
for all $x \in F_{+}$, $f \in \mathcal{C}_{0}((0,1])_{+}$ and $\tau \in T_{\infty}(A)$.
\end{nprop}

\begin{nproof}
Let 
\[
(\mathcal{G}_{q})_{q\in \mathbb{N}}
\]
be a nested sequence of finite subsets of the unit ball of $\mathcal{C}_{0}((0,1])_{+}$ with dense union. For each $q \in \mathbb{N}$, find $ L_{q} \in \mathbb{N}$ such that the  standard partition of unity
\begin{equation}
\label{w4-23}
f_{q,1}^{(i)}, \ldots, f_{q,L_{q}}^{(i)} \in \mathcal{C}_{0}((0,1]) , \; i=0,1,
\end{equation}
of the interval $(0,1]$ as in \ref{standard-partition} satisfies 
\begin{eqnarray}
\left\| \bigoplus_{r,s=1}^{L_{q}} \left( f_{q,r}^{(i)}(d) f^{(j)}_{q,s}(\varphi(\be_{F})) f(d) \varphi(x) - f\left( \frac{r}{L_{q}}\right) \cdot \frac{s}{L_{q}} \cdot f^{(i)}_{q,r}(d) f^{(j)}_{q,s}(\varphi)(x) \right) \right\| && \nonumber \\
 \le  \frac{1}{q+1}&& \label{w4-21}
\end{eqnarray}
for $q \in \mathbb{N}$,   $i,j\in \{0,1\}$, $f \in \mathcal{G}_{q}$ and $x$ in the unit ball of $F_{+}$.

For each $q\in \mathbb{N}$, $r,s \in \{1,\ldots,L_{q}\}$, $i,j \in \{0,1\}$, the map
\begin{equation}
\label{ww4-6-1}
x \mapsto f^{(i)}_{q,r}(d) f^{(j)}_{q,s}(\varphi)(x)
\end{equation}
is c.p.c.\ order zero, so that we may  apply Proposition~\ref{Dnew} to obtain c.p.c.\ order zero maps
\begin{equation}
\label{w4-6-2}
\Phi^{(i,j)}_{q,r,s}:M_{k} \to A_{\infty} \cap (f^{(i)}_{q,r}(d) (f^{(j)}_{q,s}(\varphi)(F) ))'.
\end{equation}
For fixed $q$, $i$ and $j$, the $f^{(i)}_{q,r}(d) f^{(j)}_{q,s}(\varphi)(\be_{F})$ are pairwise orthogonal positive contractions by \eqref{w4-6-1} and \eqref{w4-23}.

We therefore have  c.p.c.\ order zero maps 
\[
\Phi^{(i,j)}_{q}:M_{k} \to A_{\infty} \cap \left(\bigoplus_{r,s} f^{(i)}_{q,r}(d) f^{(j)}_{q,s}(\varphi) (F) \right)'
\]
given by
\begin{equation}
\label{w4-6-3}
\Phi^{(i,j)}_{q}( \, . \,):=  \bigoplus_{r,s} \Phi^{(i,j)}_{q,r,s}(\, .\,) f^{(i)}_{q,r}(d) f^{(j)}_{q,s}(\varphi(\be_{F})).
\end{equation}
For fixed $q, i,j$, one checks that
\begin{eqnarray}
\left\| \Phi^{(i,j)}_{q} (y) f(d) - \bigoplus_{r,s} \Phi^{(i,j)}_{q,r,s}(y) \left(  f\left(\frac{r}{L_{q}}\right) \cdot f^{(i)}_{q,r} (d) f^{(j)}_{q,s}(\varphi)(\be_{F}) \right) \right\| && \nonumber \\
 \le \frac{1}{q+1} \cdot \|y\|   &&  \label{w4-19}
\end{eqnarray}
for $y \in (M_{k})_{+}$ and $f \in \mathcal{G}_{q}$, from which (in connection with \eqref{w4-6-2}) it follows that
\begin{equation}
\label{w4-25}
\|[\Phi^{(i,j)}_{q}(y) ,f(d) ]\| \le \frac{2}{q+1} \cdot \|y\|
\end{equation}
for $y \in (M_{k})_{+}$, $f \in \mathcal{G}_{q}$, $q \in \mathbb{N}$ and $i,j\in \{0,1\}$.

Similarly, we have
\begin{eqnarray}
\label{w4-20}
 \left\| \Phi^{(i,j)}_{q}(y) \varphi(x) - \bigoplus_{r,s} \Phi^{(i,j)}_{q,r,s}(y) f^{(i)}_{q,r}(d) \left( \frac{s}{L_{q}} \cdot f^{(j)}_{q,s}(\varphi)(x)\right) \right\| && \nonumber \\
\le   \frac{1}{q+1} \cdot \|x\| \|y\|   &&
\end{eqnarray}
for $x \in F_{+}$ and $y \in (M_{k})_{+}$, from which we see that
\begin{equation}
\label{w4-24}
\|[\Phi^{(i,j)}_{q}(y),\varphi(x) ]\| \le \frac{2}{q+1} \cdot \|x\| \|y\|
\end{equation}
for $x \in F_{+}$, $y \in (M_{k})_{+}$, $q \in \mathbb{N}$ and $i,j \in \{0,1\}$. 

Next, we check that for $q \in \mathbb{N}$, $x$ in the unit ball of $F_{+}$, $f \in \mathcal{G}_{q}$ and $\tau \in T_{\infty}(A)$ 
\begin{eqnarray}
\lefteqn{\tau\left(  \sum_{i,j \in \{0,1\}} \Phi^{(i,j)}_{q}(\be_{k}) f(d) \varphi(x)  \right)} \nonumber  \\
& \stackrel{\eqref{w4-6-3},\eqref{w4-21}}{\ge} & \tau\left(  \sum_{i,j \in \{0,1\}} \bigoplus_{r,s =1}^{L_{q}} \Phi^{(i,j)}_{q,r,s} (\be_{k}) \left(f\left( \frac{r}{L_{q}} \right) \cdot \frac{s}{L_{q}} \cdot f^{(i)}_{q,r}(d) f^{(j)}_{q,s}(\varphi)(x) \right)  \right) \nonumber \\
&&  - \frac{4}{q+1} \nonumber  \\
& \stackrel{\eqref{ww4-6-1},\eqref{w4-22}}{\ge} & \frac{1}{4(m+1)^{2}} \cdot \tau \left( \sum_{i,j} \bigoplus_{r,s} f\left( \frac{r}{L_{q}} \right) \cdot \frac{s}{L_{q}} \cdot f^{(i)}_{q,r}(d) f^{(j)}_{q,s}(\varphi)(x)  \right) \nonumber \\
&& - \frac{4}{q+1} \nonumber \\
& \stackrel{\eqref{w4-23},\eqref{w4-21}}{\ge} & \frac{1}{4(m+1)^{2}} \cdot \tau(f(d) \varphi(x)) - \frac{6}{q+1}. \label{Dprime2} 
\end{eqnarray}
For each $q,i,j$ let 
\[
\bar{\Phi}^{(i,j)}_{q}:M_{k} \to \prod_{\mathbb{N}}A 
\]
be a c.p.c.\ order zero lift for $\Phi^{(i,j)}_{q}$ (cf.\ \cite[Remark~2.4]{KirWinter:dr}); let 
\[
\bar{\Phi}^{(i,j)}_{q,l}:M_{k} \to A, \; l \in \mathbb{N} 
\]
denote the components. Let 
\[
(\bar{d}_{l})_{l \in \mathbb{N}} \in \prod_{\mathbb{N}} A
\]
be a positive contraction lifting $d$, and let 
\[
(\bar{\varphi}_{l})_{l \in \mathbb{N}}:F \to \prod_{\mathbb{N}} A
\]
be a c.p.c.\ order zero lift for $\varphi$. 

We claim that for each $q \in \mathbb{N}$ there is $R_{q}\in \mathbb{N}$ such that
\begin{equation}
\label{Dprime1}
\tau\left(  \sum_{i,j} \bar{\Phi}^{(i,j)}_{q,l}(\be_{k}) f(\bar{d}_{l}) \bar{\varphi}_{l}(x)  \right) \ge \frac{1}{4(m+1)^{2}} \cdot \tau(f(\bar{d}_{l}) \bar{\varphi}_{l}(x)) - \frac{7}{q+1}
\end{equation} 
for each $l \ge R_{q}$, $\tau \in T(A)$, $f \in \mathcal{G}_{q}$ and $x$ in the unit ball of $F_{+}$. 

To prove the claim, note that otherwise for some $q \in \mathbb{N}$ there were $f \in \mathcal{G}_{q}$, an increasing  sequence $(l_{t})_{t \in \mathbb{N}} \subset \mathbb{N}$, a sequence $(x_{l_{t}})_{t \in \mathbb{N}}$ in the unit ball of $F_{+}$, and a sequence $(\tau_{l_{t}})_{t \in \mathbb{N}} \subset T(A)$ such that 
\[
\tau_{l_{t}}\left(  \sum_{i,j} \bar{\Phi}^{(i,j)}_{q,l_{t}}(\be_{k}) f(\bar{d}_{l_{t}}) \bar{\varphi}_{l}(x_{l_{t}})  \right) < \frac{1}{4(m+1)^{2}} \cdot 
\tau_{l_{t}}(f(\bar{d}_{l_{t}}) \bar{\varphi}_{l_{t}}(x_{l_{t}})) - \frac{7}{q+1}
\]
for each $t \in \mathbb{N}$. By compactness and passing to a subsequence of $(l_{t})_{t \in \mathbb{N}}$ we may assume that $(x_{l_{t}})_{t\in \mathbb{N}}$ converges to some $\bar{x}$ in the unit ball of $F_{+}$. Let $\omega \in \beta\mathbb{N} \setminus \mathbb{N}$ be a free ultrafilter along the sequence $(l_{t})_{t \in \mathbb{N}}$, then
\[
\lim_{\omega} \tau_{l} \left(  \sum_{i,j} \bar{\Phi}^{(i,j)}_{q,l}(\be_{k}) f(\bar{d}_{l}) \bar{\varphi}_{l}(\bar{x})  \right) \le \frac{1}{4(m+1)^{2}} \cdot \lim_{\omega} \tau_{l}(f(\bar{d}_{l}) \bar{\varphi}_{l}(\bar{x})) - \frac{7}{q+1},
\]
a contradiction to \eqref{Dprime2}, since the map $[(a_{l})_{l\in \mathbb{N}} ] \mapsto \lim_{\omega} \tau_{l}(\, .\,)$ is in $T_{\infty}(A)$. This establishes the claim \eqref{Dprime1}.

We may also assume that the sequence $(R_{q})_{q \in \mathbb{N}}$ is strictly increasing and that, for each $q \in \mathbb{N}$, $R_{q}$ is so large that
\[
\|[\bar{\Phi}_{q,l}^{(i,j)}(y),f(\bar{d}_{l}) ]\| \le \frac{3}{q+1} \|y\|
\]
and
\[
\|[\bar{\Phi}^{(i,j)}_{q,l}(y),\bar{\varphi}_{l}(x) ] \| \le \frac{3}{q+1} \|x\| \|y\|
\]
for all $i,j \in \{0,1\}$, $x \in F_{+}$, $y \in (M_{k})_{+}$, $f \in \mathcal{G}_{q}$ and for all $l \ge R_{q}$ (using \eqref{w4-25} and \eqref{w4-24}). 

Now for $l \in \mathbb{N}$ we may define c.p.c.\ order zero maps 
\[
\tilde{\Phi}^{(i,j)}_{l}:M_{k} \to A
\]
by
\[
\tilde{\Phi}^{(i,j)}_{l}:= 
\begin{cases}
\bar{\Phi}^{(i,j)}_{q,l} & \mbox{ if } R_{q} \le l < R_{q+1} \\
0 & \mbox{ if } l < R_{0},
\end{cases}
\]
and let 
\[
\Phi^{(i,j)}: M_{k} \to A_{\infty}
\]
be the c.p.c.\ order zero map induced by the $\tilde{\Phi}^{(i,j)}_{l}$, $i,j \in \{0,1\}$. After relabeling the $\Phi^{(i,j)}$, it is easy to check that we indeed have maps $\Phi^{(0)}, \ldots, \Phi^{(3)}$ with the desired properties.
\end{nproof}
\en

\bn
\label{Anew}
Some of the ideas of the next two propositions we have already encountered in the proof of Proposition~\ref{2-d-i}, which was used to show that finite nuclear dimension implies tracial $\tilde{m}$-almost divisibility (Proposition~\ref{dimnuc-tracial-divisibility}). In the present context we assume tracial $\bar{m}$-almost divisibility to arrive at a statement similar to that of \ref{2-d-i}.

\begin{nprop}
For any $m, \bar{m} \in \mathbb{N}$, there is $\beta_{m,\bar{m}}>0$ such that the following holds:

Let $A$ be a separable, simple, unital  $\mathrm{C}^{*}$-algebra with  tracial $\bar{m}$-almost divisibility, and let $B \subset A_{\infty}$ be a separable, unital $\mathrm{C}^{*}$-subalgebra with $\dimnuc B \le m$ and 
\begin{equation}
\label{w4-7-4}
d \in A_{\infty} \cap B'
\end{equation}
a positive contraction. 

Then, there are orthogonal positive contractions 
\[
d^{(0)},d^{(1)} \in A_{\infty} \cap B' \cap \{d\}' 
\]
satisfying
\begin{equation}
\label{ww4-7-1}
\tau(d^{(j)}f(d)b) \ge \beta_{m,\bar{m}} \cdot \tau(f(d)b)
\end{equation}
for all $b \in B_{+}$, $f \in \mathcal{C}_{0}((0,1])_{+}$, $\tau \in T_{\infty}(A)$ and $j=0,1$.
\end{nprop}

\begin{nproof}
Take
\begin{equation}
\label{w4-7-3}
\beta_{m,\bar{m}}:= \frac{1}{32(m+1) (\bar{m}+1)^{2}}.
\end{equation}
Let $A,B,d$ be as in the proposition. Let $(\mathcal{G}_{p})_{p \in \mathbb{N}}$ be a nested sequence of finite subsets (each containing $\be_{B}$) of the unit ball of $B_{+}$ with dense union.

By Proposition~\ref{central-cutout}, there is a system of $m$-decomposable c.p.\ approximations for $B$
\[
(F_{p}= F^{(0)}_{p} \oplus \ldots \oplus F^{(m)}_{p}, \psi_{p}, \varphi_{p}= \varphi^{(0)}_{p} + \ldots + \varphi^{(m)}_{p})_{p\in \mathbb{N}}
\]
such that the 
\[
\psi_{p}:B \to F_{p}
\]
are c.p.c.\ maps,  the 
\[
\varphi^{(i)}_{p}:F^{(i)}_{p} \to B
\]
are c.p.c.\ order zero maps, and such that 
\begin{equation}
\label{w4-28}
\|\varphi_{p} \psi_{p}(b) - b\| < \frac{1}{p+1} ,
\end{equation}
\begin{equation}
\label{w4-29}
\|[\varphi_{p} \psi_{p}(b) , \varphi^{(i)}_{p}(\psi^{(i)}_{p}(\be_{B})) ]\| < \frac{1}{p+1}
\end{equation}                                          
and
\begin{equation}
\label{w4-30}
\|b \varphi^{(i)}_{p}(\psi^{(i)}_{p}(\be_{B})) - \varphi^{(i)}_{p} \psi^{(i)}_{p}(b) \| < \frac{1}{p+1}
\end{equation}
for all $b \in \mathcal{G}_{p}$, $i \in \{0, \ldots,m\}$ and $p \in \mathbb{N}$. 

For each $i \in \{0, \ldots,m\}$ and $p \in \mathbb{N}$ apply Proposition~\ref{Dprime} (with $\bar{m}$ in place of $m$, $\varphi^{(i)}_{p}$ in place of $\varphi$ and $8 (m+1)$ in place of $k$) to obtain c.p.c.\ order zero maps 
\begin{equation}
\label{w4-26}
\Phi^{(i,j)}_{p}: M_{8(m+1)} \to A_{\infty} \cap \varphi^{(i)}_{p}(F^{(i)}_{p})' \cap \{d\}'
\end{equation}
for $j=0,1,2,3$ satisfying
\begin{equation}
\label{w4-27}
\tau\left( \sum_{j=0}^{3} \Phi_{p}^{(i,j)} (\be_{8(m+1)})f(d) \varphi_{p}^{(i)}(x)    \right) \ge \frac{1}{4(\bar{m}+1)^{2}} \cdot \tau(f(d) \varphi^{(i)}_{p}(x))
\end{equation}
for all $x \in (F^{(i)}_{p})_{+}$, $f \in \mathcal{C}_{0}((0,1])_{+}$ and $\tau \in T_{\infty}(A)$.

But then one checks that, for all $i\in \{0,\ldots,m\}$, $j\in \{0,\ldots,3\}$ and $p\in \mathbb{N}$,
\begin{equation}
\label{w4-31}
\check{\Phi}^{(i,j)}_{p}(\, .\,) := \Phi^{(i,j)}_{p}(\, .\,) \varphi^{(i)}_{p}(\psi^{(i)}_{p}(\be_{B})) :M_{8(m+1)} \to A_{\infty} \cap \{d\}'
\end{equation}
define c.p.c.\ order zero maps (using \eqref{w4-26} and \eqref{w4-7-4}).

Furthermore, we have
\begin{eqnarray*}
\lefteqn{\| \check{\Phi}_{p}^{(i,j)}(y)b - b \check{\Phi}_{p}^{(i,j)}(y) \|  } \\
&\stackrel{\eqref{w4-31},\eqref{w4-30},\eqref{w4-26}}{\le}&\| \Phi_{p}^{(i,j)}(y)\varphi_{p}^{(i)}\psi_{p}^{(i)}(b) - \varphi_{p}^{(i)}\psi_{p}^{(i)}(b) \Phi_{p}^{(i,j)}(y)\| + \frac{2}{p+1} \|y\| \\
& \stackrel{\eqref{w4-26}}{=} & \frac{2}{p+1} \|y\|
\end{eqnarray*}
for $b \in \mathcal{G}_{p}$ and $y \in M_{8(m+1)}$, whence
\[
\|[\check{\Phi}^{(i,j)}_{p}(y),b ]\| \stackrel{p \to \infty}{\longrightarrow} 0 
\]
for  $b \in B$ and $ y \in M_{8(m+1)}$, and
\begin{eqnarray*}
\lefteqn{
\tau\left(  \sum_{j=0}^{3} \sum_{i=0}^{m} \check{\Phi}^{(i,j)}_{p} (\be_{8(m+1)}) f(d)b   \right)  } \\
& \stackrel{\eqref{w4-31},\eqref{w4-26}, \eqref{w4-30}}{\ge} & \tau\left(  \sum_{j=0}^{3} \sum_{i=0}^{m} \Phi^{(i,j)}_{p} (\be_{8(m+1)}) f(d) \varphi_{p}^{(i)} \psi_{p}^{(i)}(b)   \right) - \frac{4(m+1) +1}{p+1} \\
& \stackrel{\eqref{w4-27},\eqref{w4-28}}{\ge} & \frac{1}{4(\bar{m}+1)^{2}} \cdot \tau(f(d)b) - \frac{4(m+1) +1}{p+1}
\end{eqnarray*}
for all $b \in \mathcal{G}_{p}$, $f$ in the unit ball of $\mathcal{C}_{0}((0,1])_{+}$  and $\tau \in T_{\infty}(A)$. 

After choosing a c.p.c.\ order zero lift 
\[
\hat{\Phi}^{(i,j)}_{p}: M_{8(m+1)} \to \prod_{\mathbb{N}} A
\]
for each $\check{\Phi}^{(i,j)}_{p}$ (cf.\ \cite[Remark~2.4]{KirWinter:dr}), a standard diagonal sequence argument yields c.p.c.\ order zero maps 
\[
\bar{\Phi}^{(i,j)}:M_{8(m+1)} \to A_{\infty} \cap B' \cap \{d\}'
\]
such that
\begin{equation}
\label{w4-32}
\tau\left(  \sum_{j=0}^{3} \sum_{i=0}^{m} \bar{\Phi}^{(i,j)} (\be_{8(m+1)}) f(d)b   \right) \ge \frac{1}{4(\bar{m}+1)^{2}} \cdot \tau(f(d)b) 
\end{equation}
for all $b \in B$, $f \in \mathcal{C}_{0}((0,1])_{+}$  and $\tau \in T_{\infty}(A)$.  

Now for each $q \in \mathbb{N}$ we may define orthogonal positive contractions 
\begin{equation}
\label{w4-33}
d^{(0)}_{q} := f_{\frac{1}{q+1}, \frac{2}{q+1}}\left( \sum_{i,j} \bar{\Phi}^{(i,j)} (e_{11})   \right) \in A_{\infty} \cap B' \cap\{d\}'
\end{equation}
and
\begin{equation}
\label{w4-34}
d^{(1)}_{q} := \be_{A_{\infty}} - g_{0, \frac{1}{q+1}}\left( \sum_{i,j} \bar{\Phi}^{(i,j)} (e_{11})   \right) \in A_{\infty} \cap B' \cap\{d\}'.
\end{equation}
We estimate
\begin{eqnarray}
\tau(d^{(0)}_{q} f(d)b) & \stackrel{\eqref{w4-33}}{\ge} & \tau\left( \sum_{i,j} \bar{\Phi}^{(i,j)} (e_{11}) f(d)b   \right) - \frac{1}{q+1} \nonumber  \\
& \stackrel{\eqref{w4-32}}{\ge} & \frac{1}{4 (\bar{m}+1)^{2} \cdot 8 (m+1)} \cdot \tau(f(d)b) - \frac{1}{q+1} \label{w4-7-2}
\end{eqnarray}
for all $b \in B_{+}$ and $f \in \mathcal{C}_{0}((0,1])_{+}$ of norm at most one and for all  $\tau \in T_{\infty}(A)$.  Moreover,
\begin{eqnarray}
\tau((\be_{A_{\infty}} - d^{(1)}_{q}) f(d)b) & \stackrel{\eqref{w4-34}}{=} & \tau\left(  g_{0,\frac{1}{q+1}}\left(  \sum_{i,j} \bar{\Phi}^{(i,j)}(e_{11})  \right) f(d)b      \right) \nonumber \\
& \le & \lim_{l \to \infty} \tau\left( \left(  \sum_{i,j} \bar{\Phi}^{(i,j)}(e_{11}) \right)^{\frac{1}{l}} f(d)b     \right) \nonumber  \\
& \le & \sum_{i=0,\ldots,m;j=0,\ldots,3} \lim_{l \to \infty} \tau\left( \left(  \bar{\Phi}^{(i,j)}(e_{11}) \right)^{\frac{1}{l}} f(d)b     \right) \nonumber \\
& \le & \frac{4(m+1)}{8(m+1)} \cdot \tau(f(d)b) \nonumber  \\
& = & \frac{1}{2} \cdot \tau(f(d)b) \label{w4-7-1}
\end{eqnarray}
for $\tau \in T_{\infty}(A)$ and $b \in B_{+}$, where for both the second and third estimate we have used that $\tau(\, .\, f(d)b)$ is a trace on $A_{\infty} \cap B' \cap\{d\}'$, whence  $\lim_{l \to \infty} \tau((\, .\,)^{\frac{1}{l}}f(d)b)$ is a dimension function on $A_{\infty} \cap B' \cap\{d\}'$ for any $\tau \in T_{\infty}(A)$.  

As a consequence of \eqref{w4-7-1}, 
\[
\tau(d^{(1)}_{q} f(d)b) \ge \frac{1}{2} \cdot \tau(f(d)b)
\]
for $b \in B_{+}$, $f \in \mathcal{C}_{0}((0,1])_{+}$ and $\tau \in T_{\infty}(A)$. 

Applying a diagonal sequence argument once more, from \eqref{w4-7-2} and \eqref{w4-7-1} we obtain orthogonal positive contractions
\[
d^{(0)}, d^{(1)} \in A_{\infty} \cap B' \cap \{d\}'
\] 
such that
\[
\tau(d^{(j)} f(d)b) \ge \frac{1}{32(m+1) (\bar{m}+1)^{2}} \cdot \tau(f(d)b)
\]
for $b \in B_{+}$, $f \in \mathcal{C}_{0}((0,1])_{+}$ and $\tau \in T_{\infty}(A)$.
\end{nproof}
\en

\bn
\label{Bnew}
\begin{nprop}
For any $m,\bar{m},k \in \mathbb{N}$, there is $\gamma_{m,\bar{m},k}>0$ such that the following holds:

Let $A$ be a separable, simple, unital  $\mathrm{C}^{*}$-algebra with  tracial $\bar{m}$-almost divisibility, and let $B \subset A_{\infty}$ be a separable, unital $\mathrm{C}^{*}$-subalgebra with $\dimnuc B \le m$.  Let $d \in A_{\infty} \cap B'$ be a positive contraction. 

Then, there are pairwise orthogonal positive contractions 
\[
d^{(0)},\ldots ,d^{(k)} \in A_{\infty} \cap B' \cap \{d\}' 
\]
satisfying
\[
\tau(d^{(i)}f(d)b) \ge \gamma_{m,\bar{m},k} \cdot \tau(f(d)b)
\]
for all $b \in B_{+}$, $f \in \mathcal{C}_{0}((0,1])_{+}$, $\tau \in T_{\infty}(A)$ and $i = 0,\ldots,k$.
\end{nprop}

\begin{nproof}
Let us first assume that $k=2^{l}-1$ for some $l \in \mathbb{N}$.  For $k$ of this specific form, we verify the statement by induction over $l$.

For $l=0$, the statement is trivial with $\gamma_{m,\bar{m},0}= 1$ and $d^{(0)}=\be_{A_{\infty}}$. Suppose now the statement has been verified for $k=2^{l}-1$ for some $l \in \mathbb{N}$, with elements 
\[
\bar{d}^{(0)}, \ldots, \bar{d}^{(2^{l}-1)} \in A_{\infty} \cap B' \cap \{d\}'. 
\]
Set
\[
\tilde{B}^{(i)}:= \mathrm{C}^{*}(B,\bar{d}^{(i)}) \subset A_{\infty}
\]
and note that for each $i \in \{0, \ldots, 2^{l}-1\}$, $\tilde{B}^{(i)}$ is a quotient of $\mathcal{C}_{0}([0,1]) \otimes B$ by \ref{order-zero-facts}(i), whence 
\[
\dimnuc \tilde{B}^{(i)} \le 2m+1
\]
by \cite{WinterZac:dimnuc}. Set
\[
\gamma_{m,\bar{m},2^{l+1}-1} :=  \beta_{2m+1,\bar{m}} \cdot  \gamma_{m,\bar{m},2^{l}-1},
\]
where $\beta_{2m+1,\bar{m}}$ comes from \ref{Anew}. 

For each $i \in \{0, \ldots, 2^{l}-1\}$, apply Proposition~\ref{Anew} (with $2m+1$ in place of $m$ and $\tilde{B}^{(i)}$ in place of $B$) to obtain orthogonal positive contractions 
\[
d^{(i,j)} \in A_{\infty} \cap B' \cap \{d\}' \cap \{\bar{d}^{(i)}\}', \; j=0,1,
\]  
satisfying 
\begin{eqnarray*}
\tau(d^{(i,j)} \bar{d}^{(i)}db) &\stackrel{\eqref{ww4-7-1}}{\ge} & \beta_{2m+1,\bar{m}} \cdot \tau( \bar{d}^{(i)}f(d)b) \\
& \ge & \beta_{2m+1,\bar{m}}\cdot \gamma_{m,\bar{m},2^{l}-1} \cdot \tau( f(d)b) \\
& = & \gamma_{m,\bar{m},2^{l+1}-1} \cdot \tau( f(d)b)
\end{eqnarray*}
for all $b \in B_{+}$, $f \in \mathcal{C}_{0}((0,1])_{+}$ and $\tau \in T_{\infty}(A)$. This yields the assertion for $k=2^{l+1}-1$, where the $2^{l+1}$ pairwise orthogonal  positive contractions are given by 
\[
d^{(i,j)}\bar{d}^{(i)} \in A_{\infty} \cap B' \cap \{d\}'
\]
for  $i=0, \ldots, 2^{l}-1$  and  $j=0,1$. 

Induction yields the statement   of the proposition for any $k$ of the form $2^{l}-1$, $l \in \mathbb{N}$. This clearly also implies the case of arbitrary $k \in \mathbb{N}$.
\end{nproof}
\en

\bn
\label{Cnew}
\begin{nprop}
For any $m,\bar{m}, k \in \mathbb{N}$, there is $\gamma_{m,\bar{m},k}>0$ such that the following holds:

Let $A$ be a separable, simple, unital  $\mathrm{C}^{*}$-algebra with  tracial $\bar{m}$-almost divisibility, and let $B \subset A_{\infty}$ be a separable, unital $\mathrm{C}^{*}$-subalgebra with $\dimnuc B \le m$.  Let $l \in \mathbb{N}$ and let  
\[
\Phi: M_{l} \to A_{\infty} \cap B'
\]
be a c.p.c.\ order zero map.

Then, there are pairwise orthogonal positive contractions 
\[
d^{(0)},\ldots ,d^{(k)} \in A_{\infty} \cap B' \cap \Phi(M_{l})' 
\]
satisfying
\[
\tau(d^{(i)} f(\Phi)(x)b) \ge \gamma_{m,\bar{m},k} \cdot \tau(f(\Phi)(x)b)
\]
for all $x \in (M_{l})_{+}$, $f \in \mathcal{C}_{0}((0,1])_{+}$, $b \in B_{+}$, $\tau \in T_{\infty}(A)$ and $i = 0,\ldots,k$.
\end{nprop}

\begin{nproof}
Apply Proposition~\ref{Bnew} (with $\Phi(e_{11})$ in place of $d$) to obtain pairwise orthogonal positive contractions
\[
\bar{d}^{(0)},\ldots, \bar{d}^{(k)} \in A_{\infty} \cap B' \cap \{\Phi(e_{11})\}'
\]
such that 
\begin{equation}
\label{w4-36}
\tau(\bar{d}^{(i)} f(\Phi(e_{11}))b) \ge \gamma_{m,\bar{m},k} \cdot \tau(f(\Phi(e_{11}))b)
\end{equation}
for $i=0,\ldots,k$, $f \in \mathcal{C}_{0}((0,1])_{+}$, $b \in B_{+}$ and $\tau\in T_{\infty}(A)$. 

For $q \in \mathbb{N}$, define pairwise orthogonal positive contractions
\begin{equation}
\label{w4-35}
\bar{d}_{q}^{(i)} := \bar{d}^{(i)} g_{0,\frac{1}{q+1}} (\Phi(e_{11})) \in A_{\infty} \cap B' \cap \{\Phi(e_{11})\}';
\end{equation}
then
\[
\tau(\bar{d}_{q}^{(i)} f(\Phi)(e_{11})b) \stackrel{\eqref{w4-35},\eqref{w4-36}}{\ge} \gamma_{m,\bar{m},k} \cdot \tau(f(\Phi)(e_{11})b) 
\]
for $i=0, \ldots,k$, for all positive contractions $b \in B$, for  all $f \in \mathcal{C}_{0}((\frac{1}{q+1},1])_{+}$ and for all $\tau \in T_{\infty}(A)$. 

Set 
\[
d_{q}^{(i)}:= \sum_{j=1}^{l} \pi_{\Phi}(e_{j1}) \bar{d}_{q}^{(i)} \pi_{\Phi}(e_{1j}), \, i=0,\ldots,k,
\]
where $\pi_{\Phi}$ denotes the supporting $*$-homomorphism of $\Phi$;  one checks that  these are pairwise orthogonal positive contractions, that 
\[
d_{q}^{(i)} \in A_{\infty} \cap B' \cap \Phi(M_{l})'
\] 
and that
\[
\tau(d_{q}^{(i)} f(\Phi)(x)b) \ge \gamma_{m,\bar{m},k} \cdot \tau(f(\Phi)(x)b) 
\]
for all positive contractions $x \in M_{l}$ and $b \in B$, for all $f \in \mathcal{C}_{0}((\frac{1}{q+1},1])_{+}$ and for all $\tau \in T_{\infty}(A)$.

After taking positive contractions
\[
\tilde{d}^{(i)}_{q} \in \prod_{\mathbb{N}} A
\]
and a c.p.c.\ order zero map
\[
\tilde{\Phi}:M_{l} \to \prod_{\mathbb{N}} A
\]
lifting the $d^{(i)}_{q}$ and $\Phi$, respectively, a standard diagonal sequence argument  yields orthogonal positive contractions 
\[
d^{(i)} \in A_{\infty}\cap B' \cap \Phi(M_{l})', \; i=0, \ldots,k,
\] 
such that the assertion of the Proposition holds.
\end{nproof}
\en

\bn
\label{3-9prime}
\begin{nprop}
For any $m,\bar{m} \in \mathbb{N}$, there is $\alpha_{m,\bar{m}}>0$ such that the following holds:

Let $A$ be a separable, simple, unital  $\mathrm{C}^{*}$-algebra with tracial $\bar{m}$-almost divisibility, and let $B \subset A$ be a unital $\mathrm{C}^{*}$-subalgebra with $\dimnuc B \le m$.  Let $k, l \in \mathbb{N}$ and let  
\begin{equation}
\label{w4-40}
\Phi: M_{l} \to A_{\infty} \cap B'
\end{equation}
be a c.p.c.\ order zero map.

Then, there is a c.p.c.\ order zero map 
\[
\Psi: M_{k} \to A_{\infty} \cap B' \cap \Phi(M_{l})'
\]
satisfying
\[
\tau(\Psi(\be_{k}) \Phi(\be_{l})b) \ge \alpha_{m,\bar{m}} \cdot \tau(\Phi(\be_{l})b)
\]
for all $b \in B_{+}$ and $\tau \in T_{\infty}(A)$.
\end{nprop}

\begin{nproof}
Set
\begin{equation}
\label{w4-50}
\alpha_{m,\bar{m}}:= \frac{\gamma_{m,\bar{m},m} \cdot \gamma_{m,\bar{m},1} }{4(\bar{m}+1)^{2}},
\end{equation}
where $\gamma_{m,\bar{m},l}$ comes from Proposition~\ref{Cnew}. 

Let $A,B,k,l,\Phi$ be as in the proposition. Let $(\mathcal{G}_{p})_{p \in \mathbb{N}}$ be a nested sequence of finite subsets (each containing $\be_{B}$) of the unit ball of $B_{+}$ with dense union. 

By Proposition~\ref{central-cutout}, there is a system of $m$-decomposable c.p.\ approximations for $B$
\[
(F_{p}=F_{p}^{(0)} \oplus \ldots \oplus F_{p}^{(m)}, \psi_{p},\varphi_{p})_{p \in \mathbb{N}}
\]
such that the 
\[
\psi_{p}:B \to F_{p}
\]
are c.p.c.\ maps,  the 
\begin{equation}
\label{w4-10-1}
\varphi^{(i)}_{p}:F^{(i)}_{p} \to B
\end{equation}
are c.p.c.\ order zero maps, and such that 
\begin{equation}
\label{w4-49}
\|\varphi_{p} \psi_{p}(b) - b\| < \frac{1}{p+1} ,
\end{equation}
\[
\|[\varphi_{p} \psi_{p}(b) , \varphi^{(i)}_{p}(\psi^{(i)}_{p}(\be_{B})) ]\| < \frac{1}{p+1}
\]
and
\begin{equation}
\label{w4-39}
\|b \varphi^{(i)}_{p}(\psi^{(i)}_{p}(\be_{B})) - \varphi^{(i)}_{p} \psi^{(i)}_{p}(b) \| < \frac{1}{p+1}
\end{equation}
for all $b \in \mathcal{G}_{p}$, $i \in \{0, \ldots,m\}$ and $p \in \mathbb{N}$.

Let
\[
h^{(j)}_{p,1},\ldots,h^{(j)}_{p,p+1} \in \mathcal{C}_{0}((0,1]), \, j=0,1,
\]
be a standard partition of unity for $(0,1]$ as in \ref{standard-partition}, and note that 
\begin{equation}
\label{w4-48}
\left\| \sum_{j=0,1;r=1,\ldots,p+1} \frac{r}{p+1} \cdot h^{(j)}_{p,r} - \id_{(0,1]} \right\| \le \frac{1}{p+1}
\end{equation}
(where $\id_{(0,1]}$ denotes the identity function on $(0,1]$) and
\begin{equation}
\label{w4-45}
\left\|  \bigoplus_{r=1}^{p+1} h^{(j)}_{p,r} (\Phi)(\be_{l}) \Phi(z) - \bigoplus_{r=1}^{p+1} \frac{r}{p+1} \cdot h^{(j)}_{p,r}(\Phi)(z)   \right\| \le \frac{1}{p+1} \cdot \|z\|
\end{equation}
for $j=0,1$ and $z \in (M_{l})_{+}$.

Use Proposition~\ref{Cnew} to find orthogonal positive contractions
\begin{equation}
\label{w4-41}
\tilde{d}^{(0)}_{p}, \tilde{d}^{(1)}_{p} \in A_{\infty} \cap B' \cap \Phi(M_{l})'
\end{equation}
such that
\begin{eqnarray}
\lefteqn{\tau\left( \tilde{d}^{(j)}_{p} \left( \bigoplus_{r=1}^{p+1} \frac{r}{p+1} \cdot h_{p,r}^{(j)}(\Phi)(z) b \right)  \right)} \nonumber \\
& \ge& \gamma_{m,\bar{m},1} \cdot \tau\left(  \bigoplus_{r=1}^{p+1} \frac{r}{p+1} \cdot h_{p,r}^{(j)}(\Phi)(z)  b \right) \label{w4-47}
\end{eqnarray}
for all $p\in \mathbb{N}$, $z \in (M_{l})_{+}$, $b \in B_{+}$, $\tau \in T_{\infty}(A)$ and $j=0,1$.

For each $p \in \mathbb{N}$, $j\in \{0,1\}$ and $r \in \{1,\ldots,p+1\}$ note that 
\[
\tilde{d}^{(j)}_{p} h^{(j)}_{p,r}(\Phi)(\, .\,)
\]
is a c.p.c.\ order zero map, so that we may apply Proposition~\ref{Cnew} once more to find pairwise orthogonal positive contractions 
\begin{equation}
\label{w4-42}
d^{(0,j)}_{p,r},\ldots,d^{(m,j)}_{p,r} \in A_{\infty} \cap B' \cap (\tilde{d}^{(j)}_{p} h^{(j)}_{p,r}(\Phi)(M_{l}))'
\end{equation}
satisfying
\begin{equation}
\label{w4-46}
\tau(d^{(i,j)}_{p,r} \tilde{d}^{(j)}_{p} h^{(j)}_{p,r}(\Phi)(z)b) \ge \gamma_{m,\bar{m},m} \cdot  \tau( \tilde{d}^{(j)}_{p} h^{(j)}_{p,r}(\Phi)(z)b)
\end{equation}
for all $i=0,\ldots,m$, $z \in (M_{l})_{+}$, $b \in B_{+}$ and $\tau \in T_{\infty}(A)$.

Use \eqref{w4-40}, \eqref{w4-10-1}, \eqref{w4-41} and \eqref{w4-42} to check that
\[
x \otimes z \mapsto \varphi^{(i)}_{p}(x) d^{(i,j)}_{p,r} \tilde{d}^{(j)}_{p} h^{(j)}_{p,r}(\Phi)(z)
\]
defines a c.p.c.\ order zero map
\[
F^{(i)}_{p} \otimes M_{l} \to A_{\infty}
\]
and apply Proposition~\ref{Dnew} to obtain c.p.c.\ order zero maps 
\begin{equation}
\label{w4-43}
\Psi^{(i,j)}_{p,r}: M_{k} \to A_{\infty} \cap (\varphi_{p}^{(i)}(F^{(i)}_{p})d^{(i,j)}_{p,r} \tilde{d}^{(j)}_{p} h^{(j)}_{p,r}(\Phi)(M_{l})  )'
\end{equation}
satisfying
\begin{eqnarray}
\lefteqn{ \tau( \varphi_{p}^{(i)}(x)d^{(i,j)}_{p,r} \tilde{d}^{(j)}_{p} h^{(j)}_{p,r}(\Phi)(z)  \Psi^{(i,j)}_{p,r} (\be_{k}))  } \nonumber  \\
& \ge & \frac{1}{4(\bar{m}+1)^{2}} \cdot  \tau(\varphi_{p}^{(i)}(x)d^{(i,j)}_{p,r} \tilde{d}^{(j)}_{p} h^{(j)}_{p,r}(\Phi)(z)) \label{w4-10-3}
\end{eqnarray}
for all $x \in (F^{(i)}_{p})_{+}$, $ z \in (M_{l})_{+}$, $i \in \{0,\ldots,m\}$, $j \in \{0,1\}$, $r \in \{1,\ldots,p+1\}$ and $\tau \in T_{\infty}(A)$.

Note that
\begin{eqnarray*}
\lefteqn{\varphi_{p}^{(i)}(\psi^{(i)}_{p}(b))  d^{(i,j)}_{p,r} \tilde{d}^{(j)}_{p} h^{(j)}_{p,r}(\Phi)(z) \Psi^{(i,j)}_{p,r}(y) }\\
& = & \Psi^{(i,j)}_{p,r}(y) h^{(j)}_{p,r}(\Phi)(z) \tilde{d}^{(j)}_{p} d^{(i,j)}_{p,r} \varphi^{(i)}_{p}(\psi^{(i)}_{p}(b))
\end{eqnarray*}
for all  $i \in \{0,\ldots,m\}$, $j \in \{0,1\}$, $r \in \{1,\ldots,p+1\}$,  $b \in B$, $ z \in M_{l}$ and $y \in M_{k}$ (using \eqref{w4-40}, \eqref{w4-10-1}, \eqref{w4-41}, \eqref{w4-42} and \eqref{w4-43}).

Define c.p.\ maps 
\[
\Psi_{p}:M_{k} \to A_{\infty}
\]
by
\begin{equation}
\label{w4-10-2}
\Psi_{p}(\, .\,) := \sum_{i,j,r} \varphi^{(i)}_{p}(\psi^{(i)}_{p}(\be_{B})) d^{(i,j)}_{p,r} \tilde{d}^{(j)}_{p} h^{(j)}_{p,r}(\Phi)(\be_{l}) \Psi^{(i,j)}_{p,r}(\, .\,).
\end{equation}
We may write $\Psi_{p}$ as a direct sum
\begin{equation}
\label{w4-44}
\Psi_{p}(\, .\,) = \bigoplus_{j=0,1}\left( \bigoplus_{i=0}^{m} \left( \bigoplus_{r=1}^{p+1} \varphi^{(i)}_{p}(\psi^{(i)}_{p}(\be_{B}))  d^{(i,j)}_{p,r}  \tilde{d}^{(j)}_{p} h^{(j)}_{p,r}(\Phi)(\be_{l}) \Psi^{(i,j)}_{p,r}(\, .\,) \right) \right)
\end{equation}
with c.p.c.\ order zero summands, which entails that $\Psi_{p}$ itself is a c.p.c.\ order zero map (for \eqref{w4-44} we have used \eqref{w4-41}, \eqref{w4-42} and the facts that
\[
\tilde{d}^{(0)}_{p} \perp \tilde{d}^{(1)}_{p},
\]
that for fixed $p$, $j$ and $r$, the $d^{(i,j)}_{p,r}$, $i=0,\ldots,m$, are pairwise orthogonal, and that for fixed $p$ and $j$ the $h^{(j)}_{p,r}$ are pairwise orthogonal).  

We estimate
\begin{eqnarray*}
\lefteqn{\tau(\Phi(\be_{l})b \Psi_{p}(\be_{k}))}\\
&\stackrel{\eqref{w4-10-2},\eqref{w4-44},\eqref{w4-39}}{\ge} & \tau\left( \Phi(\be_{l}) \sum_{i,j,r} \varphi^{(i)}_{p}(\psi^{(i)}_{p}(b)) d^{(i,j)}_{p,r} \tilde{d}^{(j)}_{p} h^{(j)}_{p,r}(\Phi) (\be_{l}) \Psi^{(i,j)}_{p,r}(\be_{k})  \right) \\
&&- \frac{1}{p+1}\\
& \stackrel{\eqref{w4-45}, \eqref{w4-42},\eqref{w4-41}, \eqref{w4-44}}{\ge} &  \tau\left( \sum_{i,j,r} \frac{r}{p+1} \cdot \varphi^{(i)}_{p}(\psi^{(i)}_{p}(b)) d^{(i,j)}_{p,r} \tilde{d}^{(j)}_{p} h^{(j)}_{p,r}(\Phi) (\be_{l}) \Psi^{(i,j)}_{p,r}(\be_{k})  \right) \\
&&- \frac{2}{p+1}\\
& \stackrel{\eqref{w4-10-3}}{\ge} & \frac{1}{4(\bar{m}+1)^{2}} \cdot \tau\left( \sum_{i,j,r} \frac{r}{p+1} \cdot \varphi^{(i)}_{p}(\psi^{(i)}_{p}(b)) d^{(i,j)}_{p,r} \tilde{d}^{(j)}_{p} h^{(j)}_{p,r}(\Phi) (\be_{l}) \right) \\
&&- \frac{2}{p+1}\\
& \stackrel{\eqref{w4-46},\eqref{w4-42}}{\ge} & \frac{\gamma_{m,\bar{m},m}}{4(\bar{m}+1)^{2}} \cdot \tau\left( \sum_{i,j,r} \frac{r}{p+1} \cdot \varphi^{(i)}_{p}(\psi^{(i)}_{p}(b))  \tilde{d}^{(j)}_{p} h^{(j)}_{p,r}(\Phi) (\be_{l}) \right) \\
&&- \frac{2}{p+1}\\
& \stackrel{\eqref{w4-47},\eqref{w4-42}}{\ge} & \frac{\gamma_{m,\bar{m},m} \cdot \gamma_{m,\bar{m},1}}{4(\bar{m}+1)^{2}} \cdot \tau\left( \sum_{i,j,r} \frac{r}{p+1} \cdot \varphi^{(i)}_{p}(\psi^{(i)}_{p}(b))   h^{(j)}_{p,r}(\Phi) (\be_{l}) \right) \\
&&- \frac{2}{p+1}\\
& \stackrel{\eqref{w4-48}}{\ge} & \frac{\gamma_{m,\bar{m},m} \cdot \gamma_{m,\bar{m},1}}{4(\bar{m}+1)^{2}} \cdot \tau\left( \sum_{i}  \varphi^{(i)}_{p}(\psi^{(i)}_{p}(b))   \Phi (\be_{l}) \right) - \frac{2+m+1}{p+1}\\
& \stackrel{\eqref{w4-49}, \eqref{w4-40}}{\ge} & \frac{\gamma_{m,\bar{m},m} \cdot \gamma_{m,\bar{m},1}}{4(\bar{m}+1)^{2}} 
\cdot \tau (    \Phi (\be_{l}) b) - \frac{4+m}{p+1} \\
& \stackrel{\eqref{w4-50}}{=} & \alpha_{m,\bar{m}} \cdot \tau (    \Phi (\be_{l}) b) - \frac{4+m}{p+1}\\
\end{eqnarray*}
for all $b \in \mathcal{G}_{p}$ and $\tau \in T_{\infty}(A)$.

Next we compute
\begin{eqnarray*}
\lefteqn{  \left\|    \Phi(z) b \Psi_{p}(y) - \sum_{i,j,r} \frac{r}{p+1} \cdot \varphi^{(i)}_{p}(\psi^{(i)}_{p}(b)) d^{(i,j)}_{p,r} \tilde{d}^{(j)}_{p} h^{(j)}_{p,r}(\Phi)(z) \Psi^{(i,j)}_{p,r}(y)    \right\| } \\
& \stackrel{\eqref{w4-44},\eqref{w4-39}}{\le} &  \left\|    \Phi(z)  \sum_{i,j,r}  \varphi^{(i)}_{p}(\psi^{(i)}_{p}(b)) d^{(i,j)}_{p,r} \tilde{d}^{(j)}_{p} h^{(j)}_{p,r}(\Phi)(\be_{l}) \Psi^{(i,j)}_{p,r}(y)  \right.  \\
& & \left. -  \sum_{i,j,r} \frac{r}{p+1} \cdot \varphi^{(i)}_{p}(\psi^{(i)}_{p}(b)) d^{(i,j)}_{p,r} \tilde{d}^{(j)}_{p} h^{(j)}_{p,r}(\Phi)(z) \Psi^{(i,j)}_{p,r}(y)\right\| \\
&&- \frac{1}{p+1}\\
& \stackrel{\eqref{w4-45},\eqref{w4-42},\eqref{w4-41},\eqref{w4-44}}{\le} & \frac{2}{p+1}
\end{eqnarray*}
for $b \in \mathcal{G}_{p}$, $y \in (M_{k})_{+}$ and $z \in (M_{l})_{+}$ of norm at most one. Taking adjoints, we see from \eqref{w4-43}, \eqref{w4-42}, \eqref{w4-41} and \eqref{w4-40} that 
\[
\|[\Phi(z)b, \Psi_{p}(y) ] \| \le \frac{4}{p+1}
\]
for $b \in \mathcal{G}_{p}$, $y \in (M_{k})_{+}$ and $z \in (M_{l})_{+}$ of norm at most one.

The usual diagonal sequence argument now yields a c.p.c.\ order zero map
\[
\check{\Psi}:M_{k} \to A_{\infty}\cap (B\Phi(M_{l}))'
\]
such that 
\[
\tau(\check{\Psi}(\be_{k}) \Phi(\be_{l})b) \ge \alpha_{m,\bar{m}} \cdot \tau(\Phi(\be_{l})b)
\]
for all $b \in B_{+}$ and $\tau \in T_{\infty}(A)$; it is not hard to check that in fact we have
\[
\check{\Psi}(M_{k}) \subset (B f(\Phi)(M_{l}))'
\]
for any $f \in \mathcal{C}_{0}((0,1])_{+}$.

For $p \in \mathbb{N}$, we may now define c.p.c.\ order zero maps 
\[
\check{\Psi}_{p}:M_{k} \to A_{\infty} \cap (B\Phi(M_{l}))'
\]
by
\[
\check{\Psi}_{p}(\, .\, ):= \check{\Psi} (\, . \,) g_{0,\frac{1}{p+1}}(\Phi(\be_{l}));
\]
one checks that 
\[
\check{\Psi}_{p} (M_{k}) \subset A_{\infty} \cap B' \cap \Phi(M_{l})'
\]
and that
\[
\tau(\check{\Psi}_{p}(\be_{k}) \Phi(\be_{l}) b) \ge \alpha_{m,\bar{m}} \cdot \tau(\Phi(\be_{l})b) - \frac{1}{p+1}
\]
for all positive contractions $b \in B_{+}$ and $\tau \in T_{\infty}(A)$. 

Another straightforward application of the diagonal sequence argument now yields a c.p.c.\ order zero map
\[
\Psi:M_{k} \to A_{\infty} \cap B' \cap \Phi(M_{l})'
\]
with the desired properties.
\end{nproof}
\en

\bn
\label{3-10prime}
We are finally prepared to prove the main technical result of this section, which will provide large almost central order zero maps into simple, unital $\mathrm{C}^{*}$-algebras with locally finite nuclear dimension and  tracial  $\bar{m}$-almost divisibility. In spirit this is very similar to \cite[Lemma~3.10]{Winter:dr-Z-stable}; the main difference is that there almost divisibility was ensured by finite decomposition rank, whereas here we make it a separate hypothesis.

\begin{nlemma}
Let $A$ be a separable, simple, unital $\mathrm{C}^{*}$-algebra with  tracial  $\bar{m}$-almost divisibility, and let $B \subset A$ be a unital $\mathrm{C}^{*}$-subalgebra with $\dimnuc B \le m$ for some $\bar{m},m \in \mathbb{N}$. Let $\mathcal{E} \subset B$ be a compact subset of  elements of norm at most one and let $k \in \mathbb{N}$ and $\gamma>0$ be given.

Then, there is a c.p.c.\ order zero map
\[
\Phi:M_{k} \to A
\]
such that
\[
\|[\Phi(x),b ]\| \le \gamma \cdot \|x\|
\]
and
\[
\tau(\Phi(\be_{k})) \ge 1 - \gamma
\]
for all $x \in M_{k}$, $b \in \mathcal{E}$ and $\tau \in T(A)$.
\end{nlemma}

\begin{nproof}
Obtain 
\[
0< \gamma_{m,\bar{m},1}<1
\]
from Proposition~\ref{Cnew} and 
\[
0<\alpha_{m,\bar{m}}<1
\]
from Proposition~\ref{3-9prime} and choose $L \in \mathbb{N}$ such that
\begin{equation}
\label{ww4-11-1}
\alpha_{m,\bar{m}} \cdot \gamma_{m,\bar{m},1} \cdot \sum_{j=0}^{L} \left(1 - \alpha_{m,\bar{m}} \cdot \gamma_{m,\bar{m},1} \right)^{j} > \left(1-\frac{\gamma}{3} \right).
\end{equation}
Set
\begin{equation}
\label{ww4-11-2}
\eta:= \frac{\gamma}{6L}.
\end{equation}
Apply Proposition~\ref{3-9prime} (with $l=1$ and $\Phi:\mathbb{C} \to A_{\infty} \cap B'$ given by $\Phi(1):= \be_{A_{\infty}}$) to obtain a c.p.c.\ order zero map 
\[
\Phi_{0}:M_{k} \to A_{\infty} \cap B'
\]
such that 
\[
\tau(\Phi_{0}(\be_{k})b) \ge \alpha_{m,\bar{m}} \cdot \tau(b) \ge \alpha_{m,\bar{m}} \cdot \gamma_{m,\bar{m},1} \cdot \tau(b)
\]
for all $b \in B_{+}$ and $\tau \in T_{\infty}(A)$. 

Now suppose that for some $i \in \{0, \ldots,L-1\}$ we have constructed a c.p.c.\ order zero map 
\[
\Phi_{i}: M_{k} \to A_{\infty} \cap B'
\]
such that
\begin{equation}
\label{w4-58}
\tau(\Phi_{i}(\be_{k})b) \ge \alpha_{m,\bar{m}} \cdot \gamma_{m,\bar{m},1} \cdot \sum_{j=0}^{i} \left( 1- \alpha_{m,\bar{m}} \cdot \gamma_{m,\bar{m},1}\right)^{j} \cdot \tau(b) - 2  i  \eta
\end{equation}
for all positive contractions $b \in B$ and for all $\tau \in T_{\infty}(A)$. 

Define c.p.c.\ order zero maps
\begin{equation}
\label{w4-57}
\hat{\Phi}_{i}:= (g_{0,\eta} - g_{\eta,2 \eta})(\Phi_{i}),\, \check{\Phi}_{i}:= (g_{\eta,2\eta} - g_{2\eta,3 \eta})(\Phi_{i}):M_{k} \to A_{\infty} \cap B'.
\end{equation}
Apply Proposition~\ref{Cnew} (with $k$ in place of $l$, $1$ in place of $k$ and $\Phi_{i}$ in place of $\Phi$) to obtain positive contractions
\begin{equation}
\label{w4-11-1}
d^{(0)} \perp d^{(1)} \in  A_{\infty}\cap B' \cap \Phi_{i}(M_{k})'
\end{equation}
such that 
\begin{equation}
\label{w4-55}
\tau(d^{(j)} f(\Phi_{i})(\be_{k})b) \ge \gamma_{m,\bar{m},1} \cdot \tau(f(\Phi_{i})(\be_{k}) b) 
\end{equation}
for all $j\in \{0,1\}$, $b \in B_{+}$, $f \in \mathcal{C}_{0}((0,1])_{+}$ and $\tau \in T_{\infty}(A)$. 

Define a c.p.c.\ order zero map 
\[
 \tilde{\Phi}_{i}:M_{k} \to A_{\infty}\cap B'
\]
by
\begin{equation}
\label{w4-54}
\tilde{\Phi}_{i}:= g_{2\eta,3\eta}(\Phi_{i}) + d^{(0)} \check{\Phi}_{i}.
\end{equation}
This is clearly c.p.c.; $\tilde{\Phi}_{i}$ has order zero since 
\[
\tilde{\Phi}_{i}(x) \le \pi_{\Phi_{i}}(x) \mbox{ for } x \in (M_{k})_{+},
\]
where $\pi_{\Phi_{i}}$ is the supporting $*$-homomorphism of $\Phi_{i}$.

Define a positive contraction
\[
h \in A_{\infty} \cap B'
\]
by
\begin{equation}
\label{w4-56}
h:= d^{(1)} \hat{\Phi}_{i}(\be_{k}) + (\be_{A_{\infty}} - g_{0,\eta}(\Phi_{i}(\be_{k}))).
\end{equation}
Observe that
\begin{equation}
\label{w4-11-2}
h \perp \tilde{\Phi}_{i}(\be_{k})
\end{equation}
by \eqref{w4-57}, \eqref{w4-11-1}, \eqref{w4-54} and \eqref{w4-56}.

Apply Proposition~\ref{3-9prime} (with $l=1$ and $\Phi:M_{1} \to A_{\infty} \cap B'$ given by $1 \mapsto h$) to obtain a c.p.c.\ order zero map 
\[
\dot{\Phi}_{i}:M_{k} \to A_{\infty} \cap B' \cap \{h\}'
\]
such that
\begin{equation}
\label{w4-52}
\tau(\dot{\Phi}_{i}(\be_{k}) hb) \ge \alpha_{m,\bar{m}} \cdot \tau(hb)
\end{equation}
for all $b \in B_{+}$ and $\tau \in T_{\infty}(A)$. 

Define a c.p.c.\ order zero map
\[
\ddot{\Phi}_{i}: M_{k} \to A_{\infty} \cap B'
\]
by
\begin{equation}
\label{w4-53}
\ddot{\Phi}_{i}(\, .\,):= h \dot{\Phi}_{i}(\, .\,).
\end{equation}
Next, define a c.p.\ map 
\[
\Phi_{i+1}:M_{k} \to A_{\infty} \cap B'
\]
by
\begin{equation}
\label{w4-51}
\Phi_{i+1}:= \ddot{\Phi}_{i} + \tilde{\Phi}_{i}
\end{equation}
and observe that $\Phi_{i+1}$ is c.p.c.\ and has order zero since $\ddot{\Phi}_{i}$ and $\tilde{\Phi}_{i}$ are c.p.c.\ order zero maps with  orthogonal images (using \eqref{w4-53} and \eqref{w4-11-2}).

For a positive contraction $b \in B_{+}$ and $\tau \in T_{\infty}(A)$ we now estimate
\begin{eqnarray*}
\lefteqn{\tau(\Phi_{i+1}(\be_{k})b)}\\
 & \stackrel{\eqref{w4-51}}{=} & \tau(\ddot{\Phi}_{i}(\be_{k})b) + \tau(\tilde{\Phi}_{i}(\be_{k})b) \\
 & \stackrel{\eqref{w4-53},\eqref{w4-52},\eqref{w4-54},\eqref{w4-55}}{\ge} & \alpha_{m,\bar{m}} \cdot \tau(hb) + \tau(g_{2\eta,3\eta}(\Phi_{i}(\be_{k}))b) \\
 && + \gamma_{m,\bar{m},1} \cdot \tau(\check{\Phi}_{i}(\be_{k})b) \\
 & \stackrel{\eqref{w4-56},\eqref{w4-55}}{\ge} & \alpha_{m,\bar{m}} \cdot \gamma_{m,\bar{m},1} \cdot \tau(\hat{\Phi}_{i}(\be_{k})b)  \\
 && + \alpha_{m,\bar{m}} \cdot \tau((\be_{A_{\infty}} - g_{0,\eta}(\Phi_{i}(\be_{k})))b) \\
 && + \tau(g_{2\eta,3\eta}(\Phi_{i}(\be_{k}))b) + \gamma_{m,\bar{m},1} \cdot \tau(\check{\Phi}_{i}(\be_{k})b) \\
& \stackrel{\eqref{w4-57}}{\ge} & \tau(g_{\eta,2\eta}(\Phi_{i}(\be_{k}))b) \\
&& +\alpha_{m,\bar{m}} \cdot \gamma_{m,\bar{m},1} \cdot \tau((\be_{A_{\infty}} - g_{\eta,2\eta}(\Phi_{i}(\be_{k})))b) \\
 & = & \alpha_{m,\bar{m}} \cdot \gamma_{m,\bar{m},1} \cdot \tau(b) \\
 && + \left( 1- \alpha_{m,\bar{m}} \cdot \gamma_{m,\bar{m},1} \right) \cdot \tau(g_{2\eta,3\eta}(\Phi_{i}(\be_{k}))b)  \\
 & \ge &  \alpha_{m,\bar{m}} \cdot \gamma_{m,\bar{m},1} \cdot \tau(b)  \\
 && + \left( 1- \alpha_{m,\bar{m}} \cdot \gamma_{m,\bar{m},1} \right) \cdot \tau(\Phi_{i}(\be_{k})b) - 2  \eta \\
 & \stackrel{\eqref{w4-58}}{\ge} &  \left( 1- \alpha_{m,\bar{m}} \cdot \gamma_{m,\bar{m},1} \right) \cdot \alpha_{m,\bar{m}} \cdot \gamma_{m,\bar{m},1} \\
 &&  \cdot \sum_{j=0}^{i} \left( 1-\alpha_{m,\bar{m}} \cdot \gamma_{m,\bar{m},1}  \right)^{j} \cdot \tau(b) \\
 && + \alpha_{m,\bar{m}} \cdot \gamma_{m,\bar{m},1} \cdot \tau(b) - 2 i  \eta - 2  \eta \\
& = &  \alpha_{m,\bar{m}} \cdot \gamma_{m,\bar{m},1}\cdot \sum_{j=0}^{i+1} \left( 1- \alpha_{m,\bar{m}} \cdot \gamma_{m,\bar{m},1}  \right)^{j} \cdot \tau(b) - 2  (i+1)  \eta.
\end{eqnarray*}
Induction yields c.p.c.\ order zero maps 
\[
\Phi_{0}, \ldots, \Phi_{L}:M_{k} \to A_{\infty} \cap B'
\]
such that 
\[
\Phi_{L}: M_{k} \to A_{\infty} \cap B'
\]
satisfies
\begin{eqnarray*}
\tau(\Phi_{L}(\be_{k})b) & \ge & \alpha_{m,\bar{m}} \cdot \gamma_{m,\bar{m},1} \cdot \sum_{j=0}^{L} (1-\alpha_{m,\bar{m}} \cdot \gamma_{m,\bar{m},1})^{j} \cdot \tau(b) -2 L \eta  \\
& \stackrel{\eqref{ww4-11-1},\eqref{ww4-11-2}}{\ge} & \left( 1- \frac{\gamma}{3} \right) \cdot \tau(b) - \frac{\gamma}{3} \ge \tau(b) - \frac{2}{3}\gamma
\end{eqnarray*}
for all positive contractions $b \in B_{+}$ and all $\tau \in T_{\infty}(A)$. 

Lift $\Phi_{L}$ to a c.p.c.\ order zero map
\[
\bar{\Phi}:M_{k} \to \prod_{\mathbb{N}}A
\]
(with c.p.c.\ order zero components $\bar{\Phi}_{\nu}, \nu \in \mathbb{N}$), then it is straightforward to show that there is $N \in \mathbb{N}$ such that 
\[
\Phi:= \bar{\Phi}_{N}
\] 
has the desired properties.
\end{nproof}
\en

\bn
\begin{nremark}
If, in the preceding result,   the elements of $\mathcal{E}$ are assumed to be positive, then the proof shows that we even obtain
\[
\tau(\Phi(\be_{k})b) \ge \tau(b) - \gamma
\]
for all  $b \in \mathcal{E}$ and $\tau \in T(A)$.
\end{nremark}
\en

\section{Almost central dimension drop embeddings}
\label{almost-central-dimension-drop-embeddings}

\noindent
The main result of  this section is Proposition~\ref{D} (which in turn refines Proposition~\ref{E}); this will provide the final missing ingredient to  prove $\mathcal{Z}$-stability in very much the same way as Proposition~4.5 did in \cite{Winter:dr-Z-stable}. In fact, Proposition~\ref{D}  is only a small modification of \cite[4.5]{Winter:dr-Z-stable}, and it is proven in almost  the same manner. The essential  difference is that in \cite{Winter:dr-Z-stable} strong tracial $m$-comparison was provided by \cite[Corollary~3.12]{Winter:dr-Z-stable}, whereas in the present setup it is built into the hypotheses. 

In the subsequent Section~{\ref{main-result}} we will then assemble our results to construct $\varphi'$ and $v'$ as in Proposition~\ref{Z-stable-relations}, which in turn gives $\mathcal{Z}$-stability.  While existence of the map $\varphi'$  essentially follows from Lemma~\ref{3-10prime}, the element $v'$ will be provided by Proposition~\ref{D}.

\bn
Let us recall Proposition~4.1 from \cite{Winter:dr-Z-stable}. 

\label{G}
\begin{nprop}
Let $A$ be a unital $\mathrm{C}^{*}$-algebra, and let $a \in A_{+}$ be a positive element of norm at most one  satisfying $\tau(a) > 0$ for any $\tau \in T(A)$. 

Then, for any $k \in \N$ and $0 < \beta < 1$ there is $\gamma >0$ such that the following holds: If  $h \in A_{+}$ is another positive element of norm at most one satisfying 
\[
\tau(h)> 1-\gamma 
\]
for any  $\tau \in T(A)$, then
\[
d_{\tau}((a^{\frac{1}{2}}(\be_{A}-h)a^{\frac{1}{2}})- \beta)_{+} < \frac{1}{k} \cdot \tau(a) 
\]
for any  $\tau \in T(A)$.
\end{nprop}
\en

\bn
The next result is essentially \cite[Proposition~4.2]{Winter:dr-Z-stable}. There, the assumption of decomposition rank at most $m$ was only used to have \cite[Corollary~3.12]{Winter:dr-Z-stable}  available. But that Corollary just says that decomposition rank at most $m$ implies strong tracial $m$-comparison, which we now make a hypothesis; the conclusions are essentially the same. In the present setup we use two dimension constants, $m$ and $\bar{m}$, to emphasize the different roles played by dimension --- these will eventually enter in terms of comparison and in terms of divisibility.

\label{E-aux}
\begin{nprop}
Let $A$ be a separable, simple, unital $\mathrm{C}^{*}$-algebra for which every quasitrace is a trace and which has strong tracial $\bar{m}$-comparison for some $\bar{m} \in \mathbb{N}$. Let $m, n \in \mathbb{N}$, let 
\[
F  = M_{r_{1}} \oplus \ldots \oplus M_{r_{s}}
\] 
be a finite-dimensional $\mathrm{C}^{*}$-algebra, and let 
\[
\psi: F \to A
\]
be a c.p.c.\ order zero map. For any $0<\bar{\gamma}_{0}$ and $0< \zeta < 1$ there are $0<\bar{\gamma}_{1}$ and $0<\bar{\beta}$ such that the following holds: 

If $ i_{0} \in \{0, \ldots,m\}$ and 
\[
\Phi: M_{m+1} \otimes M_{n} \otimes M_{2} \to A
\]
is a c.p.c.\ order zero map satisfying
\[
\| [\psi(x), \Phi(y)]\| \le \bar{\beta} \|x\| \|y\|
\]
for any   $x \in F$ and  $ y \in M_{m+1} \otimes M_{n} \otimes M_{2}$ and
\begin{equation}
\label{E-auxw1}
\tau(\Phi(\be_{M_{m+1} \otimes M_{n} \otimes M_{2}})) > 1 - \bar{\beta},
\end{equation}
then there are $\bar{v}^{(j)} \in A$, $j=1, \ldots,s$, of norm at most one with the following properties:

\begin{enumerate}
\item $\|(\bar{v}^{(j)})^{*} \bar{v}^{(j)} - \psi^{(j)}(f_{11}^{(j)})^{\halb}(\be_{A}- \Phi(\be_{M_{m+1} \otimes M_{n} \otimes M_{2}}))  \psi^{(j)}(f_{11}^{(j)})^{\halb}\| < \bar{\gamma}_{0}$
\item $\|[\bar{v}^{(j)}, \psi^{(j)}(f_{11}^{(j)})]\| < \bar{\gamma}_{0}$
\item $\|\bar{v}^{(j)} - g_{0,\bar{\gamma}_{1}}(\psi^{(j)}(f_{11}^{(j)})) \bar{v}^{(j)}\| < \bar{\gamma}_{0}$
\item $\bar{v}^{(j)} = g_{\zeta,\zeta+\bar{\gamma}_{1}}(\Phi(d_{i_{0}i_{0}}\otimes e_{11} \otimes \be_{M_{2}})) \bar{v}^{(j)}$.
\end{enumerate}
Here, $\{f_{ii'}^{(j)} \mid j=1, \ldots,s; \, i,i' =1, \ldots,r_{j}\}$ and $\{d_{ii'} \mid i,i' = 1, \ldots,m+1\}$ denote sets of matrix units for $F$ and $M_{m+1}$, respectively; the canonical matrix units for $M_{n}$ and $M_{2}$ are both denoted by $\{e_{ii'}\}$.  
\end{nprop}

\begin{nproof}
The proof is  the same as that of \cite[4.2]{Winter:dr-Z-stable} almost verbatim, with only three changes: 

In (200) of \cite{Winter:dr-Z-stable} we take $(m+1)(\bar{m}+1)$ in place of $(m+1)^{2}$, 
in (212) of \cite{Winter:dr-Z-stable} we take $\frac{1}{\bar{m}+1}$ in place of $\frac{1}{m+1}$, and instead of using \cite[Corollary~3.12]{Winter:dr-Z-stable} on page 285 in \cite{Winter:dr-Z-stable}, we use strong tracial $\bar{m}$-comparison of $A$. 
\end{nproof}
\en

\bn
The next result is a slight variation of \cite[Proposition~4.3]{Winter:dr-Z-stable}, in pretty much the same manner as Proposition~\ref{E-aux} varies \cite[Proposition~4.2]{Winter:dr-Z-stable}. 

\label{E}
\begin{nprop}
Let $A$ be a  separable, simple, unital  $\mathrm{C}^{*}$-algebra for which every quasitrace is a trace and which has strong tracial $\bar{m}$-comparison for some $\bar{m} \in \mathbb{N}$. Let $m \in \mathbb{N}$, let $F$ be a finite-dimensional $\mathrm{C}^{*}$-algebra, $n \in \N$ and
\[
\psi: F \to A
\]
a c.p.c.\ order zero map. For any $0<\theta$ and $0< \zeta < 1$ there is $\beta >0$ such that the following holds: 

If $ i_{0} \in \{0, \ldots,m\}$ and 
\[
\Phi: M_{m+1} \otimes M_{n} \otimes M_{2} \to A
\]
is a c.p.c.\ order zero map satisfying
\begin{equation}
\label{c}
\| [\psi(x), \Phi(y)]\| \le \beta \|x\| \|y\| 
\end{equation}
for any   $x \in F$  and  $y  \in M_{m+1} \otimes M_{n} \otimes M_{2}$ and
\[
\tau(\Phi(\be_{M_{m+1} \otimes M_{n} \otimes M_{2}})) > 1 - \beta \; \forall \, \tau \in T(A),
\]
then there is $v \in A$ of norm at most one such that
\[
\|v^{*}v - (\be_{A}- \Phi(\be_{M_{m+1} \otimes M_{n} \otimes M_{2}}))^{\frac{1}{2}} \psi(\be_{F})  (\be_{A}- \Phi(\be_{M_{m+1} \otimes M_{n} \otimes M_{2}}))^{\frac{1}{2}}\| < \theta,
\]
\[
v v^{*} \in \overline{(\Phi(d_{i_{0}i_{0}} \otimes e_{11} \otimes \be_{M_{2}})- \zeta)_{+}A(\Phi(d_{i_{0}i_{0}} \otimes e_{11} \otimes \be_{M_{2}})- \zeta)_{+}}
\]
(where $\{d_{kl} \mid k,l=0, \ldots,m\}$  and $\{e_{kl} \mid k,l=1, \ldots, n\}$ denote  sets of matrix units for $M_{m+1}$ and $M_{n}$, respectively) and
\[
\|[\psi(x),v]\| \le \theta \|x\|
\]
for all   $x \in F$.
\end{nprop}

\begin{nproof}
The proof is essentially that of \cite[Proposition~4.3]{Winter:dr-Z-stable}. The main difference is that  Proposition~\ref{E-aux} is used instead of \cite[Proposition 4.2]{Winter:dr-Z-stable}, but I mention another necessary adjustment, pointed out to me by D.\ Archey: The elements $v^{(j)}$ and $v$ defined in \cite[equations~(241) and (242)]{Winter:dr-Z-stable} are not necessarily in $A$; they should be replaced by 
\[
\tilde{v}^{(j)}:= \sum_{k=1}^{r^{(j)}} \pi^{(j)}(f^{(j)}_{k1}) g_{0,\bar{\gamma}_1}(\psi^{(j)}(f_{11}^{(j)}))\bar{v}^{(j)} g_{0,\bar{\gamma}_1}(\psi^{(j)}(f_{11}^{(j)})) \pi^{(j)}(f^{(j)}_{1k}  )
\] 
and
\[
\tilde{v}:= g_{\zeta,\zeta+\bar{\gamma}_{1}}(\Phi(d_{i_{0}i_{0}} \otimes e_{11} \otimes \be_{M_{2}})) \sum_{j=1}^{s} \tilde{v}^{(j)}. 
\]
Then, the $\tilde{v}^{(j)}$ and $\tilde{v}$ clearly are in $A$, and one checks that $\|v^{(j)} - \tilde{v}^{(j)}\| < 2 \bar{\gamma}_0^{1/2}$ and $\|v-\tilde{v}\| < 2s \bar{\gamma}_{0}^{\frac{1}{2}}$ (using \cite[4.2(i) and 4.2(iii)]{Winter:dr-Z-stable}). The rest of the proof works in the same manner, provided one chooses $0<\bar{\gamma}_{0}<1$ such that $16 s \max\{r^{(j)}\} \bar{\gamma}_{0}^{\frac{1}{2}} < \theta$ (cf.\ \cite[inequality~(236)]{Winter:dr-Z-stable}).
\end{nproof}
\en

\bn
\label{H}
Below we essentially repeat \cite[Proposition~4.4]{Winter:dr-Z-stable}, with the additional observation that there it would not have been quite necessary to assume all of $A$ (as opposed to a $\mathrm{C}^{*}$-subalgebra $B$) to have finite decomposition rank, and that one can replace decomposition rank by nuclear dimension with only a little bit of extra effort.  

\begin{nprop}
Let $A$ be a simple, separable, unital $\mathrm{C}^{*}$-algebra, and let $B \subset A$ be a unital $\mathrm{C}^{*}$-subalgebra with $\dimnuc B \le m < \infty$. Given a finite subset $\Fh \subset B$, a positive normalized function $\bar{h} \in \mathcal{C}_{0}((0,1])$ and $\delta>0$, there are a finite subset $\Gh \subset B$ and $\alpha > 0$ such that the following holds:

Suppose 
\[
(F=F^{(0)} \oplus \ldots \oplus F^{(m)},\sigma,\varrho)
\]
is an $m$-decomposable c.p.\ approximation (for $B$) of $\Gh$ to within  $\alpha$, with c.p.c.\ maps $\sigma:B \to F$ and $\varrho^{(i)}:F^{(i)}\to B$, and such that the composition $\varrho \sigma$ is contractive. For $i=0, \ldots, m$, let $v_{i} \in A$ be normalized elements satisfying
\begin{equation}
\label{w5}
\|[\varrho^{(i)}(x),v_{i}]\| \le \alpha \|x\| 
\end{equation}
for all  $x \in F^{(i)}$. 

Then, 
\begin{equation}
\label{wH5}
v:= \sum_{i=0}^{m} v_{i} \bar{h}(\varrho^{(i)})(\sigma^{(i)}(\be_{A}))
\end{equation}
satisfies
\[
\|[v,a]\| < \delta
\]
for all  $a \in \Fh$.
\end{nprop}

\begin{nproof}
The proof is essentially the same as that of  \cite[Proposition~4.4]{Winter:dr-Z-stable}, but since there are several  modifications, we revisit the full proof in detail.

We may assume the elements of $\mathcal{F}$ to be positive and normalized. For convenience, we set
\begin{equation}
\label{wdH4}
\bar{\delta}:= \frac{\delta}{9(m+1)}.
\end{equation}
Take $\tilde{h} \in \Ch([0,1])$ such that
\begin{equation}
\label{wd2}
\| \id_{[0,1]} \cdot \tilde{h} - \bar{h}\| < \bar{\delta}.
\end{equation}
By using \cite[Proposition~1.8]{Winter:dr-Z-stable}, there is 
\begin{equation}
\label{w5-4-1}
0< \beta< \frac{\bar{\delta}}{\|\tilde{h}\|}
\end{equation}
such that, whenever $b,c \in A$ are elements of norm at most one with $c$ positive, and satisfying
\begin{equation}
\label{wdH1}
\|[b,c]\| \le \beta,
\end{equation}
we have 
\begin{equation}
\label{wdH2}
 \|[b,\tilde{h}(c)]\| \le \bar{\delta}.
\end{equation}
Using \cite[Lemma~1.10]{Winter:dr-Z-stable}, we find $\alpha>0$ and a finite subset $\mathcal{G} \subset B$ such that, whenever $(F,\hat{\sigma},\hat{\varrho})$ is a c.p.c.\ approximation (for $B$) of $\mathcal{G}$ to within  $\alpha$, we have 
\begin{equation}
\label{wd4}
\|\hat{\varrho}(x) \hat{\varrho} \hat{\sigma}(a) - \hat{\varrho}(x \hat{\sigma}(a))\| \le \frac{\beta}{2} \|x\| 
\end{equation}
for all $x \in F$ and $a \in \mathcal{F}$. We may assume that 
\begin{equation}
\label{wdH3}
\alpha < \beta, \bar{\delta},  \frac{\bar{\delta}}{\|\tilde{h}\|}. 
\end{equation}
Now let a c.p.\ approximation $(F,\sigma,\varrho)$ (for $B$) and $v_{i} \in A$ as in the proposition be given. Define maps
\[
\hat{\sigma}(\, . \,):=\sigma(\be_{A})^{-1/2} \sigma(\, . \,) \sigma(\be_{A})^{-1/2}
\]
and
\[
\hat{\varrho}(\, .\,):=\varrho(\sigma(\be_{A})^{1/2} \, . \, \sigma(\be_{A})^{1/2}) 
\]
(with inverses taken in the respective hereditary subalgebra), then 
\begin{equation}
\label{wdd1}
\hat{\varrho}^{(i)}\hat{\sigma}^{(i)} = \varrho^{(i)} \sigma^{(i)}
\end{equation}
for all $i$ and $(F,\hat{\sigma},\hat{\varrho})$ is a c.p.c.\ approximation of $\mathcal{G}$ to within $\alpha$. 

Note that  we have
\begin{eqnarray}
\lefteqn{\| \bar{h}(\varrho^{(i)})(\sigma^{(i)}(\be_{A})) - \varrho^{(i)}\sigma^{(i)} (\be_{A}) \tilde{h}(\varrho^{(i)}(\be_{F^{(i)}}))  \| } \nonumber\\ 
& \stackrel{\ref{order-zero-notation}}{=} &  \|  \pi_{\varrho^{(i)}} (\sigma^{(i)}(\be_{A})) \bar{h}(\varrho^{(i)}(\be_{F^{(i)}}))  - \varrho^{(i)}\sigma^{(i)} (\be_{A}) \tilde{h}(\varrho^{(i)}(\be_{F^{(i)}})) \| \nonumber \\
& \stackrel{\eqref{wd2}}{<}&   \|  \pi_{\varrho^{(i)}} (\sigma^{(i)}(\be_{A}))   \varrho^{(i)}(\be_{F^{(i)}})   \tilde{h}(\varrho^{(i)}(\be_{F^{(i)}})) \nonumber \\
&& - \varrho^{(i)}\sigma^{(i)} (\be_{A}) \tilde{h}(\varrho^{(i)}(\be_{F^{(i)}})) \|  + \bar{\delta}\nonumber \\
&  \stackrel{\ref{order-zero-notation}}{=}&  \bar{\delta}     \label{w5-4-*}
\end{eqnarray}
%
and 
\begin{equation}
\label{wd3}
\| \varrho^{(i)} \sigma^{(i)}(\be_{A}) \varrho \sigma (a) - \varrho^{(i)} \sigma^{(i)}(a) \| \stackrel{\eqref{wd4}, \eqref{wdH3},\eqref{wdd1},\eqref{w5-4-1}}{\le} \frac{\bar{\delta}}{\|\tilde{h}\|}
\end{equation}
for $i = 0, \ldots, m$ and $a \in \mathcal{F}$ (the elements of $\mathcal{F}$ are normalized). We obtain
\begin{eqnarray*}
\lefteqn{ \| [v_{i} \bar{h}(\varrho^{(i)}) (\sigma^{(i)}(\be_{A})),a  ] \| } \\
& \stackrel{\eqref{w5-4-*}}{\le} & \| [v_{i} \varrho^{(i)} \sigma^{(i)} (\be_{A}) \tilde{h}(\varrho^{(i)} (\be_{F^{(i)}})), \varrho \sigma(a) ] \| + 2 \alpha +\bar{\delta} \\
& \stackrel{\eqref{w5},\eqref{wdH1},\eqref{wdH2},\eqref{wdH3}}{\le} & \| v_{i} \varrho^{(i)}\sigma^{(i)} (\be_{A}) \tilde{h}(\varrho^{(i)}(\be_{F^{(i)}})) \varrho \sigma(a)  \\
&& - \varrho \sigma (a) \tilde{h}(\varrho^{(i)}(\be_{F^{(i)}})) \varrho^{(i)} \sigma^{(i)} (\be_{A}) v_{i} \| \\
&& + 2 \alpha + 3 \bar{\delta}  \\
& \stackrel{\eqref{wd3}}{\le} & \| v_{i} \tilde{h} (\varrho^{(i)}(\be_{F^{(i)}})) \varrho^{(i)} \sigma^{(i)}(a)   \\
&& - \varrho^{(i)} \sigma^{(i)}(a) \tilde{h}(\varrho^{(i)} (\be_{F^{(i)}})) v_{i} \| \\
&& + 2 \alpha + 5 \bar{\delta}   \\ 
& \stackrel{\eqref{w5},\eqref{wdH1}, \eqref{wdH2}, \eqref{wdH3}, \eqref{w5-4-1}}{<} & 9 \bar{\delta} \\
& \stackrel{\eqref{wdH4}}{\le} &  \frac{\delta}{m+1} 
\end{eqnarray*}
for $i = 0, \ldots, m$ and $ a \in \mathcal{F}$.  It follows that
\[
\|[v,a]\| \stackrel{\eqref{wH5}}{=} \left\|\left[\sum_{i=0}^{m} v_{i} \bar{h}(\varrho^{(i)})(\sigma^{(i)}(\be_{A})  ),a\right]\right\| < (m+1) \frac{\delta}{m+1} = \delta
\]
for $a \in \mathcal{F}$.
\end{nproof}
\en

\bn
\label{D}
We are now prepared to state the main technical result of this section, which again is a modification of a result from \cite{Winter:dr-Z-stable}, Proposition~4.5.

\begin{nprop}
Let $A$ be a separable, simple, unital $\mathrm{C}^{*}$-algebra with strong tracial $\bar{m}$-comparison for some $\bar{m} \in \mathbb{N}$ and such that every quasitrace on $A$ is a trace.  Let $B \subset A$ be a unital $\mathrm{C}^{*}$-subalgebra with $\dimnuc B \le m < \infty$. Given a finite subset $\Fh \subset B$, $0<\delta$, $0<\zeta<1$ and $n \in \N$, there are $\Gh \subset B$ finite and $\alpha>0$ such that the following holds: 
 
 Whenever $(F,\sigma,\varrho)$ is an $m$-decomposable c.p.\ approximation (for $B$) of $\Gh$ to within  $\alpha$, with c.p.c.\ maps $\sigma:B \to F$ and $\varrho^{(i)}:F^{(i)} \to B$, and such that the composition $\varrho \sigma$ is contractive, there is $\gamma >0$ satisfying the following:  If 
\[
\Phi: M_{m+1} \otimes M_{n} \otimes M_{2} \to A
\]
is a c.p.c.\ order zero map satisfying 
\begin{equation}
\label{wwD7}
\|[\varrho(x), \Phi(y)]\| \le \gamma \|x\| \|y\| 
\end{equation}
for all  $ x \in F, \,  y \in M_{m+1} \otimes M_{n} \otimes M_{2}$ and 
\begin{equation}
\label{wwD12}
\tau( \Phi(\be_{M_{m+1} \otimes M_{n} \otimes M_{2}})) \ge 1 - \gamma 
\end{equation}
for all  $\tau \in T(A)$, then  
\[
\|[\Phi(y),a]\| \le \delta \|y\| 
\]
for all  $ y \in M_{m+1} \otimes M_{n} \otimes M_{2}, \, a \in \Fh$, and there is $v \in A$ such that
\[
\|v^{*} v - (\be_{A} - \Phi(\be_{M_{m+1} \otimes M_{n} \otimes M_{2}}))\| < \delta,
\]
\[
vv^{*} \in \overline{(\Phi(\be_{M_{m+1}} \otimes e_{11} \otimes \be_{M_{2}})-\zeta)_{+} A (\Phi(\be_{M_{m+1}} \otimes e_{11} \otimes \be_{M_{2}})-\zeta)_{+}}
\]
and
\[
\|[a,v]\| < \delta  
\]
for all  $a \in \Fh$.
\end{nprop}

\begin{nproof}
Again, the proof is essentially the same as that of \cite[Proposition~4.5]{Winter:dr-Z-stable}, but since there are some modifications and to avoid confusion, we state the full proof explicitly.

We may clearly assume that the elements of $\Fh$ are normalized. Set 
\begin{equation}
\label{wwD9}
\bar{h}:= g_{0, \delta/(6(m+1))}
\end{equation}
and take $\tilde{h} \in \mathcal{C}([0,1])$ such that
\begin{equation}
\label{w5-5-2}
\id_{[0,1]} \cdot \tilde{h} = \bar{h};
\end{equation}
note that
\[
\|\tilde{h}\| = \frac{6(m+1)}{\delta}.
\]
From Proposition~\ref{H}, obtain a finite subset $\Gh \subset B$ and $\alpha>0$ such that the assertion of \ref{H} holds. We may assume that 
\begin{equation}
\label{wwD10}
\alpha < \frac{\delta}{6}
\end{equation} 
and that 
\begin{equation}
\label{wwD6}
\be_{A} \in \mathcal{G}.
\end{equation}  

Now suppose 
\begin{equation}
\label{wwD11}
(F,\sigma,\varrho)
\end{equation} 
is a c.p.\ approximation (for $B$) of $\Gh$ to within  $\alpha$ as in the proposition, such that $\varrho$ is $m$-decomposable with respect to $F=F^{(0)} \oplus \ldots \oplus F^{(m)}$.

Similarly as in \eqref{w5-4-*}, observe that
\begin{equation}
\label{w5-5-3}
\bar{h}(\varrho^{(i)})(\sigma^{(i)}(\be_{A})) = \varrho^{(i)} \sigma^{(i)} (\be_{A}) \tilde{h}(\varrho^{(i)}(\be_{F^{(i)}})).
\end{equation}

Choose some 
\begin{equation}
\label{2}
0< \theta<\frac{\alpha}{3(m+1)}. 
\end{equation}
Obtain $\beta >0$ from Proposition~\ref{E} so that the assertion of \ref{E} holds for each 
\[
\varrho^{(i)}:= \varrho|_{F^{(i)}},
\]
$i=0,\ldots,m$, in place of $\psi$. Using \cite[Proposition~1.8]{Winter:dr-Z-stable}, we may choose some 
\begin{equation}
\label{wwD1}
0< \gamma < \min\{\theta,\beta\}
\end{equation}
such that, if $b$ and $c$ are elements of norm at most one in some $\mathrm{C}^{*}$-algebra, with $c \ge 0$ and 
\[
\|[b,c]\|< \gamma,
\]
then 
\begin{equation}
\label{wwD8}
\|[b^{\frac{1}{2}},c ]\| , \|[b^{\frac{1}{2}},\tilde{h}(c)]\|, \|[b^{\frac{1}{2}},c^{\frac{1}{2}} ]\| < \theta \cdot \frac{\delta}{6(m+1)}.
\end{equation} 

Now if $\Phi$ is as in the assertion of Proposition~\ref{D}, then (note \eqref{wwD1}) Proposition~\ref{E} yields  elements $v_{i} \in A$ of norm at most one, $i=0,\ldots,m$, such that
\begin{eqnarray}
\|v_{i}^{*} v_{i} - (\be_{A} - \Phi(\be_{M_{m+1} \otimes M_{n} \otimes M_{2}}))^{\frac{1}{2}} \varrho^{(i)}(\be_{F^{(i)}}) (\be_{A} - \Phi(\be_{M_{m+1} \otimes M_{n} \otimes M_{2}}))^{\frac{1}{2}} \| && \nonumber \\
<  \theta, && \label{wwD5}
\end{eqnarray}
\begin{equation}
\label{wwD2}
v_{i}v_{i}^{*} \in \overline{(\Phi(d_{ii} \otimes e_{11} \otimes \be_{M_{2}}) - \zeta)_{+} A (\Phi(d_{ii} \otimes e_{11} \otimes \be_{M_{2}}) - \zeta)_{+}}
\end{equation}
and
\begin{equation}
\label{1}
\|[\varrho^{(i)}(x),v_{i}]\| \le \theta \|x\| \mbox{ for all }  x \in F^{(i)}
\end{equation}
for each $i$.  Note that 
\begin{equation}
\label{wwD3}
v_{i} v_{i}^{*} \perp v_{i'} v_{i'}^{*}
\end{equation}
by \eqref{wwD2} if $i \neq i'$, since $\Phi$ has order zero.   Set
\begin{equation}
\label{wwD4}
v:= \sum_{i=0}^{m} v_{i} \bar{h}(\varrho^{(i)})(\sigma^{(i)}(\be_{A})),
\end{equation}
then
\begin{eqnarray*}
\lefteqn{\|v^{*}v - (\be_{A}- \Phi(\be_{M_{m+1} \otimes M_{n} \otimes M_{2}})) \|}\\
& \stackrel{\eqref{wwD4},\eqref{wwD3}}{=} & \left\| \sum_{i=0}^{m} \bar{h}(\varrho^{(i)})(\sigma^{(i)}(\be_{A}))  v_{i}^{*} v_{i}  \bar{h}(\varrho^{(i)})(\sigma^{(i)}(\be_{A})) \right. \\
&& \left. -  (\be_{A}- \Phi(\be_{M_{m+1} \otimes M_{n} \otimes M_{2}}))  \right\| \\
& \stackrel{\eqref{wwD5},\eqref{wwD6}}{\le} & \left\| \sum_{i=0}^{m} \bar{h}(\varrho^{(i)})(\sigma^{(i)}(\be_{A}))   (\be_{A}- \Phi(\be_{M_{m+1} \otimes M_{n} \otimes M_{2}}))^{\frac{1}{2}} \varrho^{(i)}(\be_{F^{(i)}}) \right. \\
&& \cdot (\be_{A}- \Phi(\be_{M_{m+1} \otimes M_{n} \otimes M_{2}}))^{\frac{1}{2}} \bar{h}(\varrho^{(i)})(\sigma^{(i)}(\be_{A})) \\
&&  -  (\be_{A}- \Phi(\be_{M_{m+1} \otimes M_{n} \otimes M_{2}}))^{\frac{1}{2}}  \varrho^{(i)}\sigma^{(i)}(\be_{A}) \\
&& \left. (\be_{A}- \Phi(\be_{M_{m+1} \otimes M_{n} \otimes M_{2}}))^{\frac{1}{2}}  \right\| \\
&&  + (m+1) \theta + \alpha \\
& \stackrel{\eqref{wwD7},\eqref{wwD8}, \eqref{w5-5-3}}{\le} & \left\| \sum_{i=0}^{m} (\be_{A}- \Phi(\be_{M_{m+1} \otimes M_{n} \otimes M_{2}}))^{\frac{1}{2}}  \right. \\
&& (\bar{h}(\varrho^{(i)})(\sigma^{(i)}(\be_{A}))   \varrho^{(i)}(\be_{F^{(i)}})    \bar{h}(\varrho^{(i)})(\sigma^{(i)}(\be_{A})) -    \varrho^{(i)}\sigma^{(i)}(\be_{A})   )   \\
&& \left.    (\be_{A}- \Phi(\be_{M_{m+1} \otimes M_{n} \otimes M_{2}}))^{\frac{1}{2}} \right\| \\
&& + (m+1) \theta + \alpha + 4 (m+1) \theta \\
&\stackrel{\eqref{w5-5-3},\eqref{w5-5-2}\eqref{wwD9}}{\le} & + (m+1) \theta + \alpha + 4 (m+1) \theta + 2(m+1) \frac{\delta}{6(m+1)} \\
& \stackrel{\eqref{2},\eqref{wwD10}}{<} & \delta.
\end{eqnarray*}
Moreover,
\begin{eqnarray*}
vv^{*} & \stackrel{\eqref{wwD4},\eqref{wwD2}}{\in} & \her \left(\sum_{i=0}^{m} (\Phi(d_{ii} \otimes e_{11} \otimes \be_{M_{2}})-\zeta)_{+}\right)  \\
&\stackrel{\ref{order-zero-facts}}{=}& \her ((\Phi(\be_{M_{m+1}} \otimes e_{11} \otimes \be_{M_{2}})- \zeta)_{+}).
\end{eqnarray*}
We have 
\[
\|[v,a]\|  < \delta \mbox{ for all } a \in \Fh
\]
by Proposition~\ref{H} and our choice of $\alpha$, using \eqref{1} and \eqref{2}. Finally, for $y \in M_{m+1} \otimes M_{n} \otimes M_{2}$ and $a \in \Fh$ we have
\begin{eqnarray*}
\|[\Phi(y),a]\| & \stackrel{\eqref{wwD11}}{\le}  & \|[\Phi(y), \varrho \sigma (a)] \| + 2 \alpha \|y\| \\
& \stackrel{\eqref{wwD7}}{\le} & \gamma \|y\| + 2 \alpha \|y\| \\
& \stackrel{\eqref{wwD1},\eqref{2},\eqref{wwD10}}{\le} & \delta \|y\|.
\end{eqnarray*}
\end{nproof}
\en

\section{The main result and its consequences}
\label{main-result}

\noindent
We are now in a position to prove the main result and its corollaries. As was already mentioned, the structure of the proof of Theorem~\ref{A} is similar to that of \cite[Theorem~5.1]{Winter:dr-Z-stable}.

\bn
\label{A}
\begin{ntheorem}
Let $A$ be a separable, simple, unital $\mathrm{C}^{*}$-algebra with locally finite nuclear dimension; suppose $A$ is tracially $(\bar{m},\tilde{m})$-pure for some $\bar{m}, \tilde{m} \in \mathbb{N}$. 

Then, $A$ is $\Zh$-stable.
\end{ntheorem}

\begin{nproof}
We check that $A$ satisfies the hypotheses of Proposition~\ref{Z-stable-relations}. This follows very much the same pattern as the proof of \cite[Theorem~5.1]{Winter:dr-Z-stable}. 

So, let $n \in \N$, $\Fh \subset A$ finite and $0<\delta<1$ and $0<\zeta<1$ be given. We have to find $\varphi':M_{n} \to A$ and $v' \in A$ as in \ref{Z-stable-relations}.

Since $A$ has locally finite nuclear dimension, we may assume that $\mathcal{F} \subset B$, where $B \subset A$ is a unital $\mathrm{C}^{*}$-subalgebra with $\dimnuc B = m <\infty$. Moreover, by Remark~\ref{quasitraces-traces} every quasitrace on $A$ is a trace.

By Proposition~\ref{D}, there are a finite subset $\be_{A} \in \Gh \subset B$ and $\alpha>0$, such that, if  $(F, \sigma, \varrho)$ is an $m$-decomposable c.p.\ approximation (for $B$) of $\Gh$ to within  $\alpha$ as in \ref{D}, there is $\gamma >0$ such that the following holds: If
\[
\Phi: M_{m+1} \otimes M_{n} \otimes M_{2} \to A
\]
is a c.p.c.\ order zero map satisfying \eqref{wwD7} and \eqref{wwD12} of \ref{D}, then there is $v' \in A$ which, together with 
\[
\varphi':= \Phi|_{\be_{M_{m+1}}\otimes M_{n} \otimes \be_{M_{2}}},
\]
satisfies the hypotheses of \ref{Z-stable-relations}. 

So, choose an $m$-decomposable c.p.\  approximation $(F, \sigma, \varrho)$ (for $B$) of $\Gh$ to within  $\alpha$ as in Proposition~\ref{D}.   
The existence of $\Phi$ now follows from Lemma~\ref{3-10prime} with $(m+1)\cdot n \cdot 2$ in place of $k$ and $\Eh:= \varrho(\Bh_{1}(F))$ (where $\Bh_{1}(F)$ denotes the unit ball of $F$). 
\end{nproof}
\en

\bn
\label{Acor}
Let us note the following combination of Corollary~\ref{pure-tracially-pure-cor} and Theorem~\ref{A} explicitly.

\begin{ncor}
Let $A$ be a separable, simple, unital $\mathrm{C}^{*}$-algebra with locally finite nuclear dimension. Suppose $A$ is $(m,\bar{m})$-pure for some $m,\bar{m} \in \mathbb{N}$. (For example, $A$ could be pure.)

Then, $A$ is $\Zh$-stable.
\end{ncor}
\en

\bn
\label{findimnuc-Z-stable}
\begin{ncor}
Let $A$ be a separable, simple, nonelementary, unital $\mathrm{C}^{*}$-algebra with  finite nuclear dimension. 

Then, $A$ is $\Zh$-stable.
\end{ncor}

\begin{nproof}
This combines Corollary~\ref{dimnuc-tracially-pure} and  Theorem~\ref{A}.
\end{nproof}
\en

\bn
\label{W-Z-stable}
\begin{ncor}
Let $A$ be a separable, simple, nonelementary, unital $\mathrm{C}^{*}$-algebra with locally finite nuclear dimension. 

Then, $A \cong A \otimes \mathcal{Z}$ if and only if $W(A) \cong W(A \otimes \mathcal{Z})$.
\end{ncor}

\begin{nproof}
If $W(A) \cong W(A \otimes \mathcal{Z})$, then $A$  is pure  by Proposition~\ref{W-Z-purely-finite}, so $A$ is $\mathcal{Z}$-stable  by Corollary~\ref{Acor}. The forward implication is trivial.
\end{nproof}
\en

\bn
\label{ASH-sdg-Z}
\begin{ncor}
Let $A$ be a separable, simple, nonelementary, unital, approximately subhomogeneous $\mathrm{C}^{*}$-algebra. 

Then, $A$ is $\mathcal{Z}$-stable if and only if it has slow dimension growth (in the sense of \cite{Toms:ash-sdg}). 
\end{ncor}

\begin{nproof}
If $A$ has slow dimension growth, we have $W(A) \cong W(A \otimes \mathcal{Z})$ by \cite[Theorem~1.2]{Toms:ash-sdg},  so $A$ is $\mathcal{Z}$-stable by the preceding Corollary. ($A$ has locally finite decomposition rank, hence locally finite nuclear dimension,  see \cite{NgWinter:subhom}.) 

Conversely, if $A$ is $\mathcal{Z}$-stable, it has slow dimension growth by \cite{Toms:ash-sdg}.
\end{nproof}
\en

\bn
\label{ASH-sdg-classification}
\begin{ncor}
Let $\mathcal{A}$ denote the class of  separable, simple, nonelementary, unital, approximately subhomogeneous $\mathrm{C}^{*}$-algebras with slow dimension growth and for which projections separate traces. 

Then, $\mathcal{A}$ is classified by ordered $\mathrm{K}$-theory. 
\end{ncor}

\begin{nproof}
The members of  $\mathcal{A}$ are $\mathcal{Z}$-stable by Corollary~\ref{ASH-sdg-Z}, so they are classified by \cite[Corollary~5.5]{LinNiu:KKlifting} (cf.\ also \cite[Corollary~8.3]{Win:localizingEC}).
\end{nproof}
\en

\bn
\label{AH-sdg-classification}
Let us note the following corollary explicitly. This is mostly repeating known results;  the only new implication is (iii)$\Rightarrow$(iv), which follows from Corollary~\ref{ASH-sdg-Z} and settles an important technical question left open in \cite{EllGongLi:simple_AH}.

\begin{ncor}
Let $A$ be a separable, simple, nonelementary, unital, approximately homogeneous $\mathrm{C}^{*}$-algebra. 

Then, the following are equivalent:
\begin{itemize}
\item[(i)] $A$ has no dimension growth (as an $\mathrm{AH}$ algebra)
\item[(ii)] $A$ has very slow dimension growth (as an $\mathrm{AH}$ algebra)
\item[(iii)] $A$ has slow dimension growth (as an $\mathrm{AH}$ algebra)
\item[(iv)] $A$ is $\mathcal{Z}$-stable
\item[(v)] $A$ has finite decomposition rank.
\end{itemize}
The class of such $\mathrm{C}^{*}$-algebras is classified by the Elliott invariant.
\end{ncor}

\begin{nproof}
Implications (i)$\Rightarrow$(ii)$\Rightarrow$(iii)  and (i)$\Rightarrow$(v) are trivial; (iii)$\Rightarrow$(iv) follows from Corollary~\ref{ASH-sdg-Z}, and (v)$\Rightarrow$(iv) is (a special case of) \cite[Theorem~5.1]{Winter:dr-Z-stable}. The classification result for $\mathcal{Z}$-stable $\mathrm{AH}$ algebras was obtained in \cite{Lin:asu-class}; this also yields (iv)$\Rightarrow$(i), cf.\  \cite[Corollary~11.13]{Lin:asu-class}.
\end{nproof}
\en



\begin{thebibliography}{10}

\bibitem{Bla:encyc}
B.~Blackadar, \emph{{Operator Algebras}}, vol. 122, Encyclopaedia of
  Mathematical Sciences. Subseries: Operator Algebras and Non-commutative
  Geometry, no. III, Springer Verlag, Berlin, Heidelberg, 2006.

\bibitem{BlaHan:quasitrace}
B.~Blackadar and D.~Handelman, \emph{{Dimension functions and traces on
  $C^*$-algebras}}, J. Funct. Anal. \textbf{45} (1982), 297--340.

\bibitem{BRTTW:Cu}
B.~Blackadar, L.~Robert, A.~Tikuisis, A.~S. Toms, and W.~Winter, \emph{{An
  algebraic approach to the radius of comparison}}, arXiv preprint
  math.OA/1008.4024; to appear in Trans. Amer. Math. Soc., 2010.

\bibitem{BroPed:realrank}
L.~G. Brown and G.~K. Pedersen, \emph{{$C^*$-algebras of real rank zero}}, J.
  Funct. Anal. \textbf{99} (1991), 131--149.

\bibitem{BroPerToms:cuntz-semigroup}
N.~Brown, F.~Perera, and A.~S. Toms, \emph{{The Cuntz semigroup, the Elliott
  conjecture, and dimension functions on $C^{*}$-algebras}}, J. Reine Angew.
  Math. \textbf{621} (2008), 191--211.

\bibitem{BroWinter:dimnuc-quasitraces}
N.~P. Brown and W.~Winter, \emph{{Quasitraces are traces: a simple proof in the
  finite-nuclear-dimension case}}, arXiv preprint math.OA/1005.2229; to appear
  in C. R. Acad. Sci. Canada, 2010.

\bibitem{CEI:Cu}
K.~Coward, G.~A. Elliott, and C.~Ivanescu, \emph{{The Cuntz semigroup as an
  invariant for $C^*$-algebras}}, J. Reine Angew. Math. \textbf{623} (2008),
  161--193.

\bibitem{Cuntz:On}
J.~Cuntz, \emph{{Simple $C^*$-algebras generated by isometries}}, Comm. Math.
  Phys. \textbf{57} (1977), 173--185.

\bibitem{Cuntz:dimension}
\bysame, \emph{{Dimension functions on simple $C^*$-algebras}}, Math. Ann.
  \textbf{233} (1978), 145--153.

\bibitem{DadToms:Z}
M.~D\u{a}d\u{a}rlat and A.~Toms, \emph{{A universal property for the Jiang--Su
  algebra}}, Adv. Math. \textbf{220} (2009), 341--366.

\bibitem{EllGongLi:simple_AH}
G.~A. Elliott, G.~Gong, and L.~Li, \emph{{On the classification of simple
  inductive limit $C^*$-algebras II: The isomorphism theorem}}, Invent. Math.
  \textbf{168} (2007), no.~2, 249--320.

\bibitem{EllToms:BullAMS}
G.~A. Elliott and A.~S. Toms, \emph{{Regularity properties in the
  classification program for separable amenable {$C\sp *$}-algebras}}, Bull.
  Amer. Math. Soc. (N.S.) \textbf{45} (2008), no.~2, 229--245.

\bibitem{Haa:quasi}
U.~Haagerup, \emph{{Every quasi-trace on an exact $C^*$-algebra is a trace}},
  unpublished manuscript, 1991.

\bibitem{HirWinZac:Rokhlin-dimension}
I.~Hirshberg, W.~Winter, and J.~Zacharias, \emph{{The Rokhlin dimension of
  dynamical systems}}, in preparation, 2010.

\bibitem{JiaSu:Z}
X.~Jiang and H.~Su, \emph{On a simple unital projectionless {$C^*$}-algebra},
  Amer. J. Math. \textbf{121} (1999), no.~2, 359--413.

\bibitem{Kir:ICM}
E.~Kirchberg, \emph{{Exact $C^*$-algebras, tensor products, and the
  classification of purely infinite $C^{*}$-algebras}}, Proceedings of the
  International Congress of Mathematicians, Z\"urich, 1994 (Basel), vol. 1,2,
  Birkh\"auser, 1995, pp.~943--954.

\bibitem{Kir:CentralSequences}
\bysame, \emph{{Central sequences in $C^*$-algebras and strongly purely
  infinite $C^*$-algebras}}, Abel Symposia \textbf{1} (2006), 175--231.

\bibitem{KirRor:pi2}
E.~Kirchberg and M.~R{\o}rdam, \emph{{Infinite non-simple $C^*$-algebras:
  absorbing the Cuntz algebra $\mathcal{O}_\infty$}}, Adv. Math. \textbf{167}
  (2002), no.~2, 195--264.

\bibitem{KirWinter:dr}
E.~Kirchberg and W.~Winter, \emph{{Covering dimension and quasidiagonality}},
  Internat. J. Math. \textbf{15} (2004), 63--85.

\bibitem{Lin:asu-class}
H.~Lin, \emph{{Asymptotic unitary equivalence and classification of simple
  amenable $C^*$-algebras}}, Invent. Math. \textbf{183} (2011), no.~2,
  385--450.

\bibitem{LinNiu:KKlifting}
H.~Lin and Z.~Niu, \emph{Lifting {$KK$}-elements, asymptotic unitary
  equivalence and classification of simple {$C\sp \ast$}-algebras}, Adv. Math.
  \textbf{219} (2008), no.~5, 1729--1769.

\bibitem{NgWinter:subhom}
P.~W. Ng and W.~Winter, \emph{{A note on subhomogeneous $C^{*}$-algebras}}, C.
  R. Acad. Sci. Canada \textbf{28} (2006), 91--96.

\bibitem{OrtPerRor:cfp-dr}
E.~Ortega, F.~Perera, and M.~R{\o}rdam, \emph{{The corona factorization
  property and stability of $C^*$-algebras with finite decomposition rank}},
  arXiv preprint math.OA/0903.2917v2, 2009.

\bibitem{Rfl:sr}
M.~Rieffel, \emph{{Dimension and stable rank in the $K$-theory of
  $C^*$-algebras}}, Proc. London Math. Soc. \textbf{46} (1983), no.~(3),
  301--333.

\bibitem{Rob:dimnuc-comparison}
L.~Robert, \emph{{Nuclear dimension and $n$-comparison}}, arXiv preprint
  math.OA/1002.2180v2, to appear in M{\"u}nster J. Math., 2010.

\bibitem{Ror:UHFII}
M.~R{\o}rdam, \emph{{On the structure of simple $C^*$-algebras tensored with a
  UHF-algebra, II}}, J. Funct. Anal. \textbf{107} (1992), 255--269.

\bibitem{Ror:encyc}
\bysame, \emph{Classification of nuclear, simple {$C\sp *$}-algebras},
  Classification of nuclear {$C\sp *$}-algebras. {E}ntropy in operator
  algebras, Encyclopaedia Math. Sci., vol. 126, Springer, Berlin, 2002,
  pp.~1--145.

\bibitem{Ror:Z-absorbing}
\bysame, \emph{The stable and the real rank of {$\mathcal{Z}$}-absorbing
  {$C^*$}-algebras}, Internat. J. Math. \textbf{15} (2004), no.~10, 1065--1084.

\bibitem{RorWin:Z-revisited}
M.~R{\o}rdam and W.~Winter, \emph{The {Jiang--Su} algebra revisited}, J. Reine
  Angew. Math. \textbf{642} (2010), 129--155.

\bibitem{Toms:ash-sdg}
A.~S. Toms, \emph{{$K$-theoretic rigidity and slow dimension growth}}, Invent.
  Math. \textbf{183} (2011), no.~2, 225--244.

\bibitem{TomsWin:ssa}
A.~S. Toms and W.~Winter, \emph{Strongly self-absorbing {$C^*$}-algebras},
  Trans. Amer. Math. Soc. \textbf{359} (2007), no.~8, 3999--4029.

\bibitem{TomsWin:ZASH}
\bysame, \emph{{$\mathcal{Z}$}-stable {ASH} algebras}, Canad. J. Math.
  \textbf{60} (2008), no.~3, 703--720.

\bibitem{TomsWinter:minhom}
\bysame, \emph{{Minimal dynamics and $K$-theoretic rigidity: Elliott's
  conjecture}}, arXiv preprint math.OA/0903.4133, 2009.

\bibitem{TomsWinter:PNAS}
\bysame, \emph{{Minimal dynamics and the classification of {$C^*$}-algebras}},
  Proc. Natl. Acad. Sci. USA \textbf{106} (2009), no.~40, 16942--16943.

\bibitem{TomsWinter:VI}
\bysame, \emph{{The {E}lliott conjecture for {V}illadsen algebras of the first
  type}}, J. Funct. Anal. \textbf{256} (2009), no.~5, 1311--1340.

\bibitem{Win:cpr}
W.~Winter, \emph{Covering dimension for nuclear {$C^*$}-algebras}, J. Funct.
  Anal. \textbf{199} (2003), no.~2, 535--556.

\bibitem{Winter:fintopdim}
\bysame, \emph{{On topologically finite dimensional simple $C^{*}$-algebras}},
  Math. Ann. \textbf{332} (2005), 843--878.

\bibitem{Win:localizingEC}
\bysame, \emph{{Localizing the Elliott conjecture at strongly self-absorbing
  $C^{*}$-algebras}}, arXiv preprint math.OA/0708.0283v3, with an appendix by
  H.\ Lin, 2007.

\bibitem{Winter:lfdr}
\bysame, \emph{{Simple $C^*$-algebras with locally finite decomposition rank}},
  J. Funct. Anal. \textbf{243} (2007), 394--425.

\bibitem{Winter:dr-Z-stable}
\bysame, \emph{{Decomposition rank and $\mathcal{Z}$-stability}}, Invent. Math.
  \textbf{179} (2010), no.~2, 229--301.

\bibitem{Win:ssa-Z-stable}
\bysame, \emph{{Strongly self-absorbing {$C^*$}-algebras are
  {$\mathcal{Z}$}-stable}}, J. Noncomm. Geom. \textbf{5} (2011), no.~2,
  253--264.

\bibitem{WinZac:order-zero}
W.~Winter and J.~Zacharias, \emph{Completely positive maps of order zero},
  M\"unster J. Math. \textbf{2} (2009), 311--324.

\bibitem{WinterZac:dimnuc}
\bysame, \emph{{The nuclear dimension of $C^{*}$-algebras}}, Adv. Math.
  \textbf{224} (2010), 461--498.

\end{thebibliography}

\providecommand{\bysame}{\leavevmode\hbox to3em{\hrulefill}\thinspace}
\providecommand{\MR}{\relax\ifhmode\unskip\space\fi MR }
\providecommand{\MRhref}[2]{%
  \href{http://www.ams.org/mathscinet-getitem?mr=#1}{#2}
}
\providecommand{\href}[2]{#2}

\end{document}